%% file: B-MS2versionseptembre2024.tex
\documentclass{amsart}

\usepackage{wasysym}
\usepackage{amsmath}
\usepackage{amssymb}
\usepackage{amsthm}
\usepackage{tikz}
\usetikzlibrary{matrix,arrows,backgrounds,shapes.misc,shapes.geometric,patterns,calc,positioning}
\usetikzlibrary{calc,shapes}
\usepackage{wrapfig}
\usepackage{epsfig}  
\usepackage[margin=1.3in]{geometry}
\usepackage{color}	 %
\usepackage[color,matrix,arrow]{xy}
\usepackage[shortlabels]{enumitem}
\usepackage{tikz-cd}
\xyoption{poly}
\xyoption{2cell}
\xyoption{all}
\newtheorem{thm}{Theorem}[section]
\newtheorem{prop}[thm]{Proposition}
\newtheorem{lemma}[thm]{Lemma}
\newtheorem{corollary}[thm]{Corollary}

\newtheorem*{thm*}{Theorem}

\usepackage[colorlinks=true, pdfstartview=FitV, linkcolor=purple, citecolor=blue,
filecolor=violet, urlcolor=violet]{hyperref}
\theoremstyle{definition}
\newtheorem{definition}[thm]{Definition}
\theoremstyle{remark}
\newtheorem{remark}[thm]{Remark}
\newtheorem{example}[thm]{Example}
\numberwithin{equation}{section}

\newcommand{\xx}{\mathbf{x}}
\newcommand{\yy}{\mathbf{y}}
\newcommand{\cals}{\mathcal{S}}

\newcommand{\calf}{\mathcal{F}}
\newcommand{\cala}{\mathcal{A}}

\newcommand{\calx}{\mathcal{X}}

\newcommand{\call}{\mathcal{L}}

\newcommand{\ssi}{\Leftrightarrow}
\newcommand{\za}{\alpha}
\newcommand{\zb}{\beta}

\newcommand{\zD}{\Delta}

\newcommand{\zg}{\gamma}

\newcommand{\zs}{\sigma}
\newcommand{\zS}{\Sigma}

\newcommand{\ot}{\leftarrow}

\DeclareMathOperator{\red}{red}
\DeclareMathOperator{\Hopf}{Hopf}

\newcommand{\Qbar}{\overline{Q}}

\newcommand{\Supp}{\textup{Supp\,}}

\newcommand{\PP}{\mathbb{P}}
\newcommand{\ZZ}{\mathbb{Z}}
\newcommand{\RR}{\mathbb{R}}

\newcommand{\QQ}{\mathbb{Q}}

\newcommand{\Axyq}{\mathcal{A}(\mathbf{x},\mathbf{y},Q)}
\newcommand{\xyq}{(\mathbf{x},\mathbf{y},Q)}
\newcommand{\dimprime}{\dim^\circ}

\begin{document}

\title[Knot theory and cluster algebras II]{Knot theory and cluster algebras II: The knot cluster} 
\author{V\'eronique Bazier-Matte}
\thanks{The first author was supported by the Discovery Grants  program from The Natural Sciences and Engineering Research Council of Canada and the Research Support for New Academics from the Fonds de Recherche du Qu\'ebec Nature et Technologie}
\address{D\'epartement de math\'ematiques et de statistique,  Universit\'e Laval, Qu\'ebec (Qu\'ebec), G1V 0A6, Canada}
\email{veronique.bazier-matte.1@ulaval.ca}
\author{Ralf Schiffler}
\thanks{The second author was supported by the National Science Foundation grants  DMS-2054561 and DMS-2348909}
\address{Department of Mathematics, University of Connecticut, Storrs, CT 06269-1009, USA}
\email{schiffler@math.uconn.edu}


\setcounter{tocdepth}{1}

\begin{abstract}
To every knot (or link) diagram $K$, we associate a cluster algebra $\cala$ that contains a cluster $\xx_t$ with the property that \emph{every} cluster variable in $\xx_t$ specializes to the Alexander polynomial $\zD_K$ of $K$.
We call $\xx_t $ the \emph{knot cluster} of $\cala$. Moreover, the specialization does not depend on the choice of the cluster variable. Furthermore, there exists a cluster automorphism of $\cala$ of order two that maps the initial cluster to the cluster $\xx_t$. 

We realize this connection between knot theory and cluster algebras in two ways. In our previous work, we constructed indecomposable representations $T(i)$ of the initial quiver $Q$ of the cluster algebra $\cala$. Modulo the removal of 2-cycles, the quiver $Q$ is the incidence quiver of the segments in $K$, and the representation $T(i)$ of $Q$ is built  by taking successive boundaries of $K$ cut open at the $i$-th segment. The relation to the Alexander polynomial stems from an isomorphism between the submodule lattice of $T(i)$ and the lattice of Kauffman states of $K$ relative to segment $i$.

In the current article, we identify the knot cluster $\xx_t$ in $\cala$ via a  sequence of mutations that we construct from a sequence of bigon reductions and generalized Reidemeister III moves on the diagram $K$. On the level of diagrams, this sequence first reduces $K$ to the Hopf link, then reflects the Hopf link to its mirror image, and finally rebuilds (the mirror image of) $K$ by reversing the reduction.
We show that every diagram of a prime link admits such a sequence.

We further prove that the cluster variables in $\xx_t$ have the same $F$-polynomials as the representations $T(i)$. This establishes the important fact that our representations $T(i)$ do indeed correspond to cluster variables in $\cala$. But it even establishes the much stronger result that these cluster variables are all compatible, in the sense that they form a cluster.

    We also prove that the representations $T(i)$ have the following symmetry property. For all vertices $i,j$ of $Q$, we have $\dim T(i)_j=\dim T(j)_i$.

    As applications, we show that in the Newton polytope of the $F$-polynomial of $T(i)$, every lattice point is a vertex. We also obtain an explicit formula for the denominator vector of each cluster variable in the knot cluster. In particular, each entry in the denominator vector is either 0 or 1.
\end{abstract}

\maketitle
\tableofcontents

\section{Introduction}
\emph{Cluster algebras} were introduced by Fomin and Zelevinsky in 2002 in \cite{FZ1} in order to provide an algebraic framework for canonical bases in Lie theory. Quickly thereafter, deep connections  developed between cluster algebras and a variety of research areas in mathematics and physics, including combinatorics, representation theory, hyperbolic geometry, algebraic geometry, dynamical systems, string theory, and scattering theory. In this article we establish a connection to knot theory.

A cluster algebra $\cala$ is a certain subring of a field or rational functions in $N$ variables $x_1,\ldots,x_N$. It is defined by specifying generators, the so-called \emph{cluster variables}, that are constructed recursively by a process called \emph{mutation} from an initial seed $(\xx,\yy,Q)$. Here $\xx=(x_1,\ldots,x_N)$ is called the initial \emph{cluster}, $\yy=(y_1,\ldots,y_N)$ the initial \emph{coefficient tuple} and $Q$ is a quiver (oriented graph) without oriented cycles of length 1 or 2. 

The mutation $\mu_k$ in direction $k=1,\ldots,N$ produces a new seed $(\xx',\yy',Q')$ whose quiver $Q'$ only depends on $Q$, $\yy'$ depends on $\yy $ and $Q$, and $\xx'$ depends on the whole seed $(\xx,\yy,Q)$. The new cluster $\xx'$ is obtained from $\xx$ be replacing the cluster variable $x_k$ by a new cluster variable $x_k'$ that is a Laurent polynomial in $\xx,\yy$. At every seed, one performs mutations in all $N$ directions, so one can visualize the mutation process as an $N$-regular tree rooted at the initial seed, whose vertices are the seeds and whose edges are the mutations. Each mutation produces one new cluster variable, unless it has been obtained already before. The cluster algebra is generated inside the field of rational functions by the set of all  cluster variables obtained by sequences of mutations. It is determined (up to coefficients) by the choice of the initial quiver $Q$. 

It has been shown that each cluster variable is a Laurent polynomial in $\xx,\yy$ \cite{FZ1} with positive (integer) coefficients \cite{LS4}.
 The \emph{$F$-polynomial} of a cluster variable is obtained by setting the initial cluster variables $x_1,\cdots,x_N$ equal to 1, and using the so-called principal coefficients for $\yy$. The $F$-polynomial is a polynomial in $\ZZ_{\ge 0}[y_1,\ldots,y_N]$ with constant term  equal to one \cite{FZ4}. 
One can also define the $F$-polynomial of a quiver representation $M$ as
$F_M=\sum_{\mathbf e} \chi(\textup{Gr}_{\mathbf e} M) \yy^{\mathbf{e}},$
where  $\mathbf{e}$ runs over all dimension vectors,  $\textup{Gr}_{\mathbf e} M $ is the quiver Grassmannian of all subrepresentations of $M$ of dimension vector $\mathbf e$, and $\chi$ denotes the Euler characteristic. If $M$ is reachable by mutation from the zero representation, then its $F$-polynomial $F_M$ is equal to the $F$-polynomial of the cluster variable that is obtained by the same mutation sequence \cite{DWZ2}.
 
 
 The number of cluster variables in $\cala$ is finite if and only if the quiver $Q$ is mutation equivalent to a simply laced Dynkin diagram \cite{FZ2}. More generally, if the quiver $Q$ is  mutation equivalent only to a finite number of quivers then the cluster algebra $\cala(Q)$ is said to be of \emph{finite mutation type}. This is the case if and only if the rank $N$ is equal to 2, or the quiver $Q$ is the adjacency quiver of a triangulation of a surface with marked points, or it is one of 11 exceptional types including the simply laced Dynkin and affine Dynkin quivers of type $\mathbb{E}$ \cite{FST}. 
 
Beyond finite mutation type,  the computation of the cluster variables along an arbitrary mutation sequence is extremely difficult. The cluster algebras in this article are \emph{not} of finite mutation type, except for a few small examples. Nevertheless, our main tool is a well-chosen mutation sequence along which we can control the computations.
 
 \smallskip
 In this work, we establish a relation to knot theory.
 A \emph{knot} is the homeomorphic image of a circle in $\RR^3$. 
 A \emph{link} with $r$ components is the homeomorphic image of $r$ disjoint circles in $\RR^3$. Thus a knot is a link with one component. A link is called \emph{prime} if it cannot be decomposed into the connected sum of two simpler, non-trivial links.. A \emph{link diagram} is a projection of a link into the plane together with the additional information at each crossing point which of the two strands is on top and which is below. 
 
 A \emph{link invariant} is a mathematical entity defined for every link that is invariant under ambient isotopy of the link. The study and creation of  invariants is an important branch of knot theory. One motivation for it is to be able to decide whether two links are isotopic. A good invariant is one that can be computed easily and that can distinguish between many different links. In addition, a good invariant may carry information about the topology of the link.

In this work we are concerned with  the Alexander polynomial, a knot invariant introduced by Alexander in 1928 \cite{Alexander}. It can be defined in several rather different ways. The original algebraic definition is as follows. The \emph{Alexander module} of a knot  is the first homology $H_1(X_{\infty})$ of the infinite cyclic cover of the knot complement $X$. Here the abelian group $H_1(X_{\infty})$ is considered as a module over the ring $\ZZ[t,t^{-1}]$ of Laurent polynomials in $t$, where $t$ acts on $X$ as the covering transformation given by the shift in the cyclic cover. 
By definition, the \emph{Alexander polynomial}  $\zD$ is a generator of the first elementary ideal of the Alexander module. Since multiplying a generator of an ideal by a unit of the ring will produce another generator, it follows that the Alexander polynomial is defined only up to multiplication by $\pm t^j$, with $j\in \ZZ$.

There are many equivalent definitions of the Alexander polynomial. Let us mention 
(a) Alexander's original construction as the determinant of a matrix that is constructed directly from a diagram of the knot \cite {Alexander},
(b) Conway's definition via skein relations \cite{Conw}, and
(c) Kauffman's realization as a weighted sum of Kauffman states \cite{K}.  In particular, these definitions also apply to links.

The Alexander polynomial is the oldest polynomial  knot invariant and is very well-studied. A generalization was given by the HOMFLY-PT polynomial, a 2-variable skein polynomial that specializes to the Alexander polynomial (as well as to the Jones polynomial) \cite{homfly, pt}.
More recently, the Alexander polynomial has been generalized in the work of Osv\'ath and Szab\'o \cite{OS}, as well as Rasmussen \cite{Ras}, on knot Floer homology.

\subsection{Statement of the main result}
Let $K$ be a diagram of a prime link with $n$ crossing points. We may assume without loss of generality that $K$ does not contain a monogon, since one can always remove it via a Reidemeister 1 move. 
Let $Q'$ be the incidence quiver of the segments in $K$. More precisely, $Q'$ is a planar quiver whose vertices are in bijection with the segments in $K$, and the arrows of $Q'$ 
are in bijection with pairs $(p,R)$, with $p$ a crossing point and $R$ a region incident to $p$.
 An example is given in Figure \ref{example 2112}. In particular,
$Q'$ has $2n$ vertices and 
\begin{itemize}
\item [-] every vertex of $Q'$  has exactly two incoming and two outgoing arrows;
\item[-] every arrow of $Q'$ lies in exactly two chordless cycles: a clockwise cycle of length 4 determined by the crossing point of the pair $(p,R)$, and a counterclockwise cycle determined by the region $R$ of the pair.
\end{itemize}

Note that every bigon in the diagram $K$ gives rise to an oriented 2-cycle in $Q'$. In order to define a cluster algebra, we let $Q$ be the quiver obtained from $Q'$ by removing the 2-cycles.

\begin{definition} Let $K$ be a diagram of a prime link and let $Q$ be its quiver as defined above. The \emph{cluster algebra $\cala(K)$ of $K$} 
   is the cluster algebra $\cala(\xx,\yy,Q)$ with principal coefficients and initial quiver $Q$.
\end{definition}

We need one more definition.  From now on we fix an orientation on each component of the link diagram $K$.
\begin{definition}
 Let $\xx=(x_1,\ldots,x_{2n}),\ \yy=(y_1,\ldots,y_{2n})$ be the cluster and coefficient tuple of the initial seed in $\cala(K)$. The \emph{Alexander specialization} is the  substitution $x_j=1$ and

\begin{equation}
\label{specializationintro}
y_j=\left\{
\begin{array}
 {ll}
 -t &\textup{if segment $j$ runs from an undercrossing to an overcrossing;}\\  
 -t^{-1} &\textup{if segment $j$ runs from an overcrossing to an undercrossing;}\\
 -1 &\textup{if segment $j$ connects two overcrossings or two undercrossings,} \\\end{array}
 \right.
\end{equation}
for all $j=1,\ldots, 2n$.
\end{definition}
\medskip
We are now ready to state our main result. Let $Q^{op}$ denote the opposite quiver of $Q$ obtained by reversing the direction of every arrow.
\begin{thm}
 \label{thm main} Let $K$ be the diagram of a prime link and let $Q$ be its quiver.
 Then  the cluster algebra $\cala(K)$ contains a seed $\zS_t=(\xx_t,\yy_t, Q_t)$ that satisfies the following. 
 
 \begin{itemize}
\item [(a)] Every cluster variable $x_{i;t}$ in this seed specializes to the Alexander polynomial of $K$ under 
 the Alexander specialization (\ref{specializationintro}).

\item [(b)] There is a  permutation $\zs\in S_{2n}$ that induces an isomorphism of quivers $\zs\colon Q \to Q_t^{op}$.
 
\item [(c)] The permutation $\zs$ induces a cluster automorphism $\phi_{\sigma}$ of the cluster algebra with trivial coefficients defined on the initial cluster variables by $\phi_{\sigma}(x_{i;0})=x_{\zs(i);t}$. Moreover $\zs^2=1$ and $\phi_{\sigma}^2 = 1$.

\end{itemize}
\end{thm}

This is a surprising result for a variety of reasons. First, the quiver $Q$ is not of finite mutation type, which means that computing the cluster variables along an arbitrary mutation sequence would be extremely difficult. So in this theorem we are dealing with a special situation within this cluster algebra.

Second, finding \emph{one} cluster variable  that specializes to the Alexander polynomial may not be too astonishing. After all there are a lot of cluster variables and they have 4n variables that we can specialize and the Alexander polynomial only has one. However, the fact that there are $2n$ cluster variables, that all specialize to the Alexander polynomial  \emph{under the same specialization} seems quite striking. 
In addition, these $2n$ cluster variables together form a single cluster in the cluster algebra!

\begin{definition}
 The cluster $\xx_t$ in Theorem \ref{thm main} is called the \emph{knot cluster} of $\cala(K)$. 
\end{definition}
\subsection{Interpretation as representations of $Q$}
Let $K$ and $Q$ be as in the previous subsection. Define $W$ to be the potential given by the sum of all clockwise chordless cycles in $Q$ minus the sum of all counterclockwise chordless cycles. Let $B$ denote the Jacobian algebra of the quiver with potential $(Q,W)$. 

In our previous work \cite{BMS}, we defined an indecomposable $B$-module $T(i)$ for every segment $i$ of $K$ and showed that its $F$-polynomial specializes to the Alexander polynomial of $K$ via the Alexander specialization (\ref{specializationintro}). The module $T(i)$ is defined as a representation of the quiver $Q$, thus it consists of a vector space $T(i)_j$ for every vertex $j$ of $Q$ and a linear map $T(i)_\za\colon T(i)_a\to T(i)_b$ for every arrow $\za\colon a \to b$ in $Q$. 

The vector spaces of $T(i) $ are defined by taking successive boundaries of the diagram $K$ cut open at the segment $i$.  This process may seem similar to peeling an onion in the sense that one removes one layer at a time. The dimension of the vector space $T(i)_j$ is equal to $d$ if and only if $j$ lies in the $d$-th boundary.
Furthermore, if $\za\colon a\to b$  is an arrow in $Q$, the linear map $T(i)_\za$ is given by one of the following matrices
\begin{itemize}
\item [-] the identity matrix or the Jordan block with eigenvalue zero, if $\dim T(i)_a=\dim T(i)_b$;
\item [-] the identity matrix with an extra row of zeros at the bottom, if $\dim T(i)_a=\dim T(i)_b-1$;
\item [-] the identity matrix with an extra column of zeros on  the left, if $\dim T(i)_a=\dim T(i)_b+1$.
\end{itemize}
Moreover, the definition is such that, for every chordless cycle $w=\za_1\cdots\za_l$ in $Q$, the composition $T(i)_{\za_l}\circ\ldots\circ T(i)_{\za_1}$ is the Jordan block with eigenvalue zero. We refer to \cite{BMS} for precise definitions.

\smallskip

The main step in the proof of Theorem~\ref{thm main} is to show that the cluster $\xx_t$ corresponds to the  direct sum $T=\oplus_{i\in Q_0} T(i)$. More precisely, we show the following. 

\begin{thm}\label{thm module main}
Let $\zs$ be the permutation of Theorem \ref{thm main}. Then for all $i=1,\ldots, 2n$, the cluster variable $x_{i;t}$ has the same $F$-polynomial as the representation $T(\zs(i))$. 
\end{thm}
 
To prove this result we exhibit an explicit mutation sequence that transforms the initial seed of $\cala(K)$ into the seed $\zS_t$ of Theorem~\ref{thm main}. This sequence of mutations is parallel to a sequence of bigon reductions and diagrammatic Reidemeister 3 moves, or RD3 moves for short, on the diagram $K$. 
An RD3 move is similar to a Reidemeister 3 move in the sense that it moves a segment $i$ across the crossing point of two other segments $j$ and $k$. However, the RD3 move does not take into account the information of over and undercrossings of the three strands $i,j,k$ involved. In particular, an RD3 move may not be an ambient isotopy and therefore it may change the underlying link of the diagram.

 This sequence first reduces $K$ to a Hopf link, then reflects the Hopf link to its mirror image, and then rebuilds the mirror image of $K$ by reversing the reduction sequence. 
The reduction of a bigon formed by the segments $i$ and $j$ corresponds to the two mutations in direction $i$ and $j$. An RD3 move that moves a segment $i$ through a crossing formed by segments $j$ and $k$ corresponds to the mutation sequence at $i,j,k,i,j,k,j,k,j$. 
 
Not every link diagram contains a bigon. The smallest example are the Borromean rings, see section~\ref{sect Borromean}. This is the reason why we have to include the RD3 moves as a means to produce bigons. 
We say that two triangular regions of the diagram $K$ are \emph{disjoint} if they do not have a common vertex. A sequence of RD3 moves is \emph{admissible} if for every pair of consecutive moves their two triangular regions are   disjoint.

With this notation, we prove the following result. 

\begin{thm} \label{thm::sequenceR3 intro}
  Every link diagram contains a bigon or it admits an admissible sequence of RD3 moves creating a bigon.
\end{thm}

\subsection{Properties of the representation $T(i)$} 
The representations $T(i)$ are interesting in their own right. 
In our previous work \cite{BMS}, we proved that $T(i)$ is an indecomposable $B$-module and has the following properties.
\begin{enumerate}
\item For every dimension vector $\mathbf e$, the quiver Grassmannian $\textup{Gr}_{\mathbf e} (T(i))$ is either a single point or the empty set \cite[Corollary 6.8]{BMS}.
In other words, the subrepresentations of $T(i)$ are uniquely determined by their dimension vector. A direct consequence is that every coefficient in the $F$-polynomial of $T(i)$ is equal to 1.
\item The lattice of submodules of $T(i)$ is isomorphic to the lattice of Kauffman states of the diagram $K$ relative to the segment $i$ \cite[Theorem 6.9]{BMS}. 
\end{enumerate}

In the current article, we also show the following result, since it is needed in the proof of the main theorem. 

\begin{thm}\label{thm dim sym intro} (Symmetry of dimension) 
Let $i,j$ be two segments of $K$ and $T(i),T(j)$ the corresponding representations of $Q$. Then
 \[\begin{array}{rcl}
\dim T(i)_j &=&\dim T(j)_i.
\end{array}\]
 
\end{thm}
This result  seems natural from the topological point of view, since it means that the segment $j$ lies in the $d$-th boundary of $K$ opened at the segment $i$ if and only if the segment $i$ lies in the $d$-th boundary of $K$ opened at the segment $j$.  
However, from a representation theoretic perspective, this property is very unusual. A similar equation in representation theory is 
$\dim P(i)_j=\dim I(j)_i$ where $P(i)$ and $I(i)$ are the indecomposable projective and injective $B$-modules, respectively. However, the $P(i)$ and $I(i)$ are \emph{two} indecomposable modules for each vertex in $Q$, whereas the equation in Theorem~\ref{thm dim sym intro} only involves \emph{one}. 

An immediate consequence is the following statement. 
\begin{corollary}
 \label{cor intro}
 Let $T=\oplus_{i\in Q_0} T(i)$ be the direct sum of all $T(i)$. Then for every vertex $j$, we have
 \[\begin{array}{rcl}
\dim T_j &=&\dim T(j).
\end{array}\]

\end{corollary}
 
\begin{proof} We have
 $\dim T_j=\sum_{i\in Q_0}\dim T(i)_j 
 =\sum_{i\in Q_0}\dim T(j)_i=\dim T(j)$, 
 where the second identity holds by Theorem~\ref{thm dim sym intro}.
\end{proof}

 \subsection{Applications}\label{sect intro appl}
 We also obtain several  applications.  
\subsubsection{Denominator vectors}
The \emph{denominator vector} of a cluster variable $u$ is the multidegree vector of the denominator of the Laurent expansion of $u$ in the initial seed.
 
\begin{thm}
 \label{thm denominator} Let $x_{i;t}$ be any cluster variable in the knot cluster $\xx_t$. Then all entries of the denominator vector of $x_{i;t}$ are  0 or 1. 

 More precisely, the denominator  is equal to 
 $\prod_{j} x_{j} $, where $j$ runs over all segments of $K$ that do not bound a common region with $\zs (i)$. 
\end{thm}

An equivalent formulation of the last statement is that the denominator of the cluster variable corresponding to the representation $T(i)$ is equal to $\prod_{j\in\Supp T(i)} x_j$, where the product runs over all vertices $j$ of $Q$ such that $T(i)_j$ is nonzero. 

This is a surprising result, since the dimension of the representation $T(i)$ can be arbitrary large.

\subsubsection{Newton polytopes} 
 This application uses results from \cite{LP}. The \emph{Newton polytope} of a multivariate polynomial $\sum c_{(a_1,\ldots,a_N)}y_1^{a_1}\cdots y_N^{a_N}$ is the convex hull of the of the lattice points ${(a_1,\ldots,a_N)}$ with $c_{(a_1,\ldots,a_N)}\ne 0$.
 A lattice point $v$ of the Newton polytope is called a \emph{vertex} if $v$ is contained in every set of lattice points whose convex hull is the whole polytope. 
 Thus $v$ does not lie on the straight line segment between any two points in the polytope.
 \begin{thm}
 \label{thm NewtonPolytop}
 In the Newton polytope of the $F$-polynomial of $T(i)$ every lattice point is a vertex. 
\end{thm}

\subsubsection{Reidemeister moves}
A natural question is how the cluster algebra $\cala(K)$ behaves under changes to the link diagram $K$. In particular, can we interpret  the effect of a Reidemeister move in the cluster algebra? 
For the first two Reidemeister moves, the cluster algebra definitely changes because the moves change the number of crossing points in $K$ and hence the number of vertices in $Q$. 
For the third Reidemeister move we have the following result. 

\begin{thm}
 \label{thm Reidemeister3}
Let $K$ and $K'$ be two link diagrams without curls that are obtained from each other by the third Reidemeister  move on segments $i,j,k$. Then their quivers are obtained from each other by the mutation sequence on vertices $i,j,k,i$.
\end{thm}

\subsubsection{A new proof of a theorem by Murasugi}  Using our representation theoretic approach, we obtain a new proof of the following  classical result of \cite{Murasugi}.
\begin{thm}
 \label{thm Murasugi:intro}
\cite{Murasugi} 
The Alexander polynomial of an alternating knot is an alternating sum \[\sum_{j=m}^M a_j t^j,\] with $a_j\ne 0$ for all $j=m,m+1,\ldots,M$.
\end{thm}

\begin{remark}
    The difference $M-m$ is equal to twice the genus of the knot, see \cite{Crowell, Murasugi genus}.
\end{remark}

\subsubsection{A subsequent specialization}
Using the classical result that the evaluation of the Alexander polynomial at $t=-1$ is equal to $\pm 1$, if $K$ is a knot; and $0$, otherwise, we obtain the following result. 
\begin{corollary}
 \label{cor intro 2}
 For every cluster variable $u$ in the knot cluster of $\cala(K)$, the specialization at $x_j=1, y_j=-1$, $j=1,\ldots, 2N$, is equal to  $\pm 1$ if $K$ is a knot; and equal to   0 if $K$ is a link with at least 2 components.
%
\end{corollary}
\begin{proof}
 We know that the Alexander polynomial $\zD(t)$ is equal to the Alexander specialization of $u$. 
 Setting $t=-1$ in the Alexander specialization is the same as specializing each $y_j$ to $-1$. 
\end{proof}
 \subsection{A few words about the proofs}
 The proof of Theorem \ref{thm main} is an inductive argument along the mutation sequence. It involves rather intricate computations in the cluster algebra that uses a lot of the cluster algebra machinery such as $F$-polynomials, separation formula, sign coherence and the duality between the $g$-matrix and the $c$-matrix. Moreover it involves a detailed analysis of the local poset structure of the Kauffman states, which turns out to have an interesting staircase structure in the presence of an RD3 move, see Figure~\ref{fig staircase}. 
 The proof of the symmetry of dimension of $T(i)$ is a topological argument, and the proof of the creation of bigons is via a combinatorial algorithm on the link diagrams.

 \subsection{Relation to other work} 
 \subsubsection{2-bridge links}
A special kind of links is the family of 2-bridge links. In \cite{LS6}, Kyungyong Lee and the second author realized  the Jones polynomial of  2-bridge links as specializations of cluster variables of Dynkin type $\mathbb{A}$. This is an ad-hoc construction that uses the fact that both the 2-bridge links as well as the cluster variables are determined by a continued fraction. This result can be reinterpreted in the setting of our work as follows, see also \cite[Section 8]{BMS}.
Each 2-bridge link $K$ is determined by an ordered sequence of braids on two strands that are connected to each other in a precise way. In particular, $K$ contains a special segment $i$ that shares a region with each one of these braids. The cluster variable used in \cite{LS6} to obtain the Jones polynomial is the same as the  variable $x_{\zs (i);t}$ at position $\zs (i)$ in our knot cluster $\xx_t$. In particular, it is obtained by our mutation sequence, which in this case does not need any RD3 moves. So, in this case, the knot cluster contains one variable that specializes to the Alexander polynomial and to the Jones polynomial (under different specializations). Unfortunately, the Jones specialization only works for this one particular segment $i$.

Morier-Genoud and Ovsienko gave another interpretation of the result of \cite{LS6} that realizes the Jones polynomial of the 2-bridge links using $q$-rationals and $q$-continued fractions \cite{MO}. 

Nagai and Terashima used ancestral triangles constructed from continued fractions to give a formula for the cluster variables of type $\mathbb{A}$ and then defined a specialization that produces the Alexander polynomial of the corresponding 2-bridge link, see \cite{NT}. Our specialization is a generalization of theirs.

 In \cite{CDR}, Cohen, Dasbach and Russel gave a realization of the Alexander polynomial for arbitrary knots as a sum over perfect matchings of the bipartite graph whose vertices are given by the crossing points and the regions of the diagram.  Their graph can be recovered from our quiver  by the methods used for plabic graphs, see for example \cite{FWZ}. In the case of 2-bridge knots, their graph is equivalent to the snake graph associated to the continued fraction in \cite{CS4} and in that case their formula seems to be a special case of the cluster variable expansion formula of \cite{MS} and therefore may be related to ours as well. However, in their approach, the weight of a perfect matching is given by edge weights, which in the cluster algebra setup corresponds to $x$-variables, whereas we use the $y$-variables instead. For arbitrary knots, it is unclear if their formula corresponds to a cluster variable.

All of the articles above  consider a single segment of the link to produce a formula for the invariant. In our approach, we rather aim at a conceptual understanding of the collection of the $2n$ objects given by all of the segments of the link inside the cluster algebra and in the module category of the Jacobian algebra. A step in this direction was obtained by David Whiting and the second author in \cite{SW}, who studied the very special case where $K$ is a 2-bridge link with only two braids. They showed that the module $T=\oplus T(i)$ 
given as the direct sum of all the $T(i)$ is a tilting module over the Jacobian algebra in this case. Their work did not consider connections to knot theory.

\subsubsection{Connection to braid varieties} 
We would like to mention here the very interesting work on braid varieties which  lies at the interface of cluster algebras and Legendrian links that  has been developed over the last years by many people, see 
\cite{Casals-Gao,CGGS,CGGLSS,CW,GL,ShWe}
and the references therein. 
This body of work uses symplectic geometry to associate a cluster algebra to a special type of knot given as the closure of a positive braid.  The vertices of the initial quiver are associated to certain regions in the knot diagram, so these cluster algebras are clearly different  from ours.
It would be very interesting to find a direct connection between their cluster algebras and ours.

 \subsection{Organization of the article} 
 After a brief preliminary section, we illustrate the proof of Theorem \ref{thm main} in two examples in section~\ref{sect examples}. The actual proof of the theorem, as well as the proof of Theorem~\ref{thm module main} follow in section \ref{sect proof}. 
 We prove Theorem~\ref{thm::sequenceR3 intro} in section~\ref{sect bigons} and Theorem~\ref{thm dim sym intro} in section~\ref{sect sym dim}.  The proofs of Theorem \ref{thm denominator},  \ref{thm NewtonPolytop}, \ref{thm Reidemeister3} and \ref{thm Murasugi:intro} are given in section~\ref{sect appl}.
 We prove Theorem~\ref{thm::sequenceR3 intro} in section~\ref{sect bigons} and Theorem~\ref{thm dim sym intro} in section~\ref{sect sym dim}.  The proofs of Theorem \ref{thm denominator},  \ref{thm NewtonPolytop} and \ref{thm Reidemeister3} are given in section~\ref{sect appl}.

\begin{example}\label{example 2112}
Consider the knot illustrated at Figure \ref{fig 2112}.
     \begin{figure}
\begin{center}
\scalebox{1}{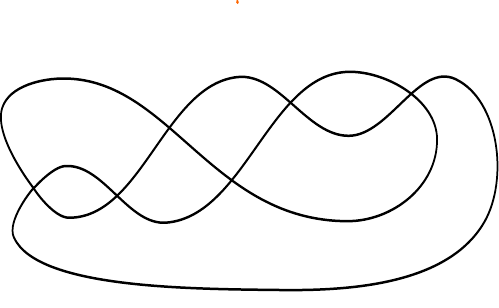 }
\caption{The knot [2,1,1,2].}
\label{fig 2112}
\end{center}
\end{figure}
Its Alexander polynomial is equal to $\zD=t^{-2}-3t^{-1}+5-3t+t^2$.
The quiver associated to this knot diagram is shown below: on the left with 2-cycles, on the right without 2-cycles.
\[ \begin{tikzpicture}
    \foreach \i/\x/\y in {1/0/3, 10/2.5/3, 4/4/3, 7/0/2, 11/1/2, 2/2/2, 5/3/2, 9/4/2, 12/0/1, 6/1.5/1, 3/4/1, 8/2/0} {
        \node (\i) at (\x,\y) {\i};
    }
    \draw[->] (1) -- (10);
    \draw[->] (10) -- (2);
    \draw[->] (2) -- (11);
    \draw[->] (11) -- (1);
    \draw[->] (10) -- (4);
    \draw[->] ([xshift=-2 pt]4.south) -- ([xshift=- 2 pt]9.north); 
    \draw[->] ([xshift=2 pt]9.north) -- ([xshift=2 pt]4.south);
    \draw[->] (9) -- (5);
    \draw[->] (5) -- (10);
    \draw[->] (7) -- (11);
    \draw[->] (11) -- (6);
    \draw[->] (6) -- (12);
    \draw[->] ([xshift=-2 pt]7.south) -- ([xshift=- 2 pt]12.north); 
    \draw[->] ([xshift=2 pt]12.north) -- ([xshift=2 pt]7.south); 
    \draw[->] (12) -- (8);
    \draw[->, looseness=1.8, out=180, in=-135] (8) to (1);
    \draw[->] (1) -- (7);
    \draw[->] (2) -- (5);
    \draw[->] (5) -- (3);
    \draw[->] (3) -- (6);
    \draw[->] (6) -- (2);
    \draw[->, looseness=1.8, out=-45, in=0] (4) to (8);
    \draw[->] (8) -- (3);
    \draw[->] (3) -- (9);    
\end{tikzpicture}  \begin{tikzpicture}
    \foreach \i/\x/\y in {1/0/3, 10/2.5/3, 4/4/3, 7/0/2, 11/1/2, 2/2/2, 5/3/2, 9/4/2, 12/0/1, 6/1.5/1, 3/4/1, 8/2/0} {
        \node (\i) at (\x,\y) {\i};
    }
    \draw[->] (1) -- (10);
    \draw[->] (10) -- (2);
    \draw[->] (2) -- (11);
    \draw[->] (11) -- (1);
    \draw[->] (10) -- (4);
    \draw[->] (9) -- (5);
    \draw[->] (5) -- (10);
    \draw[->] (7) -- (11);
    \draw[->] (11) -- (6);
    \draw[->] (6) -- (12);
    \draw[->] (12) -- (8);
    \draw[->, looseness=1.8, out=180, in=-135] (8) to (1);
    \draw[->] (1) -- (7);
    \draw[->] (2) -- (5);
    \draw[->] (5) -- (3);
    \draw[->] (3) -- (6);
    \draw[->] (6) -- (2);
    \draw[->, looseness=1.8, out=-45, in=0] (4) to (8);
    \draw[->] (8) -- (3);
    \draw[->] (3) -- (9);    
\end{tikzpicture}\]
The representations $T(i)$ are given below by their composition series. 
\[\begin{array}
    {cccccc}
    T(1)=\begin{smallmatrix}
        &&3\\
        &6&& 9\\
        12&\,2\\
        && 5
    \end{smallmatrix}
&
T(2)=\begin{smallmatrix}
    12&&4\\&8\\1&&3\\7&&9
\end{smallmatrix}
&
T(3)=\begin{smallmatrix}
&&1\\&10&&7\\&2\\4&&11
\end{smallmatrix}
&
T(4)=\begin{smallmatrix}
7&&2\\&11&&5\\&&&3\\&&6\\&&12

\end{smallmatrix}
&
T(5)=\begin{smallmatrix}
   7\\11\\6\\12&&4\\&8\\&1\\&7
\end{smallmatrix}
\\[30pt] 
T(6)=\begin{smallmatrix}
  9\\5&&1\\&10&&7\\&4
\end{smallmatrix}

&
T(7)=\begin{smallmatrix}
 &9 \\&5\\&10\\&4\\&8\\&3\\6&&9\\2\\&5
\end{smallmatrix}
&
T(8)=\begin{smallmatrix}
  7&&2&&9\\&11&&5
\end{smallmatrix}
&
T(9)=\begin{smallmatrix}
&7\\&11\\&6\\&12\\&8\\&1\\7&&10\\&&2\\&11
\end{smallmatrix}
&
T(10)=\begin{smallmatrix}
  &12\\&6&&9\\11&&3\\7
\end{smallmatrix}
\\ [30pt]
T(11)=\begin{smallmatrix}
  9\\5\\10\\4&&12\\&8\\&3\\&9
\end{smallmatrix}
&
T(12)=\begin{smallmatrix}
  &2&&9\\11&&5\\1\\&10\\&4
\end{smallmatrix}

\end{array}\]

The $F$-polynomials of these representations are quite different from each other. For example
\[\begin{array}{rl}
F_{T(1)}=& 1+y_5+y_{12}
+y_5y_{12}+y_2y_5+y_5y_9
+y_2y_5y_{12}+y_5y_9y_{12}
+y_2y_5y_9+y_2y_5y_9y_{12}
\\&+y_2y_5y_6y_{12}
+y_2y_5y_6y_9y_{12}
+y_2y_3y_5y_6y_9y_{12}
\\
F_{T(7)}=&
1+y_5+y_2y_5+y_5y_9+y_2y_5y_6+y_2y_5y_9+y_2y_5y_6y_9+y_2y_3y_5y_6y_9
+y_2y_3y_5y_6y_8y_9\\&
+y_2y_3y_4y_5y_6y_8y_9
+y_2y_3y_4y_5y_6y_8y_9y_{10}
+y_2y_3y_4y_5^2y_6y_8y_9y_{10}
+y_2y_3y_4y_5^2y_6y_8y_9^2y_{10}
\\
F_{T(8)}=&1+y_5+y_{11}
+y_5y_{11}+y_5y_9+y_7y_{11}
+y_5y_7y_{11}+y_5y_9y_{11}+y_2y_5y_9
+y_2y_5y_7y_{11}\\&+y_5y_7y_9y_{11}+y_2y_5y_9y_{11}
+y_2y_5y_7y_9y_{11}
\end{array}\]
However, each $F$-polynomial specializes to the Alexander polynomial  (up to a unit in $\mathbb{Z}[t^{\pm 1}]$) under the Alexander specialization (\ref{specializationintro}).

The knot cluster $\xx_t$ is obtained from the initial seed by the mutation sequence at the vertices
$7,12,4,9,1,11,3,5,6,10,2,8,3,5,1,11,4,9,7,12$ in that order. The permutation is given by \[\zs= (7\ 12)(2\ 9)(1\ 11)(3\ 5)(6\ 10)(2\ 8),\]
and the corresponding cluster automorphism $\phi_{\sigma}$ is given by $\phi_{\sigma}(x_{i,0})=x_{\sigma(i);t}$.
\end{example}


\section{Preliminaries}\label{sect preliminaries}
\subsection{Knot theory}\label{sect knots}

 A \emph{knot} is a subset of $\mathbb{R}^3$ that is homeomorphic to a circle.
A \emph{link with $r$ components} is a subset  of $\mathbb{R}^3$  that is homeomorphic to $r$ disjoint circles.   Thus a knot is a link with one component. Links are considered up to ambient isotopy. 
A  link is said to be \emph{prime} if it is not the connected sum of two nontrivial links. 

A \emph{link diagram} $K$ is a projection of the link into the plane, that is injective  except for a finite number of  double points that are called \emph{crossing points}. In addition, the diagram carries the information at each crossing point which of the two strands is on top and which is below. 
A diagram is called \emph{alternating} if traveling along a strand alternates between overcrossings and undercrossings.  A link is called \emph{alternating} if it has an alternating diagram.
A link is said to be  \emph{oriented} if for each component a direction of traveling along the strand is fixed. 

\subsubsection{Reidemeister moves} 
\begin{figure}[ht!] 
\begin{center}
\scalebox{1.2}{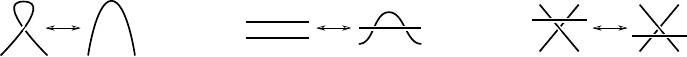 }
\caption{The three Reidemeister moves.}
\label{fig Reidemeister}
\end{center}
\end{figure}

The three Reidemeister moves illustrated in Figure~\ref{fig Reidemeister} are local moves on a link diagram.  It was shown in \cite{Reidemeister, AlexanderBriggs} that two link diagrams represent the same link
 if and only if they are related by a finite number of Reidemeister moves. 

\subsubsection{Curls, bigons and triangular regions}  A \emph{curl} is a monogon in the diagram. We usually assume without loss of generality that our link diagrams are without curls, because one can always remove them (by a Reidemeister I move) without changing the link.

The following result is well-known. We include a proof for convenience.
\begin{prop}
 \label{prop 8 triangles}
 Let $K$ be a link diagram without curls and bigons. Then $K$ contains at least 8 triangular regions.
\end{prop}
\begin{proof} 
Let $n$ be the number of crossing points in $K$ and let $a_i$ denote the number of regions in $K$ that have exactly $i$ sides. Thus $a_1=a_2=0$ by our assumption. Since the total number of regions  is $n+2$, we have
$ 4n+8=  4\sum_{i\ge 3} a_i.$
On the other hand, since every crossing point is incident to exactly four segments, we have
$4n=\sum_{i\ge 3} ia_i$ and therefore
\[8=4\sum_{i\ge 3} a_i-\sum_{i\ge 3} ia_i
=a_3-a_5-2a_6-3a_7\dots\]
Hence $a_3 $ is at least 8.
\end{proof}
\subsubsection{The Alexander polynomial}\label{sect Alexander}
The Alexander polynomial $\zD$ of an oriented link is a polynomial invariant of the link $\zD\in\ZZ[t^{\pm \frac{1}{2}}]$ that can be defined in homological terms, see \cite[Chapter 7]{Rolfsen} or \cite[Chapter 6]{Lick}. For the original definition of Alexander, see \cite{Alexander}.
The Alexander polynomial is defined up to multiplication by a signed power of $t$. 
  In \cite{Conw}, Conway showed that the Alexander polynomial $\zD$ of an oriented link $K$ can be defined recursively using skein relations.

In \cite{K}, Kauffman gave a description of the Alexander polynomial as a state sum relative to a segment $i$.  
A Kauffman state relative to  segment $i$ is a set of pairs $(x,R)$ where  $R$ is a region not incident to the segment $i$, $x$ is a crossing point incident to $R$ such that
\begin{itemize}
    \item  [(i)] each crossing point appears exactly once;
    \item  [(ii)] each region, except for the two regions incident to segment $i$, appears exactly once.
\end{itemize}
This approach is crucial for us and we reviewed it in detail in \cite[Section 3]{BMS}.

\smallskip

The Alexander polynomial $\zD$ has the following properties, see for example \cite[Chapter 6]{Lick}.
\begin{itemize}
\item [(i)] For any link, $\zD(t) \doteq \zD(t^{-1})$, where the symbol $\doteq$ means ``equal up to a signed power of $t$''. 
\item [(ii)]  $\zD(1) = \pm 1$ for any knot, and $\zD(1)=0$ for any link with at least 2 components.
\item [(iii)] For any knot 
\[\zD\doteq a_0+a_1(t^{-1}+t)+a_2(t^{-2}+t^2)+\dots \]
with $a_0$ odd.
\item[(iv)] If a knot has genus $g$ then $2g\ge \textup{breadth}(\zD)$, where the breadth is the difference between the maximal and the minimal degree of the polynomial.
\end{itemize}

\subsection{Cluster algebras} \label{sect cluster algebras}
We review basic notions and results on cluster algebras. The cluster algebras in this paper will always be skew-symmetric and defined over a tropical semifield. 
\subsubsection{Definition} We start by recalling the definition of skew-symmetric cluster algebras. We follow the exposition of \cite{Sln}.

Let $(\PP,\cdot)$ be a free abelian group (written
multiplicatively) on variables $y_1,\ldots,y_N$ 
and define an addition $\oplus$ in $\PP$ by
\begin{equation}
\label{eq:tropical-addition}
\prod_j y_j^{a_j} \oplus \prod_j y_j^{b_j} =
\prod_j y_j^{\min (a_j, b_j)} .
\end{equation}
Then  $(\PP,\oplus,\cdot)$ is a semifield
and is  called \emph{tropical
 semifield}. 
 Let $\ZZ\PP$ denote the group ring of $\PP$. 
Then  $\ZZ\PP$ is the ring of Laurent
polynomials in the variables~$y_1,\ldots, y_N$.
Let $\QQ\PP$ denote the field of fractions of $\ZZ\PP$ and 
let $\calf=\QQ\PP(x_1,\ldots,x_N)$ be the field of rational functions in $n$ variables and coefficients in $\QQ\PP$.

The cluster algebra is determined by the choice of an initial seed $\zS_0=\xyq$, which consists of the following data.
\begin{itemize}
\item $Q$ is a quiver  without loops 
and 2-cycles, 
and with $N$ vertices;
\item  $\yy=(y_1,\ldots,y_N)$ is the $N$-tuple of generators of $\PP$, called {\em initial coefficient tuple};
\item $\xx=(x_1,\ldots,x_N)$  is the $N$-tuple of variables of $\calf$, called {\em initial cluster}.
\end{itemize}


The \emph{seed mutation} $\mu_k$ in direction~$k$ transforms
$(\xx, \yy, Q)$ into the seed
$\mu_k(\xx, \yy, Q)=(\xx', \yy', Q')$ defined as follows:

\begin{itemize}
\item $\mathbf{x}'$ is obtained from $\mathbf{x}$ by replacing one cluster variable, $\mathbf{x}'=\mathbf{x}\setminus\{x_k\}\cup \{x_k'\}$, and $x_k'$ is defined by the following
\emph{exchange relation}

\begin{equation}\label{exch}
{\,x_k} x'_k = {\frac{1}{y_k\oplus 1}\left( y_k\prod_{i\to k} x_i \ 
+ \ \prod_{i\ot k} x_i \right)}
\end{equation}
 where 
 the first product runs over all arrows in $Q$ that end in $k$ and the second product runs over all arrows that start in $k$.
\item $\yy'=(y_1',\dots,y_N')$ is 
a new coefficient $N$-uple, where 
 \begin{equation*}
\label{eq:y-mutation}
y'_j =
\begin{cases}
y_k^{-1} & \text{if $j = k$};\\[.05in]
\displaystyle y_j \prod_{k\to j} y_k (y_k\oplus 1)^{-1} \prod_{k\ot j} (y_k\oplus 1)
& \text{if $j \neq k$}.
\end{cases}
\end{equation*}
Note that one of the two products is always trivial, hence equal to 1,  since $Q$ has no oriented 2-cycles. Also note
that $\yy'$  depends only on $\yy$ and $Q$.

\item
The quiver $Q'$ is obtained from $Q$ in three steps: 
\begin{enumerate}
\item for every path $i\to k \to j$ add one arrow $i\to j$,
\item reverse all arrows at $k$,
\item delete 2-cycles.
\end{enumerate} 
\end{itemize}

 Mutations are involutions, that is, $\mu_k\mu_k \xyq=\xyq$. 

\begin{definition}
 Let $\calx$ be the set of all cluster variables obtained by mutation from $\xyq$.
The \emph{cluster algebra} $\cala=\Axyq$ is the $\ZZ\PP$-subalgebra of $\calf$ generated by $\calx$.
\end{definition}

Following \cite{FZ4}, we label the seeds $\zS_t=(\xx_t,\yy_t,Q_t)$ by vertices $t$ of an $n$-ary tree whose edges correspond to mutations between the seeds. 
The clusters and coefficient tuples are denoted by
$$\xx_t=(x_{1;t},\ldots, x_{N;t}), \qquad\yy_t=(y_{1;t},\ldots, y_{N;t}),$$ and the Laurent expansion of the cluster variable $x_{i;t}$ in the cluster $\xx_{t'}$ is denoted by $x_{i;t}^{t'}$.
 The initial seed is denoted by $\zS_{t_0}$, or simply by $\zS_0$.

\subsubsection{Laurent phenomenon and positivity}
In the remainder of this subsection, we recall several results on cluster algebras.  We start with Laurent phenomenon and positivity.
\begin{thm}
  Let $u\in \calx$ be a cluster variable and $\zS=\xyq$ any seed. Then \\
\textup{(a)}  \cite{FZ1}\label{thm LP} 
$u$ is a Laurent polynomial in $\xx$, that is,
$u=
f(x_1,\ldots,x_N)/(x_1^{d_1}\cdots x_N^{d_N})$, where $f\in\ZZ\PP[x_1,\ldots,x_N], d_i\in \ZZ$.
\\
\textup{(b)}
\cite{LS6}\label{positivity} the coefficients of $f$ are positive, that is,
$f\in\ZZ_{\ge 0}\PP[x_1,\ldots,x_N]$.
\end{thm}

\subsubsection{Principal coefficients, $F$-polynomials, $g$-vectors and $c$-vectors}
The cluster algebra $\Axyq$ is said to have \emph{principal coefficients} in the initial seed $\zS_0=\xyq$ if $\yy=(y_1,\ldots,y_N)$, that is, the initial coefficient tuple consists of the generators of $\PP$. 
Let $x_{i;t}^{t_0}$ be the cluster expansion of the cluster variable $x_{i;t}$ in the initial seed. This is a Laurent polynomial in $x_1,\ldots,x_N,y_1,\ldots,y_N$. Then the \emph{$F$-polynomial} $F_{i;t}$ is obtained from $x_{i;t}^{t_0}$ by specializing all $x_j$ to 1, that is,
\[F_{i;t}=x_{i;t}^{t_0}|_{x_{j}=1}.\]

\begin{thm}\cite[Theorem 1.7]{DWZ2}
\label{thm fconstant term 1}
Each $F$-polynomial $F_{i;t}$ has constant term 1 and a unique monomial of maximal degree that has coefficient 1 and is divisible by all the other occurring monomials. 
\end{thm}

In particular, the Laurent polynomial $x_{i;t}^{t_0}$ contains a unique monomial without $y$-variables. The degree vector of this monomial is called 
the \emph{$g$-vector} $\mathbf{g}_{i;t}$ of the cluster variable $x_{i;t}$.  

On the other hand, the elements of the coefficient tuple $\yy_t=(y_{1;t},\ldots,y_{N;t})$ of the seed $\zS_t$ can be expressed as monomials in $\PP$ as
\[y_{i;t}=y_1^{c_1}\cdots y_N^{c_N},\]
for some $c_j\in \ZZ$.  The vector  $c_{i;t}=(c_1,\ldots,c_N)$ is called the $i$-th \emph{$c$-vector} of the seed $\zS_t$.

\subsubsection{Separation formulas}
The following formula is sometimes called the separation of addition or the change of coefficients formula. It says that the cluster variable in an arbitrary cluster algebra can be simply computed from the corresponding cluster variable in the cluster algebra with principal coefficients by substituting the new coefficients and clearing the denominators.

The notation $f_{\mathbb{P}}$ stands for the polynomial $f$ over the semifield $\mathbb{P}$. For example if $f=1+y_1y_2^{-1}+y_2y_3^{-3}$  then  $f|_{\mathbb{P}} =1\oplus y_1y_2^{-1}\oplus y_2y_3^{-3}=y_2^{-1}y_3^{-3}$. In our situation, the constant term of the polynomial $f$ is equal to 1, and therefore we see that dividing by $f_{\mathbb{P}}$ is the same as clearing the denominators. In the  example above, $f/f_{\mathbb{P}}=y_2y_3^3+y_1y_3^3+y_2^2.$

\begin{thm}\cite[Theorem 3.7]{FZ4}\label{thmFZ4}
Let $\cala=\cala(\xx,\yy,Q)$ be a cluster algebra with principal coefficients in the initial seed $\zS_{t_0}=(\xx,\yy,Q)$ and let 
 $\overline{\cala}=\cala({\xx},\overline{\yy},Q)$ be a cluster algebra with initial seed $\overline{\zS}_{t_0}=(\xx,\overline{\yy},Q)$ sharing the same cluster $\xx$ and the same quiver $Q$ but with a different coefficient tuple $\overline{\yy}=(\bar{y}_1,\ldots,\bar{y}_N)$. Let $ x_{i;t}^{t_0}$ and $\bar x_{i;t}^{t_0}$ be the Laurent expansions of the $i$-th cluster variable of the seeds $\zS_t$ and $\overline{\zS}_t$ in  $\cala$ and $\overline{\cala}$, respectively. Then 
\[\bar{x}_{i;t}^{t_0}= \frac{x_{i;t}^{t_0}(x_1,\ldots,x_N;\bar{y}_1,\ldots,\bar{y}_N)}
{F_{i;t}(\bar{y}_1,\ldots,\bar{y}_N)|_{\mathbb{P}}}.\]
\end{thm}

Moreover, the cluster variable is determined by its $F$-polynomial and $g$-vector as follows.
\begin{thm}\cite[Corollary 6.3]{FZ4} \label{thmFZ4b}
  \[\bar{x}_{i;t}^{t_0}=\xx^{\mathbf{g}_{i;t}} 
  \frac{F_{i;t}(\widehat{y}_1,\ldots,\widehat{y}_N)}
  {F_{i;t}({\bar y}_1,\ldots,{\bar y}_N)|_{\mathbb{P}}}, \] 
 where $\mathbf{g}_{i;t}$ is the $g$-vector and 
\[\widehat{y}_i=\bar y_i \,\frac{\prod_{i\to j }\, x_{j}}{ \prod_{i\ot j  }\, x_{j}}.\]
\end{thm}
\subsubsection{Sign-coherence and tropical duality}
\begin{thm}
 \cite[Theorem 1.7]{DWZ2}
The $c$-vectors  are sign-coherent, that is, each vector $\mathbf{c}_{i;t}$ has either all entries nonnegative or all entries nonpositive.
\end{thm}

The $c$-vectors and $g$-vectors are dual to each other in the following sense.
Let $C_t^{Q;t_0}$ be the matrix with columns $\mathbf{c}_{1;t},\ldots,\mathbf{c}_{N;t}$ and 
$G_t^{Q;t_0}$ the matrix with columns $\mathbf{g}_{1;t},\ldots,\mathbf{g}_{N;t}$, both taken relative to the quiver $Q$. 
\begin{thm}\label{thm NZ}
 \cite[Theorem 1.2]{NZ}
 \begin{equation*}
 \left(G_{t}^{{Q;t_0}}\right)^{\!\!\mathsf{T}} = \left(C_{t}^{{Q;      t_0}}\right)^{-1}
= C_{t_0}^{Q_{t}^{op};t},\end{equation*}
where $^\mathsf{T}$ denotes the transpose.
\end{thm}

\subsubsection{Generalizations}
Although we do not need it in the current paper, let us mention here that many of the results  above hold in the more general settings of skew-symmetrizable or the even more general setting of totally sign-skew symmetric cluster algebras.

The Laurent phenomenon was proved in \cite{FZ1} for sign-skew symmetric algebras. Positivity and sign-coherence were shown in \cite{GHKK} in the skew-symmetrizable and in \cite{LP} in the sign-skew symmetric case. Theorems~\ref{thmFZ4} and~\ref{thmFZ4b} hold in the skew-symmetrizable case by \cite{FZ4} and a variation of Theorem~\ref{thm NZ} holds in this case as well \cite{NZ}. 

\section{Recollections from Knot Theory and Cluster Algebras I}

We review several concepts and results from our previous article  \cite{BMS}. A \emph{quiver} (or oriented graph) $Q$ is a quadruple $Q=(Q_0,Q_1,s,t)$, where $Q_0,Q_1$ are finite sets and $s,t\colon Q_1\to Q_0$ are maps. The set $Q_0$ is the set of vertices, $Q_1$ the set of arrows, $s(\za)$ is the source of an arrow $\za$ and $t(\za)$ is its target. Let $k$ be a field. 
A \emph{representation} $M=(M_i,\varphi_\za)$ is a family of finite dimensional $k$-vector spaces $M_i,i\in Q_0$ and $k$-linear maps $\varphi_\za\colon M_{s(\za)}\to M_{t(\za)}$. For further details about quivers and their representations, we refer to the textbooks \cite{ASSbook,Sbook}.

\subsection{Quiver and Jacobian algebra} Let $K$ be a link diagram. Denote by $K_0$ the set of crossing points, $K_1$ the set of segments and by $K_2$ the set of regions of the complement of $K$ in the plane, including the unbounded region. Let $n=|K_0|$ be the number of crossing points.  Then $|K_1|=2n$, because every crossing point is incident to four segments, and every segment is incident to two crossing points. Moreover $|K_2|=n+2$, by Euler's formula. 

We associate a quiver $Q=(Q_0,Q_1)$ to $K$ whose vertices are the segments of $K$, thus $Q_0=K_1$. Moreover, the arrows of $Q$ are given by the following rule: \begin{itemize}
    \item [(i)] for every crossing point in $K_0$, whose segments are $a,b,c,d\in K_1$ in clockwise order, draw a clockwise oriented  4-cycle $a\to b\to c\to d\to a$;
    \item[(ii)] remove all 2-cycles.
\end{itemize} 
 See Example~\ref{example 2112} for an illustration.

 The quiver $Q$ thus has clockwise oriented cycles that arise from the crossing points  and counterclockwise oriented cycles from  the regions in $K$. We define a potential $W$ as the sum of the clockwise cycles minus the sum of the counterclockwise cycles. 

 The completed Jacobian algebra $B=\textup{Jac}(Q,W)$ is called the  \emph{Jacobian algebra of the link diagram $K$}. 
\begin{remark}
    For the Jacobian algebra, we may also use the quiver obtained by step (i) only, thus not removing the 2-cycles. In that quiver, the clockwise cycles are all of length four and they are in bijection with the crossing points, and the counterclockwise cycles are in bijection with the regions in $K$. Moreover every arrow lies in precisely two chordless cycles. We have shown in \cite{BMS} that the Jacobian algebras are isomorphic. 
    However, for the definition of the cluster algebra, we need to remove the 2-cycles.  
\end{remark}
\begin{remark}
    The algebra $B$ may be infinite dimensional.
\end{remark}

\subsection{Link diagram module and Alexander polynomial}
To every segment $i\in K_1$, we associate a partition of $K_1=\sqcup_{d\ge 0} K(d)$, where, roughly speaking, $K(d)$ is obtained recursively as a kind of boundary of $K\setminus \cup_{d'<d} K(d')$. We refer to \cite[Section 5.1]{BMS} and it's Addendum \cite{BMS Add} for details.

We further associate to each segment $i\in K_1$, a representation $T(i)$ of the quiver $Q$ and show that $T(i)$ is an indecomposable $B$-module. This representation has the $k$-vector spaces $k^d$ at vertex $j$ precisely if $j\in K(d)$. The linear maps of $T(i)$ are such that their composition along each  chordless cycle in $Q$ is  equal to the Jordan block with diagonal $0$. See \cite[Section 5.3]{BMS} for details. 

Our first result states that the F-polynomial has only coefficients $0$ or $1$.

\begin{thm}\cite[Corollary 6.8]{BMS}
 \label{thm grin}
 For every dimension vector $\mathbf{e}$ the quiver Grassmannian $\textup{Gr}_{\mathbf{e}}(T(i))$ is either empty or a point. 
 In particular,  
the $F$-polynomial of $T(i)$ is \[F_{T(i)}= \sum_{L\subset T(i)} \mathbf{y}^{\underline{\dim}\, L},\]  where the sum is over all submodules of $T(i)$ and $\mathbf{y}^{\underline{\dim}\, L}=\prod_{i=1}^{2n}y_i^{\dim L_i}$.  
\end{thm}

We then show that the Kauffman states correspond to the submodules of $T(i)$ as follows.

\begin{thm}\cite[Theorem 6.9]{BMS}
 \label{thm lattice iso} For every segment $i$ of the link diagram $K$,
 there is a lattice isomorphism from the lattice of Kauffman states of $K$ relative to the segment $i$ to the lattice of submodules of the module $T(i)$. 
\end{thm}
The next theorem is the main result of \cite{BMS}. It follows from the above result using Kauffman's realization of the Alexander polynomial as a weighted sum over all Kauffman states.

We let $T=\oplus_{i\in Q_0} T(i)$ and call it the \emph{link diagram module} of $K$. 
\begin{thm}\label{thmAlexPol}
 The Alexander polynomial of $K$ is equal to the Alexander specialization (\ref{specializationintro}) of the $F$-polynomial of every indecomposable summand $T(i)$ of the link diagram module $T$. That is, for all $i\in Q_0$ \[\zD=F_{T(i)}|_t.\]
\end{thm}

\begin{remark} Note that the specialization in the theorem does not depend on the choice of the segment $i$. Note also that the construction of the quiver does not depend on the overcrossing and undercrossing information in the link diagram. In particular, $Q$ and $\mathcal{A}(Q)$ are not invariants of the link. This information is only recovered in the specialization (\ref{specializationintro}).
\end{remark}

\section{Proof of the main theorem in two examples} \label{sect examples}
Let us start with two examples to illustrate the idea of the proof. Precise definitions will be given in the later sections.

We prove Theorem~\ref{thm main} and Theorem~\ref{thm module main} by constructing a sequence of mutations $\mu$ that transforms the initial seed $\zS_{t_0}=(\xx_0,\yy_0,Q)$ into the seed $\zS_{t}=(\xx_t,\yy_t,Q_t)$ with the desired properties.
The mutation sequence corresponds to a topological construction on the link diagram $K$ that we call \emph{bigon reduction}. A bigon reduction replaces a bigon in $K$ by a single crossing. Let $j,k$ be the labels of the two segments that form the bigon. Then the corresponding step in the cluster algebra is the mutation sequence $\mu_j,\mu_k$. 
We continue performing these  bigon reduction until we obtain the diagram of a Hopf link, that is, a diagram with four segments $a,b,c,d$ consisting of two circles that cross twice. Denote the sequence of mutations so far by $\mu_{\red}$. Note that there are many sequences of bigon reductions leading to the Hopf link. We then perform the mutations in the directions $a,b,c,d$ corresponding to the segments of the Hopf link. Finally we perform the sequence $\mu_{\red}$ in the reverse order. 

The reader will have noticed that this method only works when each diagram in the reduction sequence has a bigon, which is not the case in general. In the general case, we therefore also need another operation on the diagram that will produce bigons. 
In the first example below, the diagrams obtained by this algorithm always admit a bigon, and in the second example, this is not the case.

\subsection{The figure-eight knot}
\begin{figure}
\begin{center}
\scriptsize
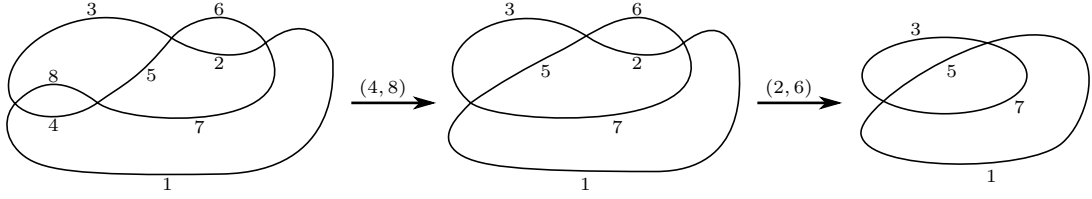
\caption{Bigon reduction on the Figure 8 knot.}
\label{figure8knot}
\end{center}
\end{figure}

Consider this procedure in the example of the figure-eight knot. The bigon reduction is illustrated in Figure~\ref{figure8knot}. The first step is the reduction of the bigon formed by the segments 4 and 8. The second step reduces the bigon formed by segments 2 and 6. After the second step, we obtain a diagram of the Hopf link.
The corresponding mutation sequence is 
\[\mu=\underbrace{\mu_4,\mu_8,\mu_2,\mu_6}_{\mu_{\red}},\underbrace{\mu_1,\mu_5,\mu_3,\mu_7}_{\mu_{\Hopf}},
\underbrace{\mu_2,\mu_6,\mu_4,\mu_8}_{\mu_{\red}^{-1}}\]
where the first four mutations form the reduction sequence $\mu_{\red}$ and the last four mutations are the inverse of that sequence. The middle sequence are the four mutations that correspond to the Hopf link.
The mutation sequence induces a permutation $\zs\in S_8$ as the product of the 2-cycles given by bigons and the pairs in the Hopf sequence. Thus \[\zs=(4\,8)(2\,6)(1\,5)(3\,7).\] 

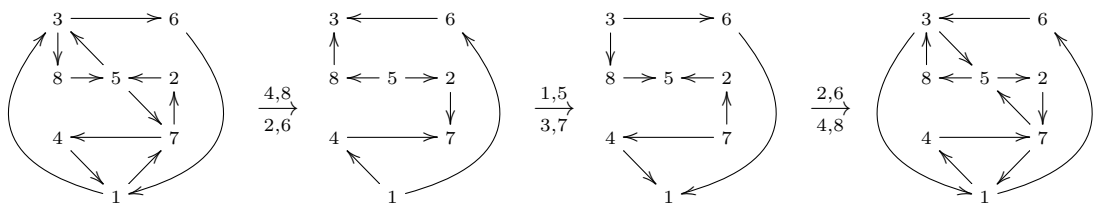
\begin{figure}
\hspace{20pt}\begin{minipage}{0.18\textwidth}
$\scriptsize
\xymatrix@R12pt@C12pt{3 \ar[rr] \ar[d] && 6 \ar@/^30pt/[lddd] \\
8\ar[r]&5\ar[rd]\ar[lu] & 2\ar[l] \\
4\ar[rd]&&7 \ar[ll]\ar[u] \\ 
&1\ar@/^30pt/[luuu] \ar[ru]
}$
\end{minipage}
$\displaystyle \mathrel{\mathop{\longrightarrow}^{\mathrm{4,8}}_{\mathrm{2,6}}} \hspace{5pt}$
\begin{minipage}{0.18\textwidth}
$\scriptsize 
\xymatrix@R12pt@C12pt{3 \ar@{<-}[rr] \ar@{<-}[d] && 6 \ar@{<-}@/^30pt/[lddd] \\
8\ar@{<-}[r]&5 & 2\ar@{<-}[l] \\
4\ar@{<-}[rd]&&7 \ar@{<-}[ll]\ar@{<-}[u] \\ 
&1
}$
\end{minipage}
$\displaystyle \mathrel{\mathop{\longrightarrow}^{\mathrm{1,5}}_{\mathrm{3,7}}} \hspace{5pt}$
\begin{minipage}{0.18\textwidth}
$\scriptsize 
\xymatrix@R12pt@C12pt{3 \ar[rr] \ar[d] && 6 \ar@/^30pt/[lddd] \\
8\ar[r]&5 & 2\ar[l] \\
4\ar[rd]&&7 \ar[ll]\ar[u] \\ 
&1
}$
\end{minipage}
$\displaystyle \mathrel{\mathop{\longrightarrow}^{\mathrm{2,6}}_{\mathrm{4,8}}} 
\hspace{20pt}$
\begin{minipage}{0.18\textwidth}
$\scriptsize 
\xymatrix@R12pt@C12pt{3 \ar@{<-}[rr] \ar@{<-}[d] && 6 \ar@{<-}@/^30pt/[lddd] \\
8\ar@{<-}[r]&5\ar@{<-}[rd]\ar@{<-}[lu] & 2\ar@{<-}[l] \\
4\ar@{<-}[rd]&&7 \ar@{<-}[ll]\ar@{<-}[u] \\ 
&1\ar@{<-}@/^30pt/[luuu] \ar@{<-}[ru]
}$
\end{minipage}
\caption{Four quivers along the mutation sequence for the figure-eight knot. The first quiver is the original quiver $Q$, the second is obtained after the first four mutations at $4,8,2,6$, the next by mutations at 1,5,3,7 and the last is the final quiver $Q_t$. }
\label{figquiver8}
\end{figure}

The effect of the mutation sequence on the quiver $Q$ is illustrated in Figure~\ref{figquiver8}. The initial quiver $Q$ is shown on the left and the quiver $Q_t=\mu Q$ is shown on the right. We also show the intermediate quivers obtained after mutation subsequences of length four. 

We see that $\zs$ induces an isomorphism of quivers $\zs\colon Q\to Q^{op}_t$. 
Actually, in this particularly symmetric example, the identity permutation also induces a quiver isomorphism, but this is not the case in general. Another exceptional feature of this example is that the quiver $Q$ is mutation equivalent to an acyclic quiver, as we can see from the two quivers in the middle of Figure~\ref{figquiver8}, which are of extended Dynkin type $\widetilde{\mathbb{A}}_{4,4}$.

In the cluster algebra $\cala(Q)$, the above sequence of mutations $\mu$ produces a seed $\zS_t$  whose $F$-polynomials are the following.
\[
\begin{array}
 {lcl}
F_1=1+y_1+y_1y_4+y_1y_6 + y_1 y_4 y_6,
&&F_6=1 + y_8 + y_3 y_8 + y_1 y_3 y_8 + y_1 y_3 y_4 y_8,
\\F_5=1 + y_5 + y_2y_5 + y_5 y_8 + y_2 y_5 y_8,
&&F_2=1 + y_4 + y_4 y_7 + y_4 y_5 y_7 + y_4 y_5 y_7 y_8,
\\F_3=1 + y_6 + y_8 + y_6 y_8 + y_3 y_6 y_8,
&&F_4=1 + y_2 + y_2 y_7 + y_1 y_2 y_7 + y_1 y_2 y_6 y_7,
\\F_7=1 + y_2 + y_4 + y_2 y_4 + y_2 y_4 y_7,
&&F_8=1 + y_6 + y_3 y_6 + y_3 y_5 y_6 + y_2 y_3 y_5 y_6.
\end{array}\]

These polynomials $F_i$ are precisely the $F$-polynomials of our representations $T(\zs (i))$ in  \cite{BMS}. We list the representations below.
\[
\begin{array}
 {lcl}
T(5)= 
\begin{smallmatrix}
 4\ 6\\1
\end{smallmatrix}
&\qquad&
T(2)= 
\begin{smallmatrix}
 4\\1\\3\\8
\end{smallmatrix}
\\ [10pt]
T(1)= 
\begin{smallmatrix}
 2\ 8\\5
\end{smallmatrix}
&&T(6)= 
\begin{smallmatrix}
 8\\5\\7\\4
\end{smallmatrix}
\\ [10pt]
T(7)= 
\begin{smallmatrix}
 3\\8\ 6
\end{smallmatrix}
&&
T(8)= 
\begin{smallmatrix}
 6\\1\\7\\2
\end{smallmatrix}
\\ [10pt]
T(3)= 
\begin{smallmatrix}
 7\\2\ 4
\end{smallmatrix}
&&T(4)= 
\begin{smallmatrix}
2\\5\\3\\6
\end{smallmatrix}
\end{array}
\]

\begin{remark}
    In this example each mutation is a green mutation. If we append the two mutations at 3 and 7 to the sequence we obtain  a maximal green sequence.
\end{remark}

\subsection{The Borromean rings}\label{sect Borromean}
\begin{figure}
\begin{center}
\scriptsize
\scalebox{.9}{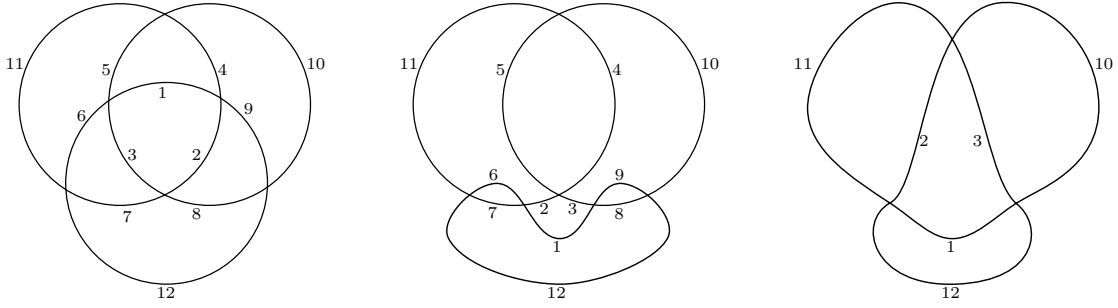}
\caption{Diagram R3 move on Borromean Rings and subsequent bigon reductions}
\label{figBorroemeanRings}
\end{center}
\end{figure}

The Borromean rings form a 3-component link whose diagram is shown in the left picture in Figure~\ref{figBorroemeanRings}. This is the smallest example of a link that does not have a bigon. Our first step is therefore to transform the diagram into one that does contain bigons. This operation is similar to the third Reidemeister move in that it moves a segment over a crossing point; however it is not a Reidemeister move, because we are ignoring the over/under crossing information of the link. It is important to note that this move does not preserve the isomorphism type of the link. We call this move a diagram R3 move or RD3 move for short.

In our example, we move the segment labeled 1 over the crossing point formed by the segments 2 and 3. The result is shown in the center picture of Figure~\ref{figBorroemeanRings}.  This diagram has three bigons $(4,5),(6,7),$ and $(8,9)$. The picture on the right of the figure shows the diagram obtained by reducing these three bigons. Next, we reduce the bigon $(1,12)$ and we obtain a Hopf link on segments 2,10,3,11. 

In the cluster algebra, we build our mutation sequence $\mu$ in several steps. First, we associate to our diagram R3 move
the mutation sequence at vertices 1,2,3,1,2,3,2,3,2. 
The subsequent bigon reductions then correspond to mutations at vertices 4,5,6,7,8,9,1,12. Denote the sequence of all these mutations by $\mu_{\red}$ so far. The Hopf link gives us the sequence $\mu_{\textup{Hopf}}$ at vertices 2,10,3,11. 
As in our previous example, we define a permutation $\zs\in S_{12}$ as the product of 2-cycles given by the bigons and the pairs in the Hopf sequence. Thus
\[\zs =(4,5)(6,7)(8,9)(1,12)(2,10)(3,11).\]

We now define the mutation sequence $\mu$ as $\mu_{\textup{red}}\,\mu_{\textup{Hopf}}\,{\overleftarrow{\mu_{\textup{red}}^\zs}}$, where superscript $\zs$ means that a mutation $\mu_j$ is replaced by $\mu_{\zs (j)}$ and the arrow denotes the mutation sequence in reverse order. 
Thus $\mu$ is the mutation sequence at the vertices
{\small \[\textstyle(1,2,3,1,2,3,2,3,2),
(4,5,6,7,8,9, 1,12),
( 2,10,3,11),
 (1,12,8,9,6,7,4,5),
 (10,11,10,11,10,12,11,10,12).
\]}

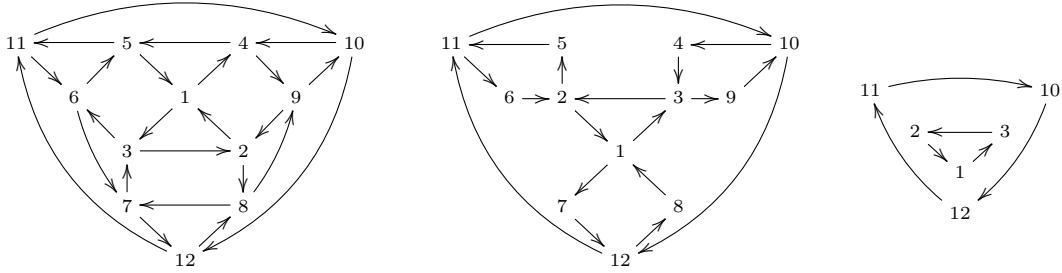
\begin{figure}
\begin{center}
\begin{minipage}{0.3\textwidth}
 $\scriptsize\xymatrix@R10pt@C10pt{11\ar[rd]\ar@/^15pt/[rrrrrr]&&5\ar[ll]\ar[rd]&&4\ar[ll]\ar[rd]&&10\ar[ll]\ar@/^15pt/[llldddd]
\\
&6\ar[ru]\ar@/_3pt/[rdd]&&1\ar[ru]\ar[ld]&&9\ar[ru]\ar[ld]
\\
&&3\ar[lu]\ar[rr]&&2\ar[lu]\ar[d]\\
&&7\ar[u]\ar[rd]&&8\ar[ll]\ar@/_3pt/[ruu]\\
&&&12\ar[ru]\ar@/^15pt/[llluuuu]
}$
\end{minipage}
\hspace{30pt}
\begin{minipage}{0.3\textwidth}
\scriptsize\xymatrix@R10pt@C10pt{11\ar[rd]\ar@/^15pt/[rrrrrr]&&5\ar[ll]&&4\ar[d]&&10\ar[ll]\ar@/^15pt/[llldddd]
\\
&6\ar[r]&2\ar[u]\ar[rd]&&3\ar[r]\ar[ll]&9\ar[ru]
\\
&&&1\ar[ld]\ar[ru]\\
&&7\ar[rd]&&8\ar[lu]\\
&&&12\ar[ru]\ar@/^15pt/[llluuuu]
}
\end{minipage}
\hspace{10pt}
\begin{minipage}{0.3\textwidth}
\scriptsize\xymatrix@R5pt@C5pt{
11\ar@/^5pt/[rrrr]&&&&10\ar@/^5pt/[llddd]\\
&2\ar[rd]&&3\ar[ll]\\
&&1\ar[ru]\\
&&12\ar@/^5pt/[lluuu]
}
\end{minipage}
\caption{The quivers of the diagrams in Figure~\ref{figBorroemeanRings}}
\label{quiverBor}
\end{center}
\end{figure}

The quivers of the diagrams of Figure~\ref{figBorroemeanRings} are shown in Figure~\ref{quiverBor}. The quiver in the center is obtained from the quiver on the left by the mutation sequence at 1,2,3,1,2,3,2,3,2, thus the sequence corresponding to the diagram R3 move. 

\begin{remark}
    One way to think about this mutation sequence is to imagine the cluster category of the full subquiver of $Q$ with vertices $1,2,3$. This quiver is an oriented 3-cycle and the corresponding cluster category is the category of all diagonals in a regular hexagon, see \cite{CCS}. The three segments $1,2,3$ correspond to an internal triangle in this hexagon.  Each mutation is a flip of the corresponding diagonal. The first four mutations in $1,2,3,1$ transform our original internal triangle into the unique other internal triangle. This is what we want. However, we do not like the positions of the labels $2$ and $3$. Therefore we add the 5 mutations in $2,3,2,3,2$ which have the effect of interchanging the diagonals $2$ and $3$. 
\[\xymatrix@C10pt@R10pt@!{
 { \xy/r1pc/: {\xypolygon6"A"{~<<{@{}}~><{@{-}}
~>>{_{}}}},
\POS"A1" \ar@{-} "A3",\POS"A3" \ar@[blue]@{-} "A5",\POS"A5" \ar@[red]@{-} "A1",\endxy}\ar[r]^{\color{blue} 1}
&
{ \xy/r1pc/: {\xypolygon6"A"{~<<{@{}}~><{@{-}}
~>>{_{}}}},
\POS"A1" \ar@{-} "A3",\POS"A1" \ar@[blue]@{-} "A4",\POS"A5" \ar@[red]@{-} "A1",\endxy}\ar[r]^{\color{red}2}
&
{ \xy/r1pc/: {\xypolygon6"A"{~<<{@{}}~><{@{-}}
~>>{_{}}}},
\POS"A1" \ar@{-} "A3",\POS"A1" \ar@[blue]@{-} "A4",\POS"A6" \ar@[red]@{-} "A4",\endxy}\ar[r]^3
&
{ \xy/r1pc/: {\xypolygon6"A"{~<<{@{}}~><{@{-}}
~>>{_{}}}},
\POS"A2" \ar@{-} "A4",\POS"A1" \ar@[blue]@{-} "A4",\POS"A6" \ar@[red]@{-} "A4",\endxy}\ar[r]^{\color{blue} 1}
&
{ \xy/r1pc/: {\xypolygon6"A"{~<<{@{}}~><{@{-}}
~>>{_{}}}},
\POS"A2" \ar@{-} "A4",\POS"A2" \ar@[blue]@{-} "A6",\POS"A6" \ar@[red]@{-} "A4",\endxy}\ar[r]^{\color{red}2}
&
{ \xy/r1pc/: {\xypolygon6"A"{~<<{@{}}~><{@{-}}
~>>{_{}}}},
\POS"A2" \ar@{-} "A4",\POS"A2" \ar@[blue]@{-} "A6",\POS"A2" \ar@[red]@{-} "A5",\endxy}\ar[r]^3
&
{ \xy/r1pc/: {\xypolygon6"A"{~<<{@{}}~><{@{-}}
~>>{_{}}}},
\POS"A3" \ar@{-} "A5",\POS"A2" \ar@[blue]@{-} "A6",\POS"A2" \ar@[red]@{-} "A5",\endxy}\ar[r]^{\color{red}2}
&
{ \xy/r1pc/: {\xypolygon6"A"{~<<{@{}}~><{@{-}}
~>>{_{}}}},
\POS"A3" \ar@{-} "A5",\POS"A2" \ar@[blue]@{-} "A6",\POS"A3" \ar@{-}@[red] "A6",\endxy}\ar[r]^3
&
{ \xy/r1pc/: {\xypolygon6"A"{~<<{@{}}~><{@{-}}
~>>{_{}}}},
\POS"A4" \ar@{-} "A6",\POS"A2" \ar@[blue]@{-} "A6",\POS"A3" \ar@[red]@{-} "A6",\endxy}\ar[r]^{\color{red}2}
&
{ \xy/r1pc/: {\xypolygon6"A"{~<<{@{}}~><{@{-}}
~>>{_{}}}},
\POS"A4" \ar@{-} "A6",\POS"A2" \ar@[blue]@{-} "A6",\POS"A2" \ar@[red]@{-} "A4",\endxy}
}\]
\end{remark}

In the cluster algebra, the mutation sequence $\mu $ produces the seed $\zS_t$ whose $F$-polynomials $F_j$ are precisely the $F$-polynomials $F_{T(\zs (j))}$ of the representations $T(\zs (j))$ constructed in \cite{BMS}. For example 
\[\begin{array}{l}
F_1=
1 + y_1 
+ y_1 y_2 
+ y_1 y_5 
+ y_1 y_2 y_5 
+ y_1 y_5 y_6  
+ y_1 y_2 y_9 
+ y_1 y_2 y_5 y_6 
+ y_1 y_2 y_5 y_9 
 + y_1 y_2 y_3 y_5 y_6 \\
\  + y_1 y_2 y_5 y_6 y_9 
 + y_1 y_2 y_4 y_5 y_9  
 + y_1 y_2 y_3 y_5 y_6 y_9 
 + y_1 y_2 y_4 y_5 y_6 y_9 
 + y_1 y_2 y_3 y_4 y_5 y_6 y_9 
 + y_1^2 y_2 y_3 y_4 y_5 y_6 y_9,
 \end{array}\]
  which is the $F$-polynomial of the representation
 \[T(\zs (1))=T(12)=
\begin{smallmatrix}
 1\\3\ 4\\6\ 9\\5\ 2\\1 
\end{smallmatrix}.\]

\begin{remark}
    In this example, the mutation is a reddening sequence. Moreover, if we remove the extra mutations in 2,3,$\zs (2),\zs (3)$, we obtain  the maximal green sequence 1,2,3,1,\ 4,5,6,7,8,9,1,12,\ 2,10,3,11,\ 1,12,8,9,6,7,4,5,\ 12,11,10,12.
\end{remark}

\section{Proof of the main theorem}\label{sect proof}
In this section we prove the Theorems~\ref{thm main} and \ref{thm module main} in the general case.
Let $K$ be a diagram of a prime link without curls. 

\subsection{The mutation sequence}\label{sect mutation} We define a sequence of mutations $\mu$ as follows. For examples, see section~\ref{sect examples}. If $K$ does not contain a bigon, then, by Theorem~\ref{thm::sequenceR3 intro}, there exists  an admissible sequence
$R^0_1,R^0_2,\ldots, R^0_{r_0}$ of diagram R3 moves
such that the resulting diagram $\overline K_0$ has a bigon $(j_0,k_0)$. If $K$ already has a bigon, let $\overline K_0=K$. 
Let $K_1=\overline K_0\setminus (j_0,k_0)$ be the diagram obtained from $
\overline K_0$ by reducing the bigon $(j_0,k_0)$.  

Recursively, given $K_i$, there exists  an admissible sequence  $R^i_1,R^i_2,\ldots, R^i_{r_i}$ of diagram R3 moves
such that the resulting diagram $\overline K_i$ has a bigon $(j_i,k_i)$. If $K_i$ already has a bigon, let $\overline K_i=K_i$. Let $K_{i+1}=\overline K_i\setminus (j_i,k_i)$ be the diagram obtained from $
\overline K_i$ by reducing the bigon $(j_i,k_i)$.  

 We need the following result. 
\begin{prop}
    \label{prop tjs premier}
   (a) A bigon reduction of a prime link $K$ is again a prime link.

(b) If a prime link $K$ with $n$ crossing points admits a sequence of $n-2$ bigon reductions then the result is a Hopf link. 
\end{prop}
\begin{proof}
 (a) Suppose $K'$ is obtained by a bigon reduction of a prime link diagram $K$ such that $K'$ is not prime. Then there exists a simple closed curve $S$ in the plane that cuts $K'$ into two nontrivial components and such that $K'$ and $S$ intersect in precisely two points. The crossing point of the reduction must lie in one of the two components of $K'$. Then we can put the bigon back in without creating additional crossings with $S$. Thus $K$ is not prime, a contradiction.

(b) Since $K$ is initially prime, if $K$ is not already a Hopf link, a bigon reduction cannot create a curl.
 Therefore, after $n-3$ bigon reductions, we reach the diagram of the Hopf link, since each step reduces the number of crossing points by one.
\end{proof}

 Denote the segments of the Hopf link obtained by the bigon reductions by $a,b,c,d$ such that the segments $a,b$ form one component and $c,d$ form the other. We define the following pairs $(j_{n-2},k_{n-2})=(a,b),(j_{n-1},k_{n-1})=(c,d)$.

We let $\zs$ be the permutation 
\[\zs=(j_0\ k_0)(j_1\ k_1)\cdots(j_{n-1}\ k_{n-1})\in S_{2n}.\]
To every diagram R3 move $R$ we associate a sequence of mutations as follows. Suppose $R$ moves a segment $a$ through a crossing point formed by two segments $b$ and $c$ such that $b$ is clockwise from $c$ when going around that crossing point. Then we let $\mu_R$ be the sequence of mutations at vertices $a,b,c,a,b,c,b,c,b$ in that order.

Thus for each diagram R3 move $R$, we have a corresponding mutation sequence $\mu_R$. Next, we define 
 for each of the sequences $R^i_1,R^i_2,\ldots, R^i_{r_i}$, $i=0,1,\ldots,n-1$ of diagram R3 moves above 
 the sequence of mutations $\mu_{R^i}$ to be the composition 
 \[\mu_{R^i}=\mu_{R^i_1}\mu_{R^i_2}\ldots \mu_{R^i_{r_i}}\]
 where we the order of the mutations is from left to right.
 Then we define a reduction sequence 
 \[\mu_{\textup {red}}=\mu_{R^0}\mu_{j_0}\mu_{k_0}
 \mu_{R^1}\mu_{j_1}\mu_{k_1} \cdots
 \mu_{R^{n-3}}\mu_{j_{n-3}}\mu_{k_{n-3}}
\] 
 and a Hopf link sequence
 \[\mu_{\textup {Hopf}} =
 \mu_{j_{n-2}}\mu_{k_{n-2}}
 \mu_{j_{n-1}}\mu_{k_{n-1}}
 \]
 and finally our mutation sequence
 \begin{equation}\label{eq mutation sequence}
 \mu=\mu_{\textup {red}}\ \mu_{\textup {Hopf}} \ 
 \overleftarrow{\mu^\zs_{\textup {red}}}
 \end{equation}
 where the superscript $\zs$ indicates that each mutation $\mu_k$ in the original sequence is replaced by the mutation $\mu_{\zs(k)}$ and
 the arrow means that we take the reverse of the sequence. 
 Thus $\nu\overleftarrow{\nu}$ is the identity for every mutation sequence $\nu$. 

\emph{Proof outline.}
Let $\zS_t$ be the seed obtained from the initial seed $\zS_0$  by the mutation sequence $\mu$. We may assume without loss of generality that the sequence $\mu$ has a minimal length, and we proceed by induction on the length.

 The smallest possible case is when $K$ is the Hopf link.
In this case, the quiver of the cluster algebra consists of four vertices and no arrows (after we have removed the 2-cycles). The sequence $\mu$ mutates at each vertex exactly once and produces the $F$-polynomials $1+y_i$, $i=1,2,3,4$. On the other hand, the module $T(\zs (i))$ is the simple module $S(i)$ at vertex $i$, and thus we have
$F_{T(\zs (i))}=1+y_i$. This proves the result for the Hopf link.

For the remainder of this section, suppose $K$ is not the Hopf link. 
There are two cases to consider depending on whether or not the diagram $K$ has a bigon. 
The case where $K$ does not contain a bigon is done in section~\ref{proof w/o bigon} and the case where $K$ has a bigon in section~\ref{proof bigon}.

\subsection{Proof in the case where \texorpdfstring{$K$}~ has no bigon}\label{proof w/o bigon}

Suppose  
that $K$ does not contain a bigon, so that 
$\mu=\mu_{R^0_1}\mu'$ 
starts with a sequence of a diagram R3 move $R^0_1$.
 
\begin{figure}
\begin{center}
{\Large\scalebox{0.8}{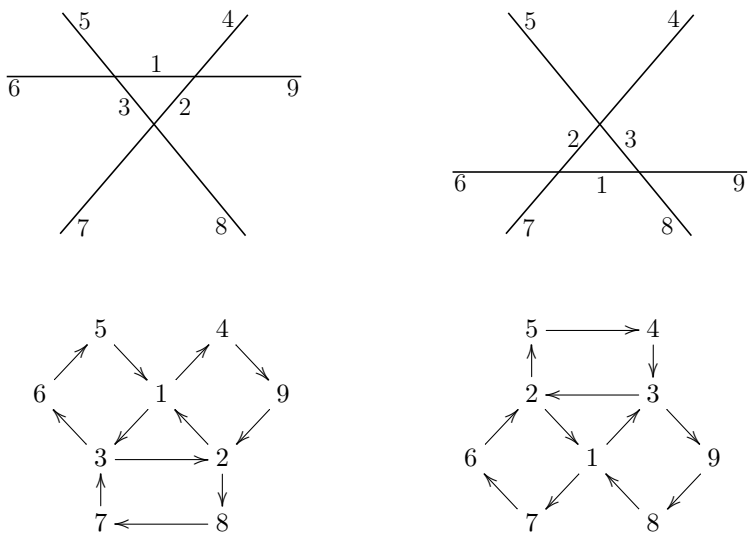}}

\vspace{10pt}
\[ \xymatrix@=1.2em{
&5\ar[rd]&&4\ar[rd]\\
6\ar[ru]&&1\ar[ld]\ar[ru]&&9\ar[ld]\\
&3\ar[lu]\ar[rr]&&2\ar[lu]\ar[d]\\
&7\ar[u]&&8\ar[ll]
}
\qquad\qquad\qquad
\xymatrix@=1.2em{
&5\ar[rr]&&4\ar[d]\\
&2\ar[u]\ar[rd]&&3\ar[ll]\ar[rd]\\
6\ar[ru]&&1\ar[ld]\ar[ru]&&9\ar[ld]\\
&7\ar[lu]&&8\ar[lu]
}
\]
\caption{The diagram R3 move from $K$ (top left) to $\overline K$ (top right) and the corresponding quivers $Q$ (bottom left) and $\overline{Q}$ (bottom right).}
\label{fig R3}
\end{center}
\end{figure}

 Let us denote the segments of $K$ that are involved in this move ${R^0_1}$ by $1,2,3,4,5,6,7,8,9$ as shown in the top left picture in Figure~\ref{fig R3}. 
Let $\overline{K}$ denote the diagram obtained from $K$ by the move ${R^0_1}$ as shown in the top right picture of the figure. We will write $\mu_0$ instead of $\mu_{R^0_1}$ for short. Thus $\mu_0$ is the sequence of mutations at vertices $1,2,3,1,2,3,2,3,2$ in order from left to right. The quivers $Q$ and $\overline{Q}$ of the two diagrams are shown in the bottom pictures of the same figure. It is easily verifiable by hand or using \cite{mutationapp} that $\overline{Q}=\mu_0 Q$. 
 
Thus the mutation sequence $\mu$ can be decomposed into $\mu=\mu_0\mu_1 {\overleftarrow{\mu^\zs_0}}$, where ${\overleftarrow{\mu^\zs_0}}$ is the mutation sequence at the vertices $\zs(2),\zs(3), \zs(2) ,\zs(3) ,\zs(2), \zs(1) ,\zs(3) ,\zs(2), \zs(1) $ in order from left to right. We illustrate this decomposition in the following diagram of seeds and mutations in the cluster algebra.
 
\begin{equation}\label{eq mutation diagram}
 \xymatrix@C55pt{\zS_0\ar[r]^{\mu_0} &\zS_{t_1}\ar[r]^{\mu_1}&\zS_{t_2}\ar[r]^{\overleftarrow{\mu^\zs_0}}&\zS_t.}
\end{equation}
 
  We start by proving that $\zs\colon Q_0\to Q_t^{\textup{op}}$ is an isomorphism of quivers. By induction, we know that $\bar\zs\colon Q_{t_1}\to Q_{t_2}^{\textup{op}}$ is an isomorphism. The quiver $Q_{t_0}$ is obtained from $Q_{t_1}$ by the mutation sequence at vertices 2,3,2,3,2,1,3,2,1 in order. 
 On the other hand, $Q_{t}$ is obtained from $Q_{t_2}$ by the mutation sequence at vertices $\zs (2),\zs (3),\zs (2),\zs (3),\zs (2),\zs (1),\zs (3),\zs (2),\zs (1)$ in order. Since $\zs=\bar\zs$, this shows that $\zs\colon Q_0\to Q_t^{\textup{op}}$ is an isomorphism. This shows part (b) of Theorem~\ref{thm main}.

Our next goal is to prove Theorem~\ref{thm module main}. 
Once this is done,  part (a) of Theorem~\ref{thm main} will follow immediately using Theorem~\ref{thmAlexPol}.
 By induction, we know that for every vertex $i$ of $\overline{Q}_0$, we have
 
\begin{equation}
 \label{Induction hypothesis} 
 F_{i;t_2}^{t_1} = F_{\overline{T}(\overline {\zs}(i))},
\end{equation}
 where
 
\begin{itemize}
 \item[-] $F_{i;t_2}^{t_1}$ is the $F$-polynomial of the cluster variable $x_{i;t_2}$ in the cluster algebra $\cala(\overline{Q})$ relative to the initial seed $\zS_{t_1}=(\xx_{t_1},\overline{\yy}, \overline{Q})$ with principal coefficients $\overline{\yy}=(\bar y_1,\ldots,\bar y_{2n})$;
 \item[-] $\overline{\zs}$ is the permutation corresponding to the mutation sequence $\mu_1$; (Note however that $\overline{\zs}=\zs$, since the sequence $\mu_0$ does not contain any bigon reductions.)
 \item [-]  ${\overline{T}(\overline {\zs}(i))}$ is the representation of $\overline{Q}$ corresponding to the segment $\overline{\zs}(i)=\zs(i)$ of $\overline{K}$; and
 \item [-] $ F_{\overline{T}(\overline {\zs}(i))}$ is the $F$-polynomial of $ {\overline{T}(\overline {\zs}(i))}$.
\end{itemize}

We want to show that
\begin{equation}
\label{eq goal}
 F_{i;t}^{t_0} = F_{T(\zs(i))}.
\end{equation}

The left hand side of this equation is the $F$-polynomial  of the cluster variable $x_{i,t}$ in the cluster algebra $\cala(Q)$ relative to the initial seed  $\zS_{t_0}=(\xx_{0},\yy_0, Q)$ with principal coefficients $\yy=( y_1,\ldots, y_{2n})$.

\medskip

\textbf{First case.} \label{sect phi def} We will first assume that $i \notin \{\zs(1),\zs(2),\zs(3)\}$.
In this case, $x_{i;t}=x_{i;t_2}$, because  the sequence $\overleftarrow{\mu_0^\zs}$ only mutates at the vertices $\zs(1),\zs(2),\zs(3)$. Thus $F_{i;t}^{t_0}=F_{i;t_2}^{t_0}$, and we need to compute $F_{i;t_2}^{t_0}$.

We obtain $F_{i;t_2}^{t_0}$ from $F_{i;t_2}^{t_1}$ in four steps. 
\begin{enumerate}
 \item Calculate the cluster variable $\overline{x}_{i;t_2}^{t_1}=\xx_{t_1}^{\overline{g}_{i;t_2}} F_{i;t_2}^{t_1}(\widehat{\bar y_i}) $ in the cluster algebra $\cala(\overline{Q})$ using the formula in Theorem~\ref{thmFZ4b}, where $\overline{g}_{i;t_2}=(\overline{g}_j)_{j=1}^{2n} $ denotes the $g$-vector and \[\widehat{\bar y_i}=\frac{\bar y_i\prod_{i\to j \in \overline{Q}} x_{j;t_1}}{\prod_{j\to i \in \overline{Q}} x_{j;t_1}};\]
 \item Calculate $x_{i;t_2}^{t_1}$ in the cluster algebra $\cala(Q)$ using the Fomin-Zelevinsky change of coefficients formula, see Theorem~\ref{thmFZ4}, 
 by replacing the coefficients $\bar y_j$ in $\bar x_{i;t_2}^{t_1}$ by the coefficients $y_{j;t_1}$ that are obtained via the mutation sequence $\mu_0$ from the seed $\zS_0$, and subsequently clearing the denominators by dividing by $F_{i;t_2}^{t_0}(y_{j;t_1}) \vert_{\mathbb{P}}$; 
 \item Calculate $x_{i;t_2}^{t_0}$ by replacing the cluster  variables $x_{j;t_1}$ in  $x_{i;t_2}^{t_1}$ by their Laurent expansions 
 $x_{j;t_1}^{t_0}$ in $\zS_0$. 
 \item Calculate  $F_{i;t_2}^{t_0}$ from $x_{i;t_2}^{t_0}$ by specializing all cluster variables $x_{j,t_0}$ to 1.
\end{enumerate}

\begin{remark}
 \label{rem steps 34}
The  steps (3) and (4) combined have the effect of replacing each cluster variable $x_{j;t_1}$ by its $F$-polynomial $F_{j;t_1}^{t_0}$.
 Moreover, since the mutation sequence $\mu_0$ only mutates in vertices 1,2 and 3, we have $F^{t_0}_{j;t_1}=1$, for $j>3$. 
 We will often use the notation $F_j=F^{t_0}_{j;t_1}$, for $j=1,2,3$.
\end{remark}

\begin{remark} \label{rem g-vector}
The term $\xx_{t_1}^{\overline{g}_{i;t_2}} = \prod_{j} (x_{j;t_1})^{\overline{g}_j}$ of step (1) is a monomial in the variables $x_{j;t_1}$.  According to Remark~\ref{rem steps 34}, in steps (3) and (4) this monomial becomes $F_1^{\overline{g}_1}F_2^{\overline{g}_2}F_3^{\overline{g}_3}$.  In particular we only need to know the first three components of the $g$-vector.
\end{remark}

In view of the above remarks, we can now describe the process of going through the steps (1)-(4) via a ring homomorphism $\varphi$ as follows.

\begin{definition}
 \label{def phi}
 Let $\varphi\colon\mathbb{Z}[\bar y_1,\ldots,\bar y_{2n}] \to \mathbb{Z}[y_1,\ldots,y_{2n}]$ be the ring homomorphism such that 
   \[ F_{i;t_2}^{t_0} =  \varphi( F_{i;t_2}^{t_1})
   \ \frac{\prod_{j =1}^3  (F_{j;t_1}^{t_0})^{\overline{g}_j} }
{F_{i;t_2}^{t_0}(y_{j;t_1}) \vert_{\mathbb{P}}}
.\]
\end{definition}

Lemma \ref{lem phi} defines  this homomorphism explicitly.

Using this notation and our {induction hypothesis} (\ref{Induction hypothesis}), we see that 
 in order to show (\ref{eq goal}), it suffices to show 
\begin{equation}
\label{eq goalphi} 
\varphi\left(F_{\overline{T}(\zs(i))}\right) 
   \ \frac{\prod_{j =1}^3  (F_{j;t_1}^{t_0})^{\overline{g}_j} }
{ F_{i;t_2}^{t_0}(y_{j;t_1}) \vert_{\mathbb{P}}}
= F_{T(\zs(i))}.
\end{equation}

\subsubsection{Computation of $\varphi$: Preparatory lemmas}
 Throughout this section we use the following notation
\[\begin{array}{rcl}F_1&=&1+y_2+y_2y_3 \\
F_2&=&1+y_3+y_1y_3 \\
F_3&=&1+y_1+y_1y_2 \end{array}\]
which we sometimes write as 
$F_a=1+y_{\rho (a)}+y_{\rho(a)}y_{\rho^2(a)}$, where $\rho$ is the permutation $(123)$. We also set $F_a=1$, for $a>3$.
It is easy to check that 
\begin{equation}
 F_a= F^{t_0}_{a;t_1}.
\end{equation}
The following relations among the $F_a$ are immediate from the definition.
\begin{lemma}
 \label{frelations}
 With the notation above, we have $y_aF_a+F_{\rho(a)}=(1+\rho^2(a))F_{\rho^2(a)}$, for all $a =1,2,3$. That is  
 \[\begin{array}{rcl}
 y _1F_1+F_2 &=&(1+y_3)F_3 \\
 y_2 F_2+F_3 &=&(1+y_1)F_1 \\
 y_3 F_3+F_1 &=&(1+y_2)F_2 .
 \end{array}\]
\end{lemma}

 The next lemma describes $\varphi $ on the variables $\bar y_j$. Note that this information completely determines $\varphi$, because $\varphi$ is a ring homomorphism.
\begin{lemma}\label{lem phi}
The values of $\varphi$ on $\bar y_j$ are as follows.
\[
\begin{array}{lll} 
 \varphi(\bar y_1)=\displaystyle\frac{F_3}{y_2F_2} \quad&   \varphi(\bar y_4)=y_4F_3 &  \varphi(\bar y_7)=\displaystyle\frac{y_2y_3y_7}{F_1} \\[15pt]
  \varphi(\bar y_2)=\displaystyle\frac{F_1}{y_3F_3 }&\varphi(\bar y_5)=\displaystyle\frac{y_1y_3y_5}{F_2} \quad& \varphi(\bar y_8)=y_8F_1\\[15pt]
   \varphi(\bar y_3)=\displaystyle\frac{F_2}{y_1F_1}
  & \varphi(\bar y_6)=y_6F_2  &   \varphi(\bar y_9)=\displaystyle\frac{y_1y_2y_9}{F_3}
  \end{array}
\]
and  $\varphi(\bar y_j)=y_j$, if $j>9$.
 \end{lemma}
\begin{proof}
 The whole proof is an exercise in mutations. 
 Recall that $\varphi $ is defined by the steps (1)-(4) above. Thus we need to compute
 $\widehat {\bar y}_j$, then replace the $\bar y_j $ by $y_{j;t_1}$, and then replace each $x_{k;t_1}$ with $k=1,2,3$ by  $F_k$, and each $x_{k;t_1}$ with $k>3$ by 1. The intermediate steps are given below for $j=1,4,$ and $5$. The other cases are symmetric.
 
 For $j=1$, we have
 $\widehat {\bar y}_1=\bar y_1 x_{3;t_1}x_{7;t_1}x_{2;t_1}^{-1}x_{8;t_1}^{-1}$ and $\bar y_1=y_2^{-1}$. Thus $\varphi(\bar y_1) =y_2^{-1} F_3F_2^{-1}$. 
 
 For $j=4$, we have
 $\widehat {\bar y}_4=\bar y_4 x_{3;t_1}$ and $\bar y_4=y_4$. Thus $\varphi(\bar y_4) =y_4 F_3$. 
 
 For $j=5$, we have
 $\widehat {\bar y}_5=\bar y_5 x_{2;t_1}^{-1}$ and $\bar y_5=y_1y_3y_5$. Thus $\varphi(\bar y_5) =y_1y_3y_5F_2^{-1}$.
\end{proof}

The next two lemmas provide formulas that will be useful in the next subsection. Here we use the notation $\bar F_a$ for the polynomial $F_a$ after substituting the variables $y_j$ by $\bar y_j$. 
\begin{lemma} 
 \label{lem phif}
 With the notation above, we have
$ \varphi(\bar F_a) = F_{\rho(a)} / (y_ay_{\rho^2(a)})$, for $a=1,2,3$. That is  
 \[\begin{array}{ccc}
 \varphi(\bar F_1) = \displaystyle\frac{F_2}{y_1y_3} &\quad
 \varphi(\bar F_2) = \displaystyle \frac{F_3}{y_1y_2}&\quad
 \varphi(\bar F_3)= \displaystyle\frac{F_1}{y_2y_3}.
 \end{array}\]
\end{lemma}
\begin{proof}
 Using the definition and Lemma~\ref{lem phi}, we have
 \[\varphi(\bar F_1)=\varphi (1+\bar y_2+\bar y_2\bar y_3)=
 1+\frac{F_1}{y_3F_3}+\frac{F_1}{y_3F_3}\frac{F_2}{y_1F_1}=
 1+\frac{F_1}{y_3F_3}\left(1+\frac{F_2}{y_1F_1}\right)\] 
 Applying Lemma~\ref{frelations} to the term in parentheses, we see that this is equal to 
 \[ 1+\frac{F_1}{y_3F_3}\frac{(1+y_3)F_3}{y_1F_1} 
 =
\frac{y_1y_3+1+y_3}{y_1y_3} =\frac{F_2}{y_1y_3}. \]
The proofs of the other two identities are similar. 
\end{proof}

We also need the following result.
\begin{lemma}
 \label{lem stepcenter}
 We have $\varphi(1+\bar y_a)=(1+y_a) F_a/(y_{\rho(a)} F_{\rho(a)})$, for all $a=1,2,3$. That is
 \[
 \varphi(1+\bar y_1)=(1+y_1)\frac{ F_1}{y_{2} F_2 }
 \qquad
 \varphi(1+\bar y_2)=(1+y_2)\frac{ F_2}{y_{3} F_3 }
 \qquad
 \varphi(1+\bar y_3)=(1+y_3)\frac{ F_3}{y_{1} F_1 }.
  \]
\end{lemma}
\begin{proof}
We have
\[ \varphi(1+\bar y_1)= 1+\frac{F_3}{y_2F_2} 
= \frac{{y_2F_2}+F_3}{{y_2F_2}}
=(1+y_1)\frac{ F_1}{y_{2} F_2 }\]
where the first equation follows from Lemma~\ref{lem phi} and the last equation  from Lemma~\ref{frelations}. The proofs of the other two identities are similar.
\end{proof}

\subsubsection{The staircase lattices $\call$ and $\overline{\call}$}
We will prove the identity (\ref{eq goalphi}) using the combinatorial structure of the  lattices  $\call$ and $\overline{\call}$ of submodules of 
$T(\zs(i))$ and 
$\overline{T}(\zs(i))$. We actually work in a projection of these lattices by setting all variables $y_k$ with $k>9$ equal to 1.
This is sufficient, because the terms of the $F$-polynomial that do not contain $y_1,y_2,\ldots,y_9$ remain unchanged under the mutation sequence $\mu_0$.

\begin{figure}
\begin{center}
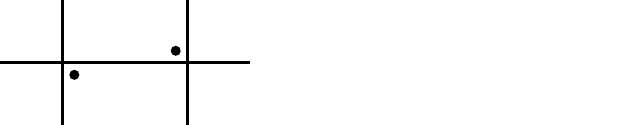
\caption{Covering relations in the poset of Kauffman states are given by transpositions along segments. The Kauffman state on the right is obtained from the one on the left by the transposition at segment $j$.}
\label{fig transposition}
\end{center}
\end{figure}

\begin{figure}
\begin{center}
\scriptsize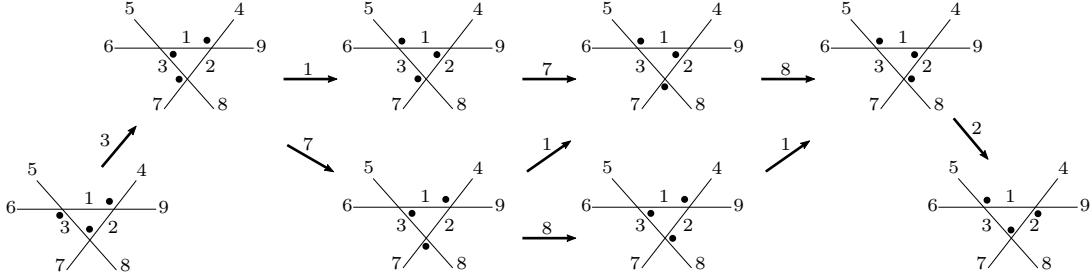
\caption{The sequence of transpositions at segments 3,1,7,8,2 in the lattice of Kauffman states of $K$, or, equivalently, in the lattice $\call$ of submodules of $T(\zs(i))$.}
\label{fig extranspositions}
\end{center}
\end{figure}

\begin{figure}
\begin{center}
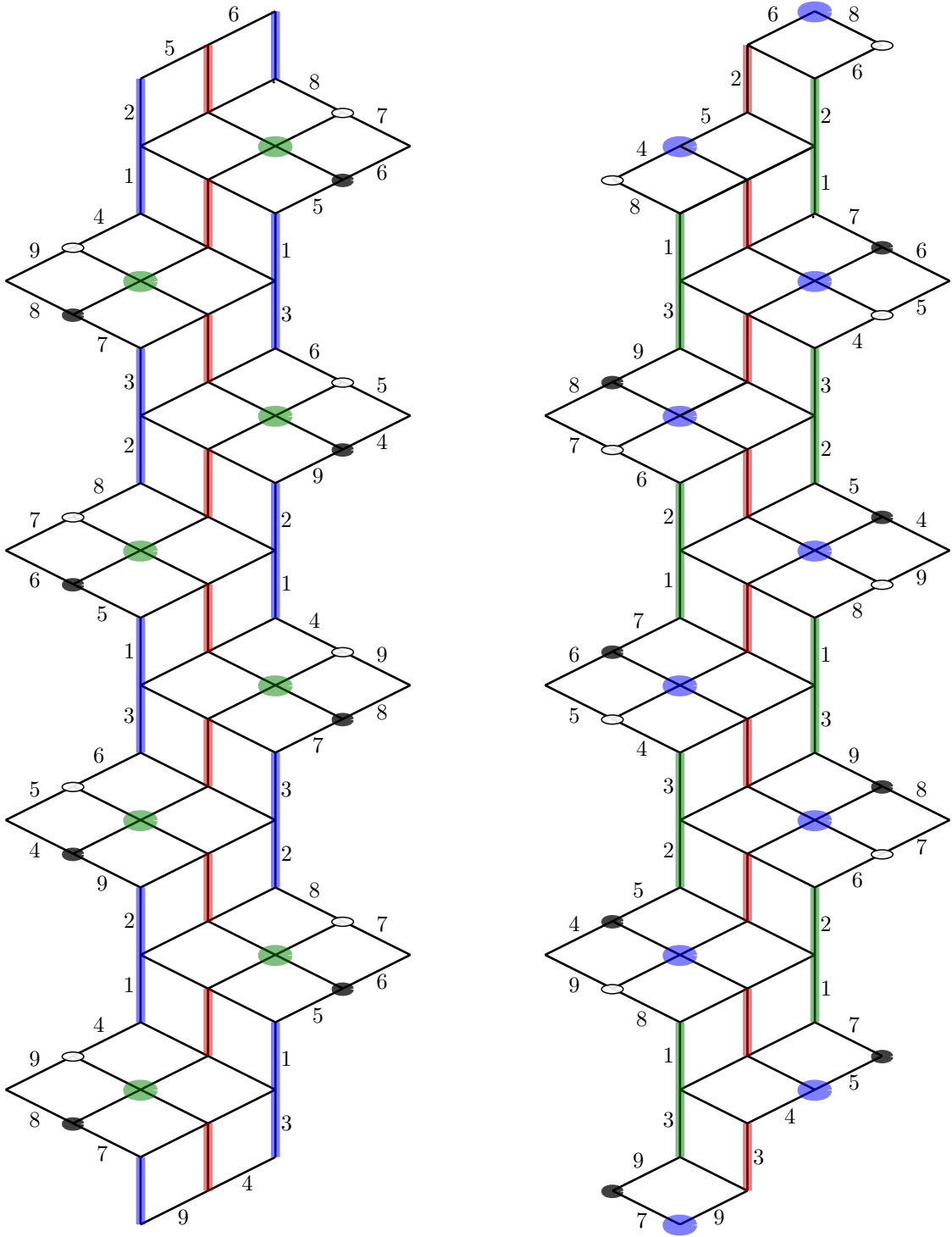
\caption{The staircase shapes of the lattice $\call$ (left) and $\overline{\call}$ (right). The blue columns on the left correspond to the blue landing centers on the right, and the green landing centers on the left to the green columns on the right. The red step centers are the same in both. The black lower boundary landing points on the left correspond to the black upper boundary landing points on the right, and the white upper boundary landing points on the left correspond to the white lower boundary landing points on the right.}
\label{fig staircase}
\end{center}
\end{figure}

By Theorem~\ref{thm lattice iso}, the submodule lattices $\overline{\call}$ and $\call$ are isomorphic to the lattices of Kauffman states of $\overline{K}$ and $K$ (with respect to the segment $i$). We may thus use the combinatorics of the Kauffman states to study the lattice structure. Recall that the covering relations for the Kauffman states are given by transpositions along segments as shown in Figure~\ref{fig transposition}.
For example, in the situation of the lattice $\call$ we have the  sequence of transpositions shown in Figure~\ref{fig extranspositions}. We may always assume that the transposition at the segments $4,5,\ldots,9$ is possible as soon as the markers at the three crossing points of the segments 1,2,3 are at the right place. Indeed, since we are working in the projection this assumption is without loss of generality. 
Therefore the lattice $\call$ has a staircase structure as illustrated in the left picture in Figure~\ref{fig staircase}. The right picture in the figure shows the lattice $\overline{\call}$ of Kauffman states of $\overline{K}$.

\begin{definition}\label{def column} We introduce the following terminology for subposets of the lattices $\call$ and $\overline{\call}$.
\begin{itemize}
\item [-]  An edge is \emph{vertical} if it is labeled by $1$, $2$ or $3$.
\item [-]
 A \emph{column} is a subposet given by two consecutive vertical edges. Columns are shaded  blue on the left in Figure~\ref{fig staircase} and green on the right.
\item [-] A \emph{step center} is a subposet consisting of exactly one vertical edge that is not incident to another vertical edge. Step centers are shaded  red on both sides in Figure~\ref{fig staircase}.
\item[-] A \emph{landing} is a maximal connected subposet with the property that none of its points is incident to a vertical edge. 
\item[-] The \emph{landing center} of  a landing is the unique point incident to four edges in the landing. Landing centers are shaded  green on the left in Figure~\ref{fig staircase} and blue on the right.
\item[-] The \emph{landing corner} is the unique point incident to exactly two edges in the landing. Landing corners are not colored in Figure~\ref{fig staircase}.
\item[-] The \emph{upper boundary point} and the \emph{lower boundary point} of a landing are the two points adjacent to the landing center. The upper boundary points are colored white and the lower boundary points are colored black on the left in Figure~\ref{fig staircase}. On the right side of the figure these colors are reversed.
\end{itemize}
\end{definition}

\begin{remark}
 \label{rem pi}
Our choice of  coloring should be understood by considering the projections $\pi $ and $\overline{\pi}$ of the two posets onto the plane. These projections are given by the contraction of all vertical edges. Notice that    $\pi(\call)=\overline\pi(\overline{\call})$, and both have exactly the same coloring. For example, each column  of $\call$ projects to the same point as a unique landing center of $\overline{\call}$. 

\end{remark}


\subsubsection{The morphism $\varphi$ on intervals of $\overline\call$}
\begin{figure}
\begin{center}
\scalebox{0.8}{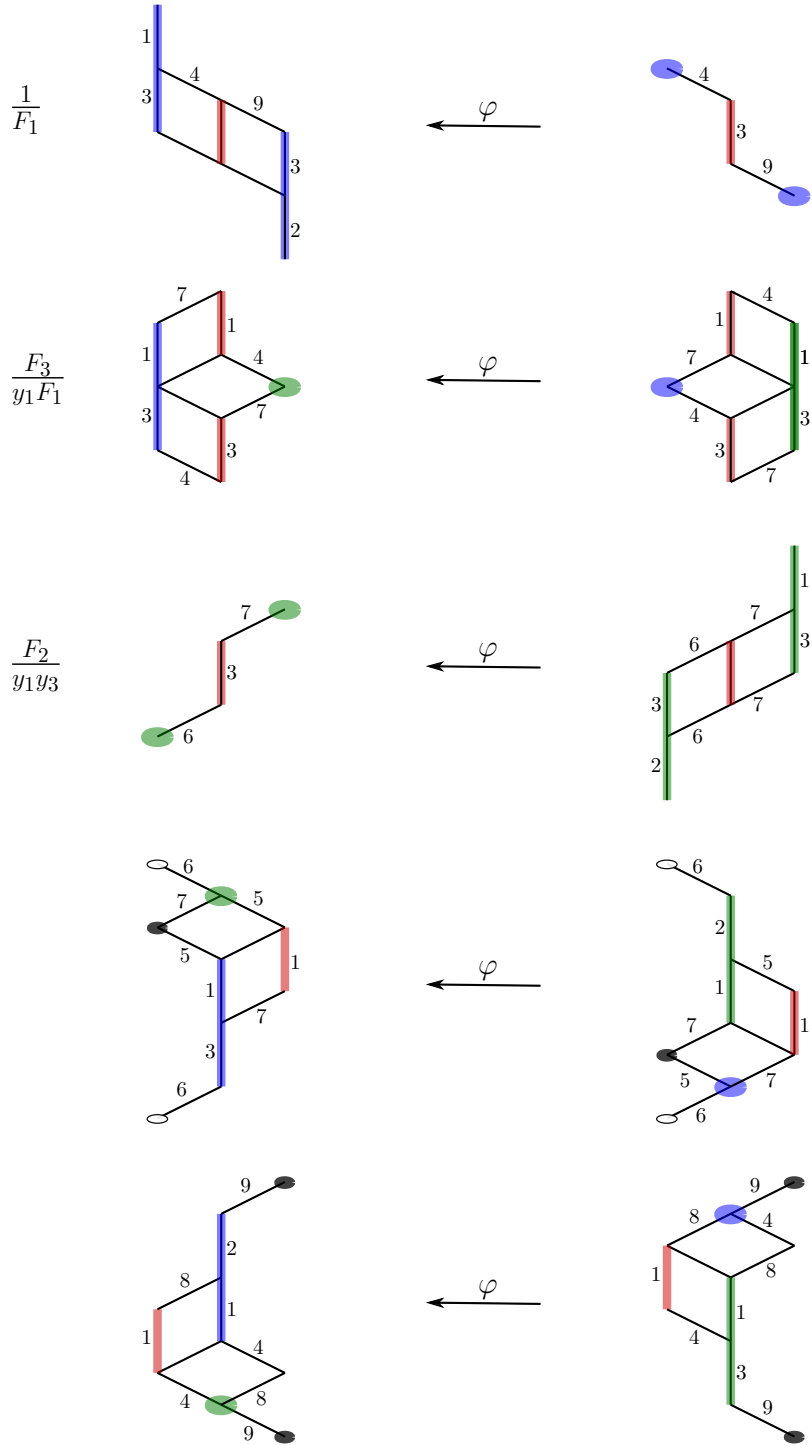}
\caption{Illustration of the five cases in Lemma \ref{lem intervals}. The pictures on the right hand side show the intervals $[\bar a,\bar b]$ in $\overline{\call}$ and the pictures on the left show their images under $\varphi$, which are intervals $[a,b]$ in $\call$.}
\label{fig lemintervals}
\end{center}
\end{figure}

We are now ready to compute $\varphi$ on the intervals of $\overline{\call}$ from the minimum of one subposet of a given type to the maximum of the next subposet of the same type, where the different types are landing centers,  step centers, columns, and upper and lower boundary landing points. 
Recall the projection maps $\pi,\overline{\pi}$ of Remark~\ref{rem pi}.
\begin{lemma} \label{lem intervals} 
Let $[\bar a,\bar b]$ be one of the intervals of $\overline{\call}$ specified in the cases (i)-(v) below. Then
\begin{equation}\label{eq intervals}\varphi\left(F\!\left([\bar a,\bar b]\right)\right)=Y(a) \,F\!\left([a,b]\right)\end{equation}
where $a, b\in\call$ are certain points specified below with the property that $\pi(a)=\overline{\pi}(\bar a)$, $\pi(b)=\overline{\pi}(\bar b)$,  $F[\bar a,\bar b], F[a,b]$ are the $F$-polynomials of the intervals, and  $Y(a)$ is a Laurent monomial in $y_1,y_2,y_3,F_1,F_2,F_3$.
%
\begin{itemize}
\item[ (i)] (Landing centers to columns) 
$\bar a<\bar b$ is a pair of consecutive landing centers in $\overline{\call}$,  $a$ is the unique minimal element of a column of $\call$ such that $\pi(a)=\overline\pi(\bar a)$ and $b$ is the unique maximal element of the following column of $\call$ going up in the poset. See the first row of Figure~\ref{fig lemintervals}. Then 
\[Y(a)=\frac{1}{F_{\rho(j)}},    \] 
where $j$ is the label of the unique vertical edge in $[\bar a,\bar b]$.
\item[(ii)] (Step centers to step centers) 
  $\bar a$ is the minimal element of a step center in $\overline{\call}$ and $\bar b$ is the maximal element of the following step center going up in the poset;
$a$ is the unique minimal element of the step center in $\call$ such that $\pi(a)=\overline\pi(\bar a)$ and $b$ is the unique maximal element of the following step center. See the second row  of Figure~\ref{fig lemintervals}. Then 
\[Y(a)=\frac{F_j}{y_{\rho(j)}F_{\rho(j)}}\]
where $j$ is the label of the lowest vertical edge in $[\bar a,\bar b]$.

\item[(iii)] (Columns to landing centers)
$\bar a$ is the minimal element of a column in $\overline{\call}$ and $\bar b$ is the maximal element of the following column going up in the poset;
$a<b$ are the unique landing centers in $\call$ such that
$\pi(a)=\overline\pi(\bar a)$ and $\pi(b)=\overline\pi(\bar b)$.
See third row  of Figure~\ref{fig lemintervals}.
 Then 
\[Y(a)=\frac{F_j}{y_{\rho(j)}y_{\rho^2(j)}}\]
where $j$ is the label of the lowest vertical edge in $[\bar a,\bar b]$.

\item[(iv)] (Lower boundary landing points to upper boundary landing points)
$\bar a<\bar b$ is a pair of consecutive lower boundary landing point of $\overline\call$; $a<b$ is the unique pair of consecutive upper boundary landing points of $\call$ such that $\pi(a)=\overline\pi(\bar a)$ and $\pi(b)=\overline\pi(\bar b)$.
See the fourth row  of Figure~\ref{fig lemintervals}.
Then
\[Y(a)=1.\]

\item[(v)] (Upper boundary landing points to lower boundary landing points)
$\bar a<\bar b$ is a pair of consecutive upper boundary landing point of $\overline\call$; $a<b$ is the unique pair of consecutive lower boundary landing points of $\call$ such that  $\pi(a)=\overline\pi(\bar a)$ and $\pi(b)=\overline\pi(\bar b)$.
See the fifth row of Figure~\ref{fig lemintervals}.
Then
\[Y(a)=1.\]
\end{itemize}
\end{lemma}

\begin{proof} In order to simplify the notation, we use the specific edge weights of the intervals in Figure~\ref{fig lemintervals}. Since the other intervals of the lattices are symmetric to the ones illustrated in that figure, the general proof is similar. 

(i) This case is illustrated in the first row of Figure~\ref{fig lemintervals}, where the interval $[\bar a, \bar b]$ is shown on the right and the interval $[a,b]$ on the left.  First note that the $F$-polynomial of the step center part of the interval $[\bar a,\bar b]$ equals $\bar y_9(1+\bar y_3)$ and applying $\varphi$ yields

\[\varphi(\bar y_9(1+\bar y_3))
=\frac{y_1y_2y_9}{F_3}\varphi(1+\bar y_3)
=\frac{y_1y_2y_9}{F_3}(1+ y_3)\frac{F_3}{y_1F_1} 
=\frac{1}{F_1}y_2y_9(1+ y_3) 
,\]
where the first identity follows from Lemma~\ref{lem phi} and the second from Lemma~\ref{lem stepcenter}. Clearly, the last expression is equal to the contribution of the step center to $F([a,b])$, shown in red in the left picture of the first row of Figure~\ref{fig lemintervals}.

Now consider the two contribution of the landing centers to $\varphi(F([\bar a,\bar b]))$. 
 \[\varphi(1+ \bar y_3\bar y_4 \bar y_9) 
 =
 1+\frac{F_2}{y_1F_1} y_4F_3 \frac{y_1y_2y_9}{F_3}
=
 1+y_2y_4y_9\frac{F_2}{F_1} 
=
\frac{1}{F_1} ( F_1+y_2y_4y_9F_2),
\]
where the first identity follows from Lemma~\ref{lem phi}. Using the definition of $F_1$ and $F_2$, we see that the last expression is equal to $\frac{1}{F_1} ( (1+y_2+y_2y_3)+y_2y_4y_9(1+y_3+y_1y_3))$, which is identical to  the part of $\frac{1}{F_1} F([a,b])$ contributed by the two columns in $[a,b]$. This shows (i).

(ii) This case is illustrated in the second row of Figure~\ref{fig lemintervals}. First consider the two step centers in $[\bar a,\bar b]$. Their contribution to $\varphi([\bar a,\bar b])$ is
 \begin{eqnarray*}
\varphi ((1+\bar y_3) +\bar y_3\bar y_4\bar y_7(1+\bar y_1))
&=&\varphi(1+\bar y_3)+\frac{F_2}{y_1F_1}y_4F_3 \frac{y_2y_3y_7}{F_1} \varphi(1+\bar y_1)
\\
&=& (1+y_3)\frac{F_3}{y_1F_1}  
+y_3y_4y_7 \frac{y_2 F_2F_3}{y_1F_1F_1} (1+y_1)\frac{F_1}{y_2F_2}
\\
&=& \frac{F_3}{y_1F_1}\left( (1+y_3) 
+y_3y_4y_7  (1+y_1)\right)\end{eqnarray*}
where the first equation follows from Lemma~\ref{lem phi}, the second from Lemma~\ref{lem stepcenter} and the third is obtained by algebraic simplification. The last expression  is equal to the contribution of the two step centers to $\frac{F_3}{y_1F_1}F([a,b])$.

Now consider the column part of $\varphi(F([\bar a,\bar b]))$.
\[\varphi(\bar y_7(1+\bar y_3+\bar y_3\bar y_1))
=\varphi(\bar y_7)\varphi(\bar F_2)
=\frac{y_2y_3y_7}{F_1}\frac{F_3}{y_1y_2} 
=\frac{F_3}{y_1F_1} (y_3y_7),
\]
where the second identity uses Lemmas \ref{lem phi} and \ref{lem phif}. The last expression is equal to the contribution of the landing center of $[a,b]$ to $\frac{F_3}{y_1F_1}F([a,b])$.

It remains to consider the contribution of the landing center of $[\bar a,\bar b]$. It is equal to 
\[\varphi(\bar y_3\bar y_4) =
\frac{F_2}{y_1F_1} y_4F_3
=\frac{F_3}{y_1F_1} y_4F_2
=\frac{F_3}{y_1F_1} y_4(1+y_3+y_1y_3),\]
where the first identity uses Lemma \ref{lem phi}.  The last expression is equal to the contribution of the column of $[a,b]$ to $\frac{F_3}{y_1F_1}F([a,b])$.
This completes the proof of part (ii).

(iii) This case is illustrated in the third row of Figure~\ref{fig lemintervals}. The contribution of the step center of $[\bar a,\bar b]$ equals
\[\varphi(\bar y_2\bar y_6(1+\bar y_3))
=\frac{F_1}{y_3F_3} y_6F_2 \,\varphi(1+\bar y_3)
=\frac{F_1}{y_3F_3} y_6F_2\, (1+y_3)\frac{F_3}{y_1F_1}
=\frac{F_2}{y_1y_3} y_6(1+y_3),\]
where the first equation follows from Lemma~\ref{lem phi}, the second from Lemma~\ref{lem stepcenter} and the third is obtained by algebraic simplification. The last expression is equal to the contribution of the step center of $[a,b]$ to $\frac{F_2}{y_1y_3}F([a,b])$.

Next we consider the two columns of $[\bar a,\bar b]$. 
The contribution of the first column equals
\[
\varphi(1+\bar y_2+\bar y_2\bar y_3)
=\varphi(\bar F_1)
=\frac{F_2}{y_1y_3}\] 
by Lemma \ref{lem phif}.  This expression is equal to the contribution of the first landing center of $[a,b]$ to $\frac{F_2}{y_1y_3}F([a,b])$.

Similarly, the contribution of the second column equals
\[
\varphi\big(\bar y_2\bar y_6\bar y_7(1+\bar y_3+\bar y_1\bar y_3)\big)
=\varphi(\bar y_2\bar y_6\bar y_7) \, \varphi(\bar F_2)
=\frac{F_1}{y_3F_3} y_6F_2 \frac{y_2y_3y_7}{F_1} \frac{F_3}{y_1y_2}
=\frac{F_2}{y_1y_3} y_3y_6y_7,
\]
where the second identity uses Lemma \ref{lem phi} and Lemma~\ref{lem phif}.  The last expression is equal to the contribution of the second landing center of $[a,b]$ to $\frac{F_2}{y_1y_3}F([a,b])$. This completes the proof of (iii).

(iv) This case is illustrated in the fourth row of Figure~\ref{fig lemintervals}. First note that the posets shown there are in fact intervals. Indeed, observe that in the right picture, the corner point of the landing is missing. It would be obtained from the minimal element by going along an edge with weight 5. This corresponds to the transposition at 5 in the corresponding Kauffman state. However, since the transposition at 6 is possible at the minimal element, the state marker of the region formed by the segments $6,2,5$ in $\overline{K}$ must be at the crossing point of 2 and 6, see the top right picture in Figure~\ref{fig R3}. Therefore the transposition at segment 5 is not possible at the minimal element. Only once we have performed the transposition at 6, the transposition at 5 becomes possible. A similar argument shows that all the posets shown in rows four and five of Figure~\ref{fig lemintervals} are intervals.

Now let's go back to our computation of $\varphi$. The contribution of the step center of $[\bar a,\bar b]$ is
\[ \varphi (\bar y_6\bar y_7(1+\bar y_1))=
y_6 F_2\frac{y_2y_3y_7}{F_1}\,\varphi(1+\bar y_1)
=y_6 F_2\frac{y_2y_3y_7}{F_1}\,(1+ y_1) \frac{F_1}{y_2 F_2}
=y_3y_6y_7 
(1+ y_1),
\]
where the first identity uses Lemma \ref{lem phi} and the second uses Lemma~\ref{lem stepcenter}. The last expression is the contribution of the step center of $[a,b]$ to $F([a,b]).$ 

The contribution of the landing center of $[\bar a,\bar b]$ is
\[\varphi(\bar y_6)=y_6 F_2 =y_6(1+y_3+y_1y_3)
\]
which is equal to the contribution of the column of $[a,b]$  to $F([a,b]).$ 

The contribution of the column of $[\bar a,\bar b]$ is
\[\varphi(\bar y_5\bar y_6\bar y_7 \bar F_3)
=\frac{y_1y_3y_5}{F_2} y_6F_2\frac{y_2y_3y_7}{F_1} \frac{F_1}{y_2y_3}
=y_1y_3y_5 y_6y_7,
\]
where the first identity uses Lemma \ref{lem phi} and Lemma~\ref{lem phif}.  The last expression is equal to the contribution of the landing center of $[a,b]$ to $F([a,b]).$ 

It remains to consider the contribution of the the three boundary landing points in $[\bar a,\bar b]$. It is equal to
\begin{eqnarray*}\varphi(1+\bar y_5\bar y_6+
\bar y_1\bar y_2\bar y_5\bar y_6^2\bar y_7)
&=&1+ \frac{y_1 y_3 y_5}{F_2} y_6 F_2 
+ \frac{F_3}{y_2F_2} \frac{F_1}{y_3F_3} \frac{y_1y_3y_5}{F_2} \\
&=&
1+y_1y_3y_5y_6+y_1y_3y_5y_6^2y_7,
\end{eqnarray*}
where the first identity uses Lemma \ref{lem phi}.  The last expression is equal to the contribution of the three boundary landing points of $[a,b]$ to $F([a,b]).$ 
This completes the proof of part (iv), and the proof of part (v) is similar.
\end{proof}

The next result  follows from the proof.
\begin{corollary}\label{cor intervals}
If $\bar a$ is the minimal element of any of the intervals in Lemma~\ref{lem intervals}, then up to a Laurent monomial in $y_1,y_2,y_3,F_1,F_2,F_3$, the morphism $\varphi$ maps landing centers to columns, step centers to step centers, columns to landing centers, lower boundary landing points to upper boundary landing points and upper boundary landing points to lower boundary landing points.
\end{corollary}

Our next goal is to show that $\varphi$ is compatible with unions of intervals in Proposition~\ref{prop gluing} below. We need a preparatory lemma that specifies the relation between the Laurent monomials $Y(a)$. We use the notation $F(a)$ for the $F$-polynomial of the subposet of all $x$ with $x\le a$.
\begin{lemma}
 \label{lem gluing}
 Let $[\bar a,\bar b],[a,b]$ be one of the pairs of intervals of Lemma~\ref{lem intervals} and $[\bar c,\bar d],[c,d]$ be another such pair and assume $\bar a<\bar c$. Then
    \[ Y(a) \, \frac{\varphi(F(\bar a))}{F(a)}=Y(c)\, \frac{\varphi(F(\bar c))}{F(c)}. \] 
    In particular, if $\bar a_0$ is the minimal element of $\overline{\call}$, then \[ Y(a_0) =Y(c)\, \frac{\varphi(F(\bar c))}{F(c)} .\] 
\end{lemma}
\begin{proof}
The second equation follows from the first one because if $\bar a $ is the minimal element of $\overline{\call}$ then $a$ is the minimal element of $\call$ and thus $F(\bar a)=F(a)=1$. Therefore it suffices to show the first equation.

Throughout the proof we call an element $\bar x$ of the poset $\overline{\call}$ an \emph{interval point} if $\bar x$ is one of the following: 
 the minimal element of a column,
 the minimal element of a step center,
 a landing center,
 a lower boundary landing point, or
 an upper boundary landing point.

Thus in the notation of the lemma, the points $\bar a, \bar c$ are interval points and $\bar a<\bar c$. We shall prove the statement by induction on the number of interval points that lie strictly between $\bar a$ and $\bar c$.

Suppose first that this number is zero. So whenever $\bar x$ is an interval point such that $\bar a\le \bar x \le\bar c$ then $\bar x=\bar a$ or $\bar x=\bar c$. We need to distinguish five cases.

1. Assume $\bar a$ is  a landing center. Without loss of generality, we may assume that the interval $[\bar a,\bar b]$ is the one shown on the right in the first row of Figure~\ref{fig lemintervals}. Then  $\bar c$ is adjacent to $\bar a$ in $\overline{\call}$, and $\bar c$ is either the minimal element of a step center  or an upper boundary landing point. 

In the former case,  the interval $[\bar c,\bar d]$
is seen in the second row  of Figure~\ref{fig lemintervals}. 
Using parts (i) and (ii) of Lemma~\ref{lem intervals}, we obtain
\[\frac{Y(a)}{Y(c)}=\frac {\left(\frac{1}{F_1} \right)} {\left(\frac{F_3}{y_1F_1}\right)} =\frac{y_1}{F_3}.\]
On the other hand, using Figure~\ref{fig staircase} and Lemma~\ref{lem phi}
\[\frac{F(a)}{F(c)}\ \varphi\!\left(\frac{F(\bar c)}{F(\bar a)}\right) 
=\frac{1}{y_2y_9}\varphi(\bar y_9) =\frac{y_1}{F_3}\]
as desired.

If $\bar c$ is the upper boundary landing point, the interval $[\overline{c},\overline{d}]$ is illustrated in the fifth row of Figure~\ref{fig lemintervals}. We have
\[\frac{Y(a)}{Y(c)}=\frac {\left(\frac{1}{F_1} \right)} {1} =\frac{1}{F_1},\]
 whereas
\[\frac{F(a)}{F(c)}\ \varphi\!\left(\frac{F(\bar c)}{F(\bar a)}\right) 
=\frac{1}{y_2y_3y_7}\varphi(\bar y_7) =\frac{1}{F_1}\]
as desired.

2. Assume $\bar a$ is the minimal element of a step center,
 Without loss of generality, we may assume that the interval $[\bar a,\bar b]$ is the one shown on the right in the second row of Figure~\ref{fig lemintervals}. Then  $\bar c$ is either the landing center or the minimal element of the column in $[\bar a,\bar b]$.
In the former case, parts (ii) and (i) of Lemma~\ref{lem intervals} 
yield
\[\frac{Y(a)}{Y(c)}=\frac  {\left(\frac{F_3}{y_1F_1}\right)}{\left(\frac{1}{F_2} \right)} =\frac{F_2F_3}{y_1F_1},\]
which is equal to 
\[\frac{F(a)}{F(c)}\ \varphi\!\left(\frac{F(\bar c)}{F(\bar a)}\right) 
=\frac{1}{y_4}\varphi(\bar y_3\bar y_4) =\frac{1}{y_4}\frac{F_2}{y_1F_1} y_4F_3.\]
In the other case,  $\bar c$ is  the minimal element of the column, and
\[\frac{Y(a)}{Y(c)}=
\frac  {\left(\frac{F_3}{y_1F_1}\right)}{\left(\frac{F_3}{y_1y_2} \right)} =\frac{y_2}{F_1},\]
which is equal to 
\[\frac{F(a)}{F(c)}\ \varphi\!\left(\frac{F(\bar c)}{F(\bar a)}\right) 
=\frac{1}{y_3y_7}\varphi(\bar y_7) 
=\frac{1}{y_3y_7}\frac{y_2y_3y_7}{F_1}.\]

3. Assume $\bar a$ is the minimal element of a column. Without loss of generality, we may assume that the interval $[\bar a,\bar b]$ is the one shown on the right in the third row of Figure~\ref{fig lemintervals}. Then  $\bar c$ is either the minimal element of the step center in $[\bar a,\bar b]$ or the lower boundary point reached from $\bar a$ following the edges labeled 2,3,4.
In the former case, parts (iii) and (ii) of Lemma~\ref{lem intervals} 
yield
\[\frac{Y(a)}{Y(c)}=\frac  {\left(\frac{F_2}{y_1y_3}\right)}{\left(\frac{F_3}{y_1F_1} \right)} 
=\frac{F_1F_2}{y_3F_3},\]
which is equal to 
\[\frac{F(a)}{F(c)}\ \varphi\!\left(\frac{F(\bar c)}{F(\bar a)}\right) 
=\frac{1}{y_6}\varphi(\bar y_2\bar y_6) =\frac{1}{y_6}\frac{F_1}{y_3F_3} y_6F_2.\]
In the other case, $\bar c$ is a lower boundary point and  parts (iii) and (iv) of Lemma~\ref{lem intervals} 
yield
\[\frac{Y(a)}{Y(c)}=\frac  {\left(\frac{F_2}{y_1y_3}\right)}{1} 
=\frac{F_2}{y_1y_3},\]
which is equal to 
\[\frac{F(a)}{F(c)}\ \varphi\!\left(\frac{F(\bar c)}{F(\bar a)}\right) 
=\frac{1}{y_4}\varphi(\bar y_2\bar y_3\bar y_4) =\frac{1}{y_4}\frac{F_1}{y_3F_3}\frac{F_2}{y_1F_1}  y_4F_3.\]

4. Assume $\bar a$ is a lower boundary landing point. Without loss of generality, we may assume that the interval $[\bar a,\bar b]$ is the one shown on the right in the fourth row of Figure~\ref{fig lemintervals}. Then  $\bar c$ is the adjacent landing center. Parts (iv) and (i) of Lemma~\ref{lem intervals} 
yield
\[\frac{Y(a)}{Y(c)}=\frac  {1}{\left(\frac{1}{F_2} \right)} 
=F_2,\]
which is equal to 
\[\frac{F(a)}{F(c)}\ \varphi\!\left(\frac{F(\bar c)}{F(\bar a)}\right) 
=\frac{1}{y_6}\varphi(\bar y_6) =\frac{1}{y_6} y_6F_2.\]

5. Assume $\bar a$ is an upper boundary landing point.
 Without loss of generality, we may assume that the interval $[\bar a,\bar b]$ is the one shown on the right in the fifth row of Figure~\ref{fig lemintervals}. Then  $\bar c$ is the minimal element of the adjacent column. Parts (v) and (iii) of Lemma~\ref{lem intervals} 
yield
\[\frac{Y(a)}{Y(c)}=\frac  {1}{\left(\frac{F_3}{y_1y_2} \right)} 
=\frac{y_1y_2} {F_3}\]
which is equal to 
\[\frac{F(a)}{F(c)}\ \varphi\!\left(\frac{F(\bar c)}{F(\bar a)}\right) 
=\frac{1}{y_9}\varphi(\bar y_9) =\frac{1}{y_9} \frac{y_1y_2y_9}{F_3}.\]
This completes the proof of the base case of our induction. 

For the induction step, assume that there is at least one interval point $\bar x$ such that  $\bar a<\bar x<\bar c$. Let $[\bar x,\bar y]$ and $[x,y]$ be the corresponding intervals in $\overline{\call}$ and $\call$, respectively.
By induction, we know 
\[\frac{  Y(a)}{Y(x)} = \frac{F(a)}{F(x)}\,\varphi\!\left(\frac{F(\bar x)}{F(\bar a)}\right)
\quad\textup{and}\quad
\frac{  Y(x)}{Y(c)} = \frac{F(x)}{F(c)}\,\varphi\!\left(\frac{F(\bar c)}{F(\bar x)}\right).
\]
Taking the product of these two equations we get
\[\frac{  Y(a)}{Y(c)} = \frac{F(a)}{F(c)}\,\varphi\!\left(\frac{F(\bar c)}{F(\bar a)}\right),
\]
and the proof is complete.
\end{proof}

We are now ready to prove the compatibility of $\varphi$ with unions of intervals. 
\begin{prop}\label{prop gluing} Let $[\bar a,\bar b],[a,b]$ be one of the pairs of intervals of Lemma~\ref{lem intervals} and $[\bar c,\bar d],[c,d]$ be another such pair. Denote by $\bar a_0$ the minimal element of the lattice $\overline{\call}$. Then
   \[\varphi\big(F(\bar a)\,F\!\left([\bar a, \bar b]\cup [\bar c, \bar d]\right)\big) =Y(a_0)
\, F(a)\, F( [a,b]\cup  [c,d] ) .\] 
In particular,
\begin{equation}
\label{eq prop 413}
\varphi\big( F(\overline{\call})\big)=Y(a_0)F(\call).
\end{equation}
\end{prop}
\begin{proof}
 By the inclusion-exclusion principle, the left hand side of the equation is equal to 
\begin{equation}\label{eq 410} \varphi\big(F( \bar a)\,F[\bar a,\bar b] +F(\bar c)\,F[\bar c,\bar d] -\bar C \big),\end{equation}
 where $\bar C$ is the sum of the terms that come from the intersection $[\bar a,\bar b]\cap [\bar c,\bar d]$. By Theorem~\ref{thm grin} we know that all terms in the $F$-polynomials have coefficient 1. Therefore $\bar C$ is precisely the sum of all the terms that appear in both $F( \bar a)\,F[\bar a,\bar b]$ and $F(\bar c)\,F[\bar c,\bar d]$.
 Because of Lemma~\ref{lem intervals},  the expression in (\ref{eq 410}) is equal to 
\begin{equation}
\label{eq 410b} 
\varphi(F(\bar a))\, Y(a)\,F[a,b] +\varphi(F(\bar c)) \,Y(c)\,F[c,d] -C,
 \end{equation}
where $C$ is the sum of all terms that appear in both of the first two expressions in (\ref{eq 410b}). From Lemma~\ref{lem gluing}, we have $\varphi(F(\bar a)) \,Y(a) =Y(c)\frac{\varphi(F(\bar c) )}{F(c)}F(a)$. Thus (\ref{eq 410b}) equals
\begin{equation}\label{eq 410c}
\textstyle Y(c)\frac{\varphi(F(\bar c) )}{F(c)} 
\big( F(a) F[a,b] + F(c) F[c,d]  -C'
\big),
\end{equation}
where $C'$ is the sum of all terms that appear in both 
$ F(a) F[a,b]$ and $ F(c) F[c,d] $. 
Again, since all terms in the $F$-polynomials have coefficient 1, the expression $C'$ is the sum of the terms that come from the intersection $[a,b]\cap[c,d]$. Moreover, the second equation of Lemma~\ref{lem gluing} shows that $Y(c)\frac{\varphi(F(\bar c) )}{F(c)} = Y(a_0)$.
 Thus (\ref{eq 410c}) equals
\begin{equation}\label{eq 410d}
Y(a_0)
\,F(a) \,F( [a,b]\cup  [c,d] ). \qedhere
\end{equation} 
\end{proof}

\subsubsection{Computation of the $g$-vector and proof of Theorem~\ref{thm module main}}
 Recall that in order to prove the theorem  we need to show equation  (\ref{eq goalphi})
\begin{equation*}
\varphi\left(F_{\overline{T}(\zs(i))}\right) 
   \ \frac{\prod_{j =1}^3  (F_{j})^{\overline{g}_j} }
{ F_{i;t_2}^{t_0}(y_{j;t_1}) \vert_{\mathbb{P}}}
= F_{T(\zs(i))}.
\end{equation*}
In view of equation~(\ref{eq prop 413}) it only remains to show that 
\begin{equation}
\label{eq 424}
\frac{Y(a_0)}{ F_{i;t_2}^{t_0}(y_{j;t_1}) \vert_{\mathbb{P}}}
= \frac{1}{ \prod_{j=1}^3 (F_{j})^{\overline{g}_j}}.
\end{equation}
Recall that the  division by $ { F_{i;t_2}^{t_0}(y_{j;t_1}) \vert_{\mathbb{P}}}$ corresponds to clearing the denominators in the $F$-polynomial; thus it removes all negative powers of $y_k$ that may occur in $Y(a_o)$.

\begin{figure}
\begin{center}
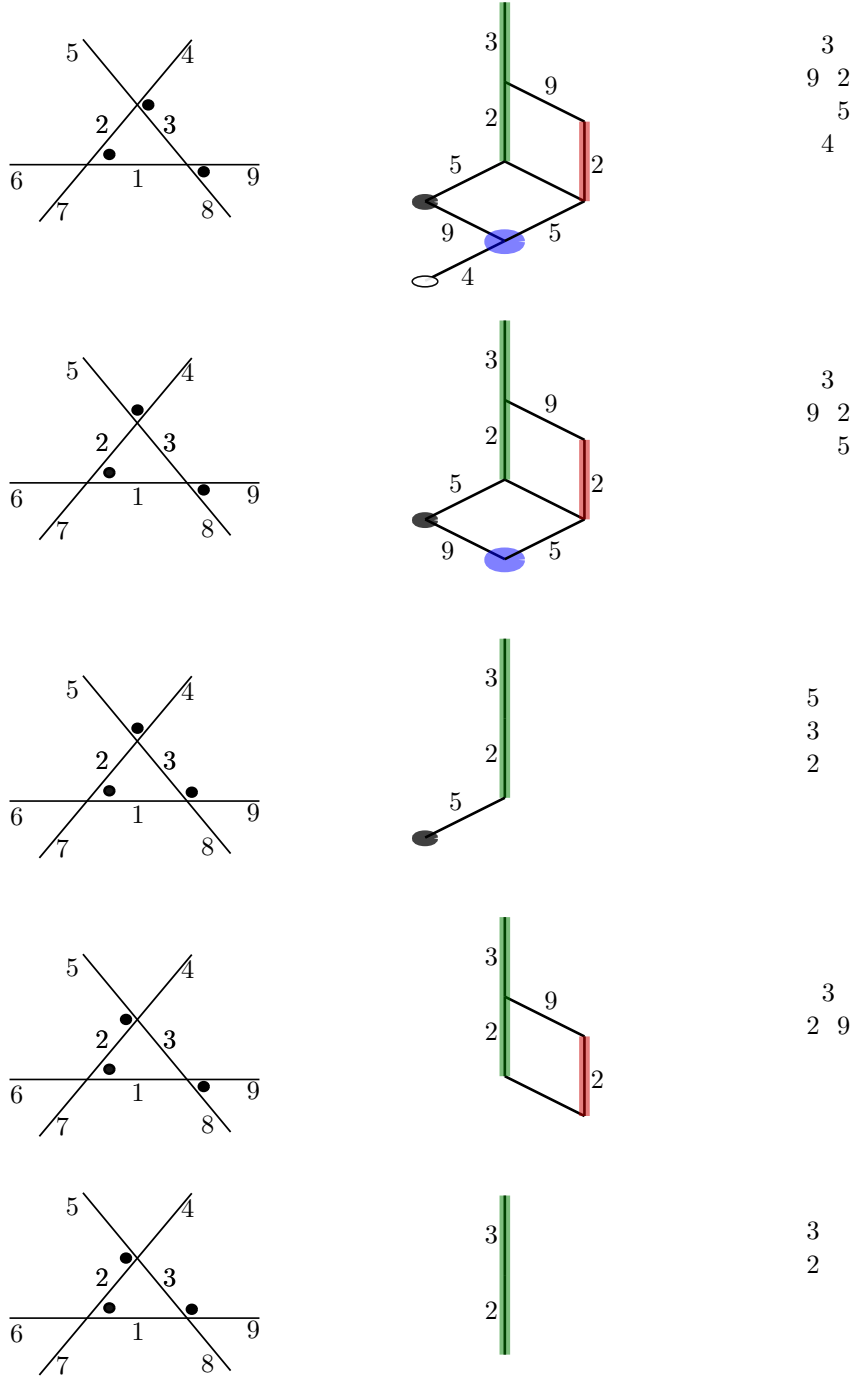
\caption{The possible minimal Kauffman states on $\overline{K}$ (left), the corresponding local configurations in the poset $\overline{\call}$ (center) and the corresponding submodule of $\overline{T}(i)$ (right).}
\label{fig minimal state}
\end{center}
\end{figure}

In the above equation, $a_0$ corresponds to the minimal element $\bar a_0$ of the lattice $\overline{\call}$, thus it corresponds to the minimal Kauffman state of $\overline{K}$ with respect to the segment $\zs(i)$. The different possibilities for this state are shown in the left column of Figure~\ref{fig minimal state}. Since $\zs(i)\ne 1,2,3$ there has to be exactly one marker in the triangular region formed by segments $1,2,3$. By symmetry, we may assume without loss of generality that this marker sits at the common endpoint of 1 and 2.

 Next, consider the crossing point where the segments $2,3,4,5$ meet. Its marker cannot sit between the segments 2 and 3, because that would be a second marker in the same region.  Each of the other three positions is allowed.
 
 Finally, consider the third crossing point, the one where the segments $1,3,9,8$ meet. As before, its marker cannot sit between the segments 1 and 3, because that would be a second marker in the same region. The marker also cannot sit between segments 1 and 8, because that would mean that the transposition at segment 1 would connect our state to a smaller state, which contradicts the minimality of our state. Moreover, if our previous marker sits between segments 3 and 4, then we also cannot put our current marker between 3 and 9, because that would produce two markers in the same region. This gives a total of five possibilities, see the left column of Figure~\ref{fig minimal state}.

The middle column  of the figure shows the local configuration of the poset $\overline{\call}$ above $\bar a_0$ in each case, and the right column shows the corresponding submodule of $\overline{T}(i)$.
In the first row, the minimal state only admits the transposition at 4. Once 4 is done, the transpositions at 5 and 9 are possible, and after 5, the transposition at 2 is possible. We see that $\bar a_0$ is a lower  boundary landing point, and thus Lemma~\ref{lem intervals} implies $Y(a_0)=1$. On the other hand, this shows that the corresponding module $\overline{T}(\zs(i))$  has a simple socle $S(4)$ and is locally of the form 
$\begin{smallmatrix} 
3\\9\ 2\\\ \ 5\\4
\end{smallmatrix}$. Its minimal injective presentation is of the form
\[\xymatrix{ 0\ar[r]&\overline{T}(\zs(i))\ar[r]&I(4)\oplus I' \ar[r]& I''},\]
for some injectives $I',I''$ that do not have direct summands $I(j)$ with $j\le 9$. Thus the $g$-vector is  $\overline{g}_j=-e_4-e'+e'' $, where $e',e''$ correspond to $I',I''$. In particular, the right hand side of equation~(\ref{eq 424}) is equal to $1$. Thus the equation holds in this case.

In the second row, the minimal state  admits the transpositions at 5 and at 9. Once 5 is done, the transpositions at 2 is possible, and after 5, 9 and 2, the transposition at 3 is possible. We see that $\bar a_0$ is a landing center, and thus Lemma~\ref{lem intervals} implies $Y(a_0)=F_3^{-1}$. On the other hand, the corresponding module $\overline{T}(\zs(i))$ has a socle $S(5)\oplus S(9)$ and is locally of the form 
$\begin{smallmatrix} 
3\\9\ 2\\\ \ \ 5
\end{smallmatrix}$. Its minimal injective presentation is of the form
\[\xymatrix{ 0\ar[r]&\overline{T}(\zs(i))\ar[r]&I(5)\oplus I(9)\oplus I' \ar[r]& I(3)\oplus I''},\]
where $I',I''$ do not have direct summands $I(j)$ with $j\le 9$. Thus the $g$-vector is  $\overline{g}_j=e_3-e_5-e_9-e'+e'' $, where $e',e''$ correspond to $I',I''$. In particular, the right hand side of equation~(\ref{eq 424}) is equal to $F_3^{-1}$, and we are done.

In the third row, the minimal state  admits the transposition at 5 only. Once 5 is done, the transpositions at 2 is possible, and after that the transposition at 3. We see that $\bar a_0$ is an upper landing boundary point, and thus Lemma~\ref{lem intervals} implies $Y(a_0)=1$. On the other hand, the corresponding module $\overline{T}(\zs(i))$ has a simple socle $S(5)$ and is locally of the form 
$\begin{smallmatrix} 
3\\ 2\\5
\end{smallmatrix}$. Its minimal injective presentation is of the form
\[\xymatrix{ 0\ar[r]&\overline{T}(\zs(i))\ar[r]&I(5)\oplus I' \ar[r]&  I''},\]
where $I',I''$ do not have direct summands $I(j)$ with $j\le 9$. Thus the $g$-vector is  $\overline{g}_j=-e_5-e'+e'' $, where $e',e''$ correspond to $I',I''$. In particular, the right hand side of equation~(\ref{eq 424}) is equal to $1$, and we are done.

In the fourth row  of Figure~\ref{fig minimal state}, the minimal state  admits the transpositions at 2 and at 9. Once 2 and 9 are done, the transpositions at 3 is possible. We see that $\bar a_0$ is the minimal element of a step center, and thus Lemma~\ref{lem intervals} implies $Y(a_0)=F_2y_3^{-1}F_3^{-1}$, and thus the left hand side of equation~(\ref{eq 424}) equals $F_2F_3^{-1}$. On the other hand, the corresponding module $\overline{T}(\zs(i))$ has  socle $S(2)\oplus S(9)$ and is locally of the form 
$\begin{smallmatrix} 
3\\9\ 2
\end{smallmatrix}$. Its minimal injective presentation is of the form
\[\xymatrix{ 0\ar[r]&\overline{T}(\zs(i))\ar[r]&I(2)\oplus I(9)\oplus I' \ar[r]& I(3)\oplus I''},\]
where $I',I''$ do not have direct summands $I(j)$ with $j\le 9$. Thus the $g$-vector is  $\overline{g}_j=e_3-e_2-e_9-e'+e'' $, where $e',e''$ correspond to $I',I''$. In particular, the right hand side of equation~(\ref{eq 424}) is equal to $F_2F_3^{-1}$, and we are done.

Lastly, in the fifth row, the minimal state  admits the transpositions at 2 only. Once 2 is done, the transpositions at 3 is possible. We see that $\bar a_0$ is the minimal element of a column, and thus Lemma~\ref{lem intervals} implies $Y(a_0)=F_2y_1^{-1}y_3^{-1}$ 
and thus the left hand side of equation~(\ref{eq 424}) equals $F_2$. On the other hand, the corresponding module $\overline{T}(\zs(i))$ has a simple socle $S(2)$ and is locally of the form 
$\begin{smallmatrix} 
3\\ 2
\end{smallmatrix}$. Its minimal injective presentation is of the form
\[\xymatrix{ 0\ar[r]&\overline{T}(\zs(i))\ar[r]&I(2)\oplus I' \ar[r]&\oplus I''},\]
where $I',I''$ do not have direct summands $I(j)$ with $j\le 9$. Thus the $g$-vector is  $\overline{g}_j=-e_2-e'+e'' $, where $e',e''$ correspond to $I',I''$. In particular, the right hand side of equation~(\ref{eq 424}) is equal to $F_2$, and we are done.

This completes the proof of Theorem~\ref{thm module main} in the case where $K$ has no bigon and $\zs(i)\ne 1,2,3$.

\medskip
{\bf Second Case.} Assume now that $\zs(i)=1,2, $ or $3$.
As before, let $\mu$ be the mutation sequence constructed in (\ref{eq mutation sequence}). Recall that $\mu$ transforms the initial cluster $\xx_{t_0}$ into the cluster $\xx_{t}.$ Denote by $\zD,\zD'$  the two triangular regions of the first two R3 moves in $\mu$, and let $1,2,3$ denote the segments of $\zD$ and  $1',2',3'$ those of $\zD'$. Thus the sequence of mutations is of the form 
$\mu=\mu_{R(\zD)}\mu_{R(\zD')} \nu$, where $\mu_{R(\zD)}$ mutates at vertices 1,2,3,1,2,3,2,3,2,  and $\mu_{R(\zD')}$ mutates at vertices  $ 1',2',3',1',2',3',2',3',2'$.
By admissibility of the sequence $\mu$ (see Definition~\ref{def admissible}), the regions $\zD,\zD'$ in $K$ are disjoint. Therefore, in the quiver $Q$, there is no arrow connecting a vertex in $\{1,2,3\}$ to a vertex in $\{1',2',3'\}$. Consequently the sequences $\mu_{R(\zD)},\mu_{R(\zD')}$ commute.

Thus we have a second mutation sequence $\mu'=\mu_{R(\zD')}\mu_{R(\zD)} \nu$, obtained from $\mu$ by switching the first two R3 sequences, that also transforms the initial cluster into the cluster $\xx_{t}$. Note that the sequence $\mu'$ may not be admissible, but this does not matter in what follows.

We have shown above, using the mutation sequence $\mu$, that every variable in $\xx_{t}$ except possibly $x_{1;t},x_{2;t},x_{3;t}$ satisfies the statement of Theorem~\ref{thm module main}. This argument did not use the fact that $\mu$ is admissible. The same argument, using the sequence $\mu'$, proves that every variable in $\xx_{t}$ except possibly $x_{1';t},x_{2';t},x_{3';t}$ satisfies the statement of Theorem~\ref{thm module main} as well. Since the sets $\{1,2,3\}$ and $\{1',2',3'\}$ are disjoint, this completes the proof in the case where the mutation sequence starts with an R3 move.


\subsection{Proof in the case where \texorpdfstring{$K$}~ has a bigon} \label{proof bigon}
\begin{figure}
\begin{center}
\scalebox{0.8}{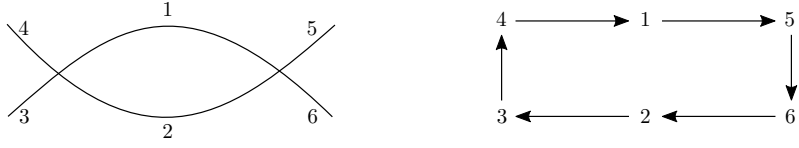}  
\caption{A bigon in $K$ with labeled segments and the corresponding subquiver of $Q$.}
\label{fig bigon}
\end{center}
\end{figure}

The structure of the proof is similar to the previous one, but simpler.  We label the segments of the bigon by 1 and 2, and the other segments incident to the bigon by 3,4,5,6 as shown in the left picture of Figure~\ref{fig bigon}. The right picture shows the corresponding subquiver of $Q$, where we removed the 2-cycle corresponding to the bigon. As before, let $\mu$ be the mutation sequence constructed in (\ref{eq mutation sequence}). Since $K$ has a bigon, the first two mutations in $\mu$ are at vertices 1 and 2. Thus we have $\mu=\mu_1\mu_2\mu'\mu_2\mu_1.$
We illustrate this decomposition in the following diagram of seeds and mutations in the cluster algebra.
 
\begin{equation}\label{eq mutation diagram2}
 \xymatrix@C55pt{\zS_0\ar[r]^{\mu_1,\mu_2} &\zS_{t_1}\ar[r]^{\mu'}&\zS_{t_2}\ar[r]^{\mu_2,\mu_1}&\zS_t.}
\end{equation}

 Let $\overline{K}$
denote the diagram obtained from $K$ by reducing the bigon $1,2.$ Let $\Qbar$ denote the corresponding quiver. It is obtained from $Q$ by mutating at the vertices 1 and 2 and then removing the vertices 1 and 2.

Again, we start by proving that $\zs\colon Q\to Q_t^{\textup{op}}$ is an isomorphism of quivers. Consider the pair $(K,Q)$ of the link diagram $K$ and the corresponding quiver $Q$. The  reduction of the bigon $(i,j)$ produces the pair $(\overline{K}, \overline{Q})$, where $\overline{K}$ is obtained from $K$ by replacing the bigon by a crossing point $p$, and  $ \overline{Q} $ is obtained from $Q$ by mutating at $i$ and $j$ and then removing the vertices $i$ and $j$.
By induction, we may assume that $\overline{\zs}\colon \overline{Q}\to \overline{Q}_{\overline{t}}^{\textup{op}}$ is an isomorphism of quivers which corresponds to replacing the diagram $\overline{K}$ by its mirror image $\overline{K}^{\textup{op}}.$ 
This map maps a local neighborhood of the crossing point $p$ with incident segments labeled $a,b,c,d$ in clockwise order to a local neighborhood of a crossing point $\zs(p)$ with incident segments labeled $\zs(a),\zs(b),\zs(c),\zs(d)$ in counterclockwise order. We can reinsert the bigon $(i,j)$ at $p$ and the bigon $(\zs (i),\zs (j))$ at $\zs(p)$, leading to quivers $Q$ and $Q'$ that locally are given by
\[\xymatrix{
    a\ar[r]&i\ar[r]&b\ar[d]\\
    d\ar[u]&j\ar[l]&c\ar[l] 
} \qquad\textup{and}\qquad
\xymatrix{
    \zs(a)\ar@{<-}[r]&\zs (i) \ar@{<-}[r]&\zs(b)\ar@{<-}[d]\\
    \zs(d)\ar@{<-}[u]&\zs (j)\ar@{<-}[l]&\zs(c)\ar@{<-}[l] 
}
\]
Moreover $\overline{Q}^{\textup{op}}$ is obtained from $Q'$ by mutating in $\zs (i),\zs (j)$ and then deleting the vertices $\zs (i),\zs (j)$. Thus $Q'=Q_t^{\textup{op}}$. 
Since $\zs$ is the composition of $\overline{\zs}$ with the transposition $(ij)$, it follows that $\zs$ is an isomorphism from $Q$ to $Q_t^{\textup{op}}$.

 Now we prove the identity on the $F$-polynomials in Theorem~\ref{thm module main}.
We may assume by induction that for every vertex $i$ of $\overline{Q}_0$, we have
 
\begin{equation}
 \label{Induction hypothesis 2} 
 F_{i;t_2}^{t_1} = F_{\overline{T}(\overline {\zs}(i))},
\end{equation}
 where
 
\begin{itemize}
 \item[-] $F_{i;t_2}^{t_1}$ is the $F$-polynomial of the cluster variable $x_{i;t_2}$ in the cluster algebra $\cala(\overline{Q})$ relative to the initial seed $\zS_{t_1}=(\xx_{t_1},\overline{\yy}, \overline{Q})$ with principal coefficients $\overline{\yy}=(\bar y_3,\ldots,\bar y_{2n})$;
 \item[-] $\overline{\zs}$ is the permutation corresponding to the mutation sequence $\mu'$ (Note that $\zs=\overline{\zs}\circ (1,2)$.); 
 \item [-]  ${\overline{T}(\overline {\zs}(i))}$ is the representation of $\overline{Q}$ corresponding to the segment $\overline{\zs}(i)=\zs(i)$ of $\overline{K}$;
 \item [-] $ F_{\overline{T}(\overline {\zs}(i))}$ is the $F$-polynomial of $ {\overline{T}(\overline {\zs}(i))}$.
\end{itemize}

We want to show that
\begin{equation}
\label{eq goal2}
 F_{i;t}^{t_0} = F_{T(\zs(i))}.
\end{equation}

The left hand side of this equation is the $F$-polynomial  of the cluster variable $x_{i,t}$ in the cluster algebra $\cala(Q)$ relative to the initial seed  $\zS_{t_0}=(\xx_{0},\yy_0, Q)$ with principal coefficients $\yy=( y_1,\ldots, y_{2n})$.

\medskip

{\bf First Case.} We will first assume that $i \notin \{1,2\}$.
Then $F_{i;t}^{t_1}=F_{i;t_2}^{t_1}$, because the last two mutations in $\mu$ are in 1 and 2. In order to calculate $F_{i;t}^{t_0}$ we need to apply a new function $\psi$, similar to the function $\varphi$ in  section~\ref{proof w/o bigon}.  To define $\psi$ we adapt Definition~\ref{def phi} to the current situation.

\begin{definition}
 \label{def psi}
 Let $\psi\colon\mathbb{Z}[\bar y_3,\ldots,\bar y_{2n}] \to \mathbb{Z}[y_1,\ldots,y_{2n}]$ be the ring homomorphism such that 
   \[ F_{i;t_2}^{t_0} = \textstyle \psi( F_{i;t_2}^{t_1})\cdot \prod_{j =1}^2  (F_{j;t_1}^{t_0})^{\overline{g}_j} .\]
\end{definition}

The function $\psi$ is determined by the following lemma.
\begin{lemma}
    \label{lem psi}
    The values of $\psi $ on the $\overline{y}_i$ are as follows. 
    \[\begin{array}{ll}
       \psi(\overline{y}_3)  = y_3(1+y_2)     & 
       \psi(\overline{y}_4)  = y_1y_4(1+y_1)^{-1}\\
       \psi(\overline{y}_5)  = y_5(1+y_1)     & 
       \psi(\overline{y}_6)  = y_2y_6(1+y_2)^{-1}
    \end{array}
    \]
       and $\psi(\overline{y_j})=y_j$ if $j>6.$
\end{lemma}

\begin{proof}
  Performing the two mutations from $t_0$ to $t_1$, one gets    \[y_{3;t_1}=y_3, \  \hat{\bar y}_3= y_{3;t_1} x_{2;t_1}x_{4;t_1}x_{6;t_1}^{-1}, \textup{ and } F_{2;t_1}^{t_0}=1+y_2,\  F_{4;t_1}^{t_0}=F_{6;t_1}^{t_0}=1.\] Thus $\psi(\overline{y}_3)  = y_3(1+y_2)$. Similarly
    \[y_{4;t_1}=y_1y_4,\  \hat{\bar y}_4= y_{4;t_1}x_{5;t_1}x_{1;t_1}^{-1}x_{3;t_1}^{-1} \textup{ and } F_{1;t_1}^{t_0}=1+y_1,\  F_{3;t_1}^{t_0}=F_{5;t_1}^{t_0}=1. \] Thus
     $\psi(\overline{y}_4)  = y_1y_4(1+y_1)^{-1}$. The other identities are proved mutatis mutandis.
\end{proof}

We are now ready to prove the theorem.

\begin{figure}
     \centering
     \huge\scalebox{0.4}{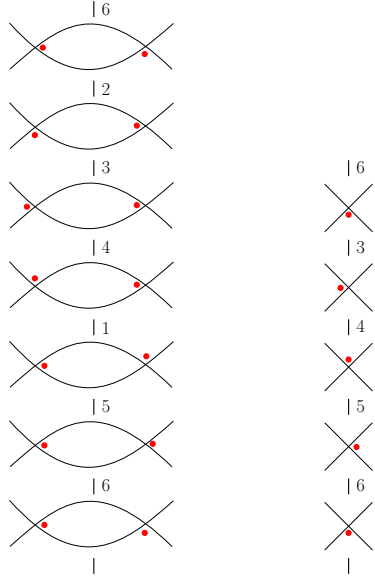}
     \caption{The local structure of the poset of Kauffman states of $K$ on the left,  and the poset of the Kauffman states of $\overline{K}$ on the right.}
     \label{fig bigonposet}
 \end{figure}

The subposet on the segments 1,2,3,4,5,6 of the poset of Kauffman states of $K$ is shown in the left picture in Figure~\ref{fig bigonposet}. It is a linear poset with edge weights $\ldots,6,5,1,4,3,2,6,5,1,\ldots$ when going up in the poset.  The right hand picture of the figure shows the subposet on the segments 3,4,5,6 of the poset of Kauffman states of $\overline{K}$. It  also is linear, but with edge weights $\ldots,6,5,4,3,6,5,\ldots$ when
 going up in the poset. Thus the poset on the right is obtained from the poset on the left by contracting the edges that are labeled by 1 or 2.

On the other hand, the function $\psi$ does exactly the same on the level of $F$-polynomials. Indeed, using Lemma~\ref{lem psi}, we have \[\psi(1+\overline{y}_3+\overline{y}_3\overline{y}_6)=1+y_3(1+y_2)+y_3(1+y_2) y_2y_6(1+y_2)^{-1} = 1+y_3+y_2y_3+y_2y_3y_6.\] Thus $\psi$ sends the linear poset with two consecutive edges with weights 3 and 6 to the linear poset with three consecutive edges with weights 3, 2 and 6. Similarly, $\psi$ sends the linear poset with two consecutive edges with weights 5 and 4 to the linear poset with three consecutive edges with weights 5, 1 and 4.

 This shows that for $i\ne 1,2$ we have $F_{i;t}^{t_0}=F_{\overline{T}(\zs(i))}$
is indeed the $F$-polynomial of the representation $T(\zs (i))$ of the segment $\zs (i)$ in $K$. \qed

\medskip

{\bf Second Case.} Now assume that $i \in \{1,2\}$.
Note that  $F_{i;t_2}^{t_0}=F_{i;t_1}^{t_0}=1+y_i$, for $i=1,2$, because the sequence $\mu'$ does not mutate at the vertices $1,2$. However, to calculate $F_{i;t}^{t_0}$ we  need to apply the last two mutations of the sequence $\mu.$ 
We will do this computation for $i=1$ only, because for $i=2$ it is entirely similar. 

Since $\zs$ is an isomorphisms between $Q_{t}$ and $Q^{op}$, we know that the local configuration at $\zs (2)=1$ of the quiver $Q_{t}$ is $\xymatrix{\zs (3)\ar[r]&1\ar[r]&\zs(6)}.$ The quiver $Q_{t_2}$ is obtained from this one by mutation in 1 and 2, thus locally $Q_{t_2}$
is $\xymatrix{\zs (3)\ar@/^8pt/[rr]\ar@{<-}[r]&1\ar@{<-}[r]&\zs(6)}.$ 
From the first part of the proof we know  that 
\begin{equation}\label{eq 50}
F_{\zs (i);t_2}^{t_0}=F_{T(i)}, \textup{ for $i\ne 1,2$.}
\end{equation}
On the other hand, $ F_{1;t_2}^{t_0}=1+y_1$.  We want to show
\begin{equation}
\label{eq 51} 
F_{1;t}^{t_0}=F_{T(2)}. 
\end{equation}

The seed $t$ is obtained from $t_2$ by mutating at 1 and 2 and these two mutations commute. Using the mutation formula for $F$-polynomials, we have
\begin{equation}
\label{eq 52}
F_{1;t}^{t_0}= 
\frac{y_{1;t_2}\, F_{\zs (3);t_2}^{t_0} + F_{\zs(6);t_2}^{t_0}}
{(1\oplus y_{1;t_2} )  \, F_{1;t_2}^{t_0} } . 
\end{equation}

Our first goal is to compute $y_{1;t_2}$. We need the following lemma. 
\begin{lemma}
 \label{lem CG}
 We have the following relation between the $C$-matrix and the $G$-matrix at $t$.
 \[ C_{t}^{B_{t_0}}=\zs \left(G_{t}^{B_{t_0}}\right)^{\!\!\mathsf{T}}\zs,\]
where $B_{t_0}$ is the skew-symmetric matrix associated to the quiver $Q=Q_{t_0}$.
\end{lemma}
\begin{proof}
 [Proof of the lemma.]
 From Theorem \ref{thm NZ}  we have
\begin{equation}
\label{eq GCformula} 
 \left(G_{t}^{B_{t_0}}\right)^{\!\!\mathsf{T}} = \left(C_{t}^{B_{t_0}}\right)^{-1}
= C_{t_0}^{-B_{t}},\end{equation}
where $B_{t}$ is the skew symmetric matrix for $Q_{t}$. Since $\zs Q_{t}=Q^{\textup{op}}$, we have ${-B_{t}}=\zs B \zs.$ Thus equation~\ref{eq GCformula} implies
\begin{equation}
\label{eq 53}
 \left(G_{t}^{B_{t_0}}\right)^{\!\!\mathsf{T}} = C_{t_0}^{\,\zs B_{t}\zs}.
\end{equation}
The right hand side of this equation is the $C$-matrix of the reverse mutation sequence $\overleftarrow{\mu}$ applied to $\zs Q_{t}$. Recall that $\mu$ is of the form 
$\mu= \mu_{\textup{red}}\,\mu_{\textup{Hopf}}\, \overleftarrow{\mu_{\textup{red}}^\zs}$ and takes $Q$ into $Q_{t}$. Thus the reverse sequence is of the form 
$\overleftarrow{\mu}
=
\overleftarrow{(\overleftarrow{\mu^\zs_{\textup{red}} }\,{\mu_{\textup{Hopf}} }\, {\mu_{\textup{red}}})}
={\mu^\zs_{\textup{red}} }\,{\mu_{\textup{Hopf}} }\,\overleftarrow{\mu_{\textup{red}}}$, where the last equality holds because the Hopf mutation sequence is its own inverse. Moreover $\overleftarrow{\mu} $ takes $Q_t$ into $Q$.
Applying the permutation $\zs$, we see that 
$\zs \overleftarrow{\mu}={\mu_{\textup{red}} }\,{\mu_{\textup{Hopf}} }\,\overleftarrow{\mu^\zs_{\textup{red}}}$, because the Hopf mutation sequence is invariant under $\zs$. 
Thus $\zs \overleftarrow{\mu}=\mu$.  Therefore the $C$-matrix $C_{t_0}^{\,\zs B_{t}\zs}$ of $\overleftarrow{\mu}$ applied to $\zs Q_{t}$ is equal to the conjugate by $\zs$ of the $C$-matrix $C_{t}^{B_{t_0}}$ of $\mu$ applied to $Q$. 
Thus equation (\ref{eq 53}) becomes 
$\left(G_{t}^{B_{t_0}}\right)^{\!\!\mathsf{T}} = \zs C_{t}^{B_{t_0}}\zs$, and the lemma is proved.
\end{proof}

\begin{lemma}
    \label{lem 519}
    Suppose $K$ is a diagram that contains the bigon in Figure~\ref{fig bigon} and $i$ is a segment in $K$ such that $T(i)$ has the same dimension $d$ at vertices 1 and 2. Then we have one of the following.
    \begin{enumerate}
        \item [(i)] The dimension of $T(i)$ is equal to $d$ at vertices 3 and 4, and it is equal to $d-1,d$, or $d+1$ and vertices 5 and 6.
        \item [(ii)] The dimension of $T(i)$ is equal to $d$ at vertices 5 and 6, and it is equal to $d-1,d$, or $d+1$ and vertices 3 and 4.
    \end{enumerate}
\end{lemma}
\begin{proof}[Proof of the lemma.]
    Consider the lattice of Kauffman states relative to $i$. Restricting the states to the vertices at the bigon produces a finite linear poset, which is equal to an interval in the poset illustrated in the left picture in Figure~\ref{fig bigonposet}. 
    Since $T(i)$ has the same dimension $d$ at vertices 1 and 2, there are precisely $d$ edges labeled 1 and $d$ edges labeled $2$ in this interval. Thus the sequence of edges in the interval is of one of the following two forms
    \[\begin{array}
        {ll}
        u,1,4,3,2,6,5,1,4,3,2,\ldots,1,4,3,2,v &\textup{with } u=\emptyset,5 \textup{ or } 6,5, \textup{ and } v=\emptyset,6 \textup{ or } 6,5;\\
        u,2,6,5,1,4,3,2,6,5,1,\ldots,2,6,5,1,v&\textup{with } u=\emptyset,3 \textup{ or } 4,3, \textup{ and } v=\emptyset,4 \textup{ or } 4,3.\\
    \end{array}\]
    In the first case, the dimension of $T(i)$ is equal to $d$ at vertices 3 and 4, and it is equal to $d-1,d,$ or $d+1$ at vertices 5 and 6 depending on the sequences $u$ and $v$. The second case is analogous.
\end{proof}

We are now ready to compute the coefficient $y_{1;t_2}$ in equation (\ref{eq 52}).  Let $e_1=(1,0,\ldots,0)$ be the first standard basis vector. 
\begin{lemma}
 \label{lem y1}
 In the situation above, we have
 \[ y_{1;t_2} = \yy^{\underline{\dim} \,T(2)-\underline{\dim}\, T(3)+e_1},\]
 and all exponents are non-negative.
\end{lemma}
\begin{proof}
 Be definition  $y_{1;t_2}=\yy^{c_{1;t_2}^{t_0}}$, where the exponent is the first $c$-vector of the seed $t_2$. Let us simply write $c $ for ${c_{1;t_2}^{t_0}}$ in this proof. Thus $c$ is the first column of the matrix $C_{t_2}^{B_{t_0}}$. Since the seed $t$ is obtained from $t_2$ by mutations at 1 and 2, we have 
 \[ c= - \textup{ first column of }  C_{t}^{B_{t_0}}.\]
 Using Lemma~\ref{lem CG} and applying the transpose, we get
 \[ c= - \textup{ first row of }  \zs G_{t}^{B_{t_0}} \zs,\]
 and since $\zs$ permutes 1 and 2,
\begin{equation}
\label{eq 54}
 c= - \textup{ second row of }  G_{t}^{B_{t_0}} \zs.
\end{equation}
The columns of $G_{t}^{B_{t_0}}$ are the $g$-vectors of the seed $t$. The negative entries of the $g$-vectors correspond to the socle of the representations $T(j)$, for $j\ge 3$. More precisely, the $i$-th entry of the $j$-th column  of $G_{t}^{B_{t_0}}$ is negative if and only if the simple module $S(i)$ is a direct summand of the socle of $T(j)$, for $j\ge 3$.
Even stronger, if the entry is negative then it is equal to the multiplicity of $S(i)$ in the socle of $T(j)$. 

We have shown in \cite[Theorem 1.2]{BMS} that the submodule lattice of $T(j)$ is isomorphic to the lattice of Kauffman states of $K$ relative to the segment $j$. This implies that the socle of $T(j)$ is the direct sum of simples $S(i)$ such that the minimal Kauffman state relative to segment $j$ has exactly two markers on the segment $i$. In particular, the multiplicity of $S(i)$ is at most one.

Going back to equation (\ref{eq 54}) and taking into account the sign and the action of $\zs$, we see that
\begin{equation}
\label{eq 55}
c_j>0 \ssi S(2) \textup{ is a submodule of $T(\zs (j))$; and in this case $c_j=1$}.
\end{equation}

 Because of the definition of the representation $T(\zs (j))$ in \cite{BMS} and using the fact that from Figure~\ref{fig bigon} there is only one arrow $2\to 3$ in $Q$ that starts at vertex $2$, we can reformulate the condition in (\ref{eq 55}) by saying that the dimension of $T(\zs (j))$ at vertex 2 is strictly larger than its dimension at vertex 3.
 Moreover, these dimensions can differ at most by one by \cite[Proposition 5.9]{BMS}. We therefore have
 \[ c_j>0 \ssi \dim T(\zs (j))_2 = \dim T(\zs (j))_3 +1,\textup{ and in this case $c_j=1$},\]
where $T(\zs(j))_x$ is the vector space of $T(\zs (j))$ at vertex $x$. 

Now we use the symmetry of the dimensions in Theorem~\ref{thm dim sym intro} to rewrite the condition as follows
 \[ c_j>0 \ssi \dim T(2)_{\zs (j)} = \dim T(3)_{\zs (j)} +1 \textup{; and in this case $c_j=1$}. \]
If we take $\zs (j)=4$, referring to the labeling in Figure \ref{fig bigon}, we see that
\begin{enumerate}
\item 
$\dim T(3)_{\zs (j)} = 0, $ because segments 3 and 4 bound the same region in $K$,
\item  
$\dim T(2)_{\zs (j)} > 0, $ because segments 2 and 4 don't.
\end{enumerate}
Here we point out that segments 3 and 4 cannot be equal, because $K$ does not contain loops. 
This shows that $c_{4}=1$. In particular, the vector $c$ has a positive entry, and by the sign coherence of $c$-vectors \cite{GHKK}, this implies that $c$ has no negative entries. Thus all entries of $c$ are either 1 or 0. Hence
\[c_j= \max (0 ,  \dim T(2)_{\zs (j)} - \dim T(3)_{\zs (j)} ).\]

This almost completes the proof; it only remains to show
\begin{equation}
\label{eq 56}
\dim T(2)_{\zs (j)} - \dim T(3)_{\zs (j)} <0 \ssi \zs (j)=1. 
\end{equation}

If $\zs (j)=1$ then $\dim T(\zs (j))_2=0$, because segments 1 and 2 bound the same region in $K$, and $\dim T(\zs (j))_3>0$, because segments 1 and 3 do not bound a region in $K$. Thus the inequality holds.

Now suppose the inequality in (\ref{eq 56}) holds for some $\zs (j)\ne 1$. By Theorem~\ref{thm dim sym intro} this means $\dim T(\zs (j))_2-\dim T(\zs (j))_3 <0$.
Note that $\zs (j)\ne 2$,  because $T(2)_2=T(2)_3=0$, contradicting our inequality.

Assume first $\dim T(\zs (j))_1=\dim T(\zs (j))_2=d$. Then we must have the scenario (i) or (ii) of Lemma~\ref{lem 519}. However, scenario (i) is impossible, since we have seen above that $\dim T(\zs (j))_2-\dim T(\zs (j))_3 <0$. Hence we are in scenario (ii) and thus the interval in the proof of the lemma is of the form
\[ u,2,6,5,1,4,3,2,6,5,1,\ldots,2,6,5,1,v \quad \textup{with } u=3 \textup{ or } 4,3, \textup{ and } v=4,3
\]  again because $\dim T(\zs (j))_2-\dim T(\zs (j))_3 <0$.
Then, locally at the bigon, the maximal Kauffman state is reached  after the last three transpositions at 1,4,3. Thus the markers of the maximal state are between segments 2 and 3 at the left point of the bigon and between segments 1 and 2 at the right endpoint of the bigon. In particular, both markers are at the segment 2, and therefore the transposition  at 2 can be applied. This is a contradiction to the maximality of the state. 

Therefore $\dim T(\zs (j))_1\ne \dim T(\zs (j))_2$. 
Because of the linear structure of the poset in Figure~\ref{fig bigonposet}  the dimensions of $T(\zs (j))$ at vertices 1,2,3,4 differ at most by one, and when
going through the cycle 1,2,3,4,1, the dimension can increase at most once.
Because our inequality requires 
$\dim T(\zs (j))_2< \dim T(\zs (j))_3$, the only possibility is that  $\dim T(\zs (j))_1=\dim T(\zs (j))_3=\dim T(\zs (j))_4=\dim T(\zs (j))_2+1$. 
Again this would mean that the maximal Kauffman state is reached after the last three transpositions at 1,4,3, and we have seen above that this is impossible.
This proves  (\ref{eq 56}) and thereby Lemma \ref{lem y1}.
\end{proof}

Let's go back to equation (\ref{eq 52}).
Using Lemma \ref{lem y1}, the formulas for the $F$-polynomials in equation (\ref{eq 50}) and 
$F_{1;t_2}^{t_0}=1+y_1$, we obtain
\begin{equation}
\label{eq 57}
(1+y_1)\, F_{1;t}^{t_0} = \yy^{\underline{\dim}T(2)-\underline{\dim}T(3) +e_1}\,F_{T(3)}+F_{T(6)}.
\end{equation}
Here we use the fact that the non-negativity of the exponents of $\yy_{1;t_2}$ implies that $1\oplus \yy_{1;t_2}=1.$
Thus in order to complete the proof of Theorem~\ref{thm module main}, we must show that this equation  holds if we replace $F_{1;t}^{t_0} $ by $F_{T(2)}$. 

\begin{figure}
\begin{center}
\LARGE\scalebox{0.6}{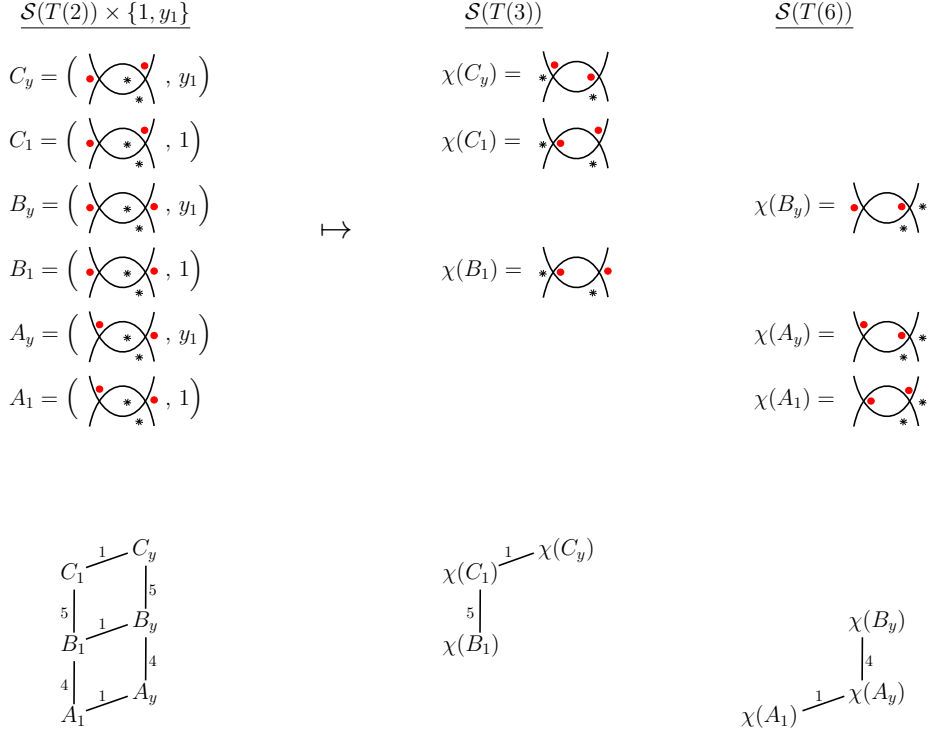}
\caption{Definition of the bijection $\chi$ (top) and the Hasse diagrams of the posets (bottom)}
\label{fig poset bij}
\end{center}
\end{figure}

Let $\cals(T(i))$ denote the set of Kauffman states of $K$ with respect to segment $i$ and consider its usual poset structure given by the Kauffman transpositions.  We define a bijection 
\[\chi\colon \cals(T(2))\times\{1,y_1\} \to \cals(T(3))\sqcup \cals(T(6)) \] as shown in Figure~\ref{fig poset bij}. 
In this figure, we only show the state markers at the bigon. All other markers are preserved under $\chi$. On the left, there are exactly three configurations for the Kauffman states at the bigon and each gives rise to two elements in the product $\cals(T(2))\times \{1,y_1\}$. On the right hand side, there are exactly three configurations for the Kauffman states at the bigon for each of $T(3)$ and $T(6)$. 
Moreover, for the five regions at the bigon, a region has a decoration (marker or star) in a state of $T(2)$ if and only if the same region has a decoration in the image of the state under $\chi$. 

Now consider the terms of the $F$-polynomials on either side. From the Hasse diagrams at the bottom of Figure~\ref{fig poset bij}, we see that on the left, we have $(1+y_4+y_4y_5)(1+y_1)$, and on the right, we have
$1+y_5+y_1y_5$ and $1+y_1+y_1y_4$. 
The minimal state $A_1$ of $T(2)$ contributes 1 to $F_{T(2)}$. Its image $\chi(A)$ is the minimal state of $T(6)$ and it contributes 1 to the right hand side of (\ref{eq 57}).
The states of $T(2)$ that restrict to states $A_y, B_y$ in Figure~\ref{fig poset bij} contribute the same monomials to $F(2)$ as their images under $\chi$ contribute to $F_{T(6)}$, because the  
 subposet on the elements $A_1,A_y,B_y$ on the left is isomorphic to the second poset on the right, with the same edge weights.

On the other hand, the maximal state $C_y$ on the left contributes $\yy^{\underline{\dim}\,T(2) +e_1}$ to $F_{T(2)}(1+y_1)$. It's image $\chi(C_y)$ is the maximal state of $T(3)$ and therefore it contributes $\yy^{\underline{\dim}\,T(3)} $ to $F_{T(3)}$. The difference between these two contributions is exactly captured by the coefficient 
$\yy^{\underline{\dim}T(2)-\underline{\dim}T(3) +e_1}$ in equation (\ref{eq 57}).  So the contribution is the same on both sides for the maximum elements of the posets. 
Again we find the poset for $T(3)$ on the right of Figure~\ref{fig poset bij} as the  subposet on elements $B_1,C_1,C_y$ of the poset on the left of that figure. Moreover $\chi$ is a weight preserving bijection between these two posets. 
This shows that 
\[(1+y_1)\, F_{T(2)} = \yy^{\underline{\dim}T(2)-\underline{\dim}T(3) +e_1}\,F_{T(3)}+F_{T(6)}.
\] as desired. Thus equation (\ref{eq 57}) is valid and this completes the proof of Theorem~\ref{thm module main}.
Part (a) of Theorem~\ref{thm main} now follows from Theorem~\ref{thmAlexPol}.

\subsection{Proof of part (c) of Theorem~\ref{thm main}} 

We need to show that $\phi_\zs$ is a cluster automorphism of order 2. The map $\phi_\zs $ sends the initial seed $(\xx_0,Q)$ to the knot seed $(\xx_t,Q_t)$ and induces an isomorphism $\zs\colon Q\to Q^{\textup{op}}_t$. Therefore \cite[Lemma 2.3]{ASS} implies that $\phi_\zs$ is a cluster automorphism. 

Recall that $\phi_\zs $ is given by the mutation sequence $\mu$ of section~\ref{sect mutation} followed by the permutation $\zs$ on the cluster variables in the knot cluster.
To show that the order of $\phi_\zs$ is two, it suffices to show that $\mu\mu^\zs$ is the identity. 
We have
 \[\mu\mu^\zs=  \mu_{\textup {red}}\ \mu_{\textup {Hopf}} \ 
 \overleftarrow{\mu^\zs_{\textup {red}}} \  \mu_{\textup {red}}^\zs\ \mu_{\textup {Hopf}}^\zs \ 
 \overleftarrow{\mu^{\zs\,\zs}_{\textup {red}}}
 =  \mu_{\textup {red}}\ \mu_{\textup {Hopf}} \ 
  \mu_{\textup {Hopf}}^\zs \ 
 \overleftarrow{\mu_{\textup {red}}}
  =  \mu_{\textup {red}}\ \overleftarrow{\mu_{\textup {red}}}
=1\qedhere\]
 \qed

\section{Creation of bigons} \label{sect bigons}
In this section we prove that any link diagram not containing a bigon admits a sequence a diagram R3 moves (RD3 moves) that creates a bigon.

Let's first introduce some vocabulary. Recall that a link is an embedding of a finite number of circles in $\mathbb{R}^3$. An \emph{arc} of a circle is a subcurve that connects two points on the circle. An arc is \emph{trivial} if it consists of a single point. 
   A \emph{strand} in a link is the image of a nontrivial arc of one of the circles under this embedding. The image of a strand in a  diagram of the link will also be called a strand. A strand in a link diagram is said to be \emph{simple} if it has no selfcrossings; however the endpoints of a simple strand may coincide.

\begin{figure}
\begin{center}
\scalebox{1}{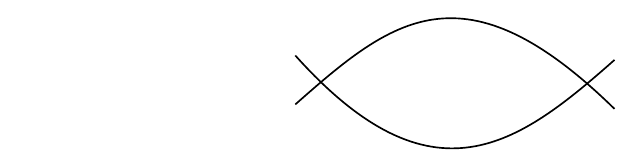}
\caption{The left picture shows an example of a generalized loop. The right picture  shows an example of a generalized bigon that is  not minimal. }
\label{figgenbigon1}
\end{center}
\end{figure}

\begin{definition} Let $K$ be a link diagram. \begin{enumerate}[(a)]
    \item A \emph{generalized loop} in $K$ is a pair $(a,s)$, where $a$ is a crossing point and $s$ is a  strand from $a$ to $a$ that does not contain a strict substrand from $a$ to $a.$
    \item  A \emph{generalized bigon} in $K$ is a region bounded by a quadruple $(a,b,s,s')$, where $a,b$ are two distinct crossing points of $K$ and $s,s'$ are two distinct simple strands with endpoints $a$ and $b$ such that:
\begin{itemize}
    \item[(i)] $s$ and $s'$ do not cross each other;
    \item[(ii)] the concatenations $ss'$ and $s's$ are not strands.
 \end{itemize}
 A generalized bigon is said to be \emph{minimal} if it does not contain another generalized bigon in its interior.
\end{enumerate}
\end{definition}
 An example of a generalized loop is given in the left picture of Figure~\ref{figgenbigon1}. The right picture of the figure shows an example of a generalized bigon. Note that this generalized bigon is not minimal, since it contains a smaller bigon.

\begin{lemma}\label{lem 0}
Let $K$   be a  nontrivial diagram of a prime link.

(a) Every generalized loop in $K$ contains a simple generalized loop as a substrand.

(b) Every generalized loop in $K$ contains a generalized bigon.
\end{lemma}
\begin{proof}
(a) Let $(a,s)$ be a generalized loop in $K$. If it is simple there is nothing to show.  Otherwise, let $a_1$ be the first self-crossing point of $s$ and let $(a_1,s_1)$ be the generalized loop given by the substrand $s_1$ of $s$. Then $(a_1,s_1)$ has fewer self-crossing points than $(a,s)$. Repeating this procedure will find a simple generalized loop as a substrand of $(a,s).$

(b) Let $(a,s)$ be a generalized loop. By part (a) we may assume it is simple.  
Because of our assumption that the link is prime and that our nontrivial diagram $K$ does not admit a Reidemeister I move, there must be a strand $t$ that crosses $s$ at two crossing points $b$ and $c$. Let's say $t$ enters the loop at $b$ and exits the loop at $c$. Let $t'$ be the substrand of $t$ from $b$ to $c$ and $s'$ be the substrand of $s$ from $b$ to $c$. 
If $t'$ is simple, then $(b,c,s',t')$ is a generalized bigon. Otherwise, $t'$ contains a generalized loop as a substrand. This generalized loop lies in the interior of the generalized loop $(a,s).$ Continuing this argument, we will find a generalized bigon, because $K$ can only contain a finite number of nested generalized loops.
\end{proof}
\begin{prop}
    Every nontrivial prime link diagram contains a generalized bigon.
\end{prop}
\begin{proof}
    Let $a$ be a crossing point. Walking in any direction along $K$ starting at $a$ will eventually lead back to $a$, thereby describing a generalized loop. Now the result follows from Lemma~\ref{lem 0}.
\end{proof}

\begin{lemma} \label{lem::generalizedbigon}
Let $B=(a,b,s,s')$ be a minimal generalized bigon. Then, any strand $s_0$ entering $B$ through one of the strands has to exit $B$ through the other strand and $s_0$ has no self-crossing inside $B$.
\end{lemma}

\begin{proof}
If $s_0$ enters and exits $B$ through the same strand, then $B$ contains a smaller generalized bigon, which contradicts its minimality.
Moreover, if $s_0$ has a self-crossing, then $B$ also contains a smaller generalized bigon by Lemma~\ref{lem 0}.
\end{proof}

This lemma allows us to define a partial order on the set of all crossing points of $K$ in the interior of a given minimal generalized bigon $B = (a,b,s,s')$.

\begin{figure}
\begin{center}
\small\scalebox{1}{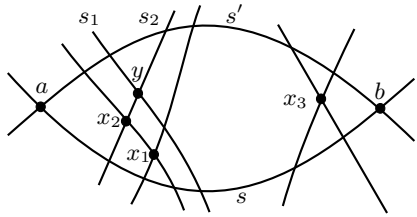}
\caption{An example of the order relation on crossing points within a generalized bigon. In this situation we have $x_1<y$ and $x_2<y$ (also $x_1<x_2$). The points $x_3$ and $y$ are not comparable. }
\label{figgenbigon}
\end{center}
\end{figure}

Let $x$ and $y$ be two such crossing points.
Then, $y$ lies on the intersection of two strands $s_1$ and $s_2$. We define $x \leq y$ if $x$ is inside the closed region bounded by $s_1$, $s_2$ and $s$. See Figure~\ref{figgenbigon} for examples.

It is easy to see that this order is well defined, thanks to Lemma \ref{lem::generalizedbigon}.

\begin{lemma} \label{lem::triangularregion}
Consider a minimal generalized bigon $B = (a, b, s, s')$. Then, either $B$ is a bigon, or there exists a triangular region in $K$ whose interior lies in the interior of $B$ and that has at least one side on the boundary of $B$. 

\end{lemma}

\begin{proof}
Suppose $B$ is not a bigon. 
If there is a crossing point in the interior of $B$, then any minimal element in the poset of crossing points defines a triangular region inside $B$ with one side on the strand $s$.
Similarly, any maximal element defines a triangular region inside $B$ with one side on the strand $s'$.

If there is no crossing point in the interior of $B$, then the crossing point $a$ is part of a triangular region bounded by $s$, $s'$ and a strand through $B$. The same applies for the crossing point $b$.
\end{proof}

\begin{definition}
 \label{def admissible}
 We say that two triangular regions of the diagram $K$ are \emph{disjoint} if they do not have a common vertex. A sequence of RD3 moves is \emph{admissible} if, for every pair of consecutive moves, their two triangular regions are   disjoint.
\end{definition}

\begin{thm} \label{thm::sequenceR3}
  Every link diagram contains a bigon or it admits an admissible sequence of RD3 moves creating a bigon.
\end{thm}

\begin{proof}
Consider a minimal generalized bigon $B = (a,b,s,s')$ in that link diagram $K$.
We proceed by induction on the number of regions of $K$ that lie inside $B$.
If that number is $1$, then $B$ is a bigon and there is nothing to prove.

Otherwise, by Lemma \ref{lem::triangularregion}, there exists a triangular region $\Delta_0$ inside $B$ with one side on the boundary of $B$.

If $\Delta_0$ is bounded both by a segment of $s$ and a segment of $s'$ and a third side $t$, then it also contains one of the crossing points $a$ and $b$. Apply the RD3 that moves $t$ across the crossing point $a$ or $b$. 
This case is illustrated in the top left picture in Figure~\ref{fig admissible}.

If $\Delta_0$ is bounded by $s$ and two strands through $B$ intersecting at a crossing point $x$ in $B$, apply the RD3 that moves $s$ across $x$.
Similarly, if that triangular region is bounded by $s'$ and two strands through $B$ intersecting at a crossing point $x'$ in $B$, apply the RD3 that moves $s'$ across $x'$.
This case is illustrated in the bottom left picture in Figure~\ref{fig admissible}.

In each case, we obtain a slightly different link diagram $K'$ in which the bigon $B$ is still minimal and contains one less region.

 Now we show that the sequence of RD3 moves can be chosen to be admissible. We distinguish two cases.

(1) First suppose that the previous RD3 move was at a triangular region $\zD_0$ that contains one of the two points $a$ or $b$; say $a$. This situation is illustrated in the top right picture of Figure~\ref{fig admissible}. Then our next RD3 move forms an admissible sequence if and only if its triangular region does not  contain the point $a$. There is at most one triangular region in $B$ that contains the point $a$. It is labeled $\zD_1$ in the top right picture in the figure. If the region labeled $A$  contains a crossing point then it contains a triangular region at one of the strands $s$ or $s'$, and we may perform the next RD3 move at that region. Otherwise, the region $A$ contains a triangular region that contains the point $b$, and we may perform the next RD3 move at that region. Note that this region may be all of $A$.

(2) Now suppose the previous RD3 move was at a triangular region $\zD_0$ that has exactly one boundary on the strand $s$ or $s'$; say $s$, see the bottom right picture  of Figure~\ref{fig admissible}. If there exists another strand that creates a crossing point in $B$, then there exists a triangular region at the strand $s'$, and we may perform the next RD3 move there.
If there exists another strand that goes through $B$ but does not produce any crossing points, then there exists a triangular region that contains the point $a$ or the point $b$, and we may perform the next RD3 move there, unless we are in the situation shown in the bottom right picture in Figure~\ref{fig admissible}. Note that $\zD_1$ and $\zD_2$ are triangular regions that both share a vertex with $\zD_0$, so we are not allowed to perform the next RD3 move at these regions. Instead we will use a triangular region that is elsewhere in $K$. 

Consider now the regions adjacent to $B$; we label them as in the bottom right picture in Figure~\ref{fig admissible}.
\begin{figure}
\begin{center}
\scalebox{0.9}{\small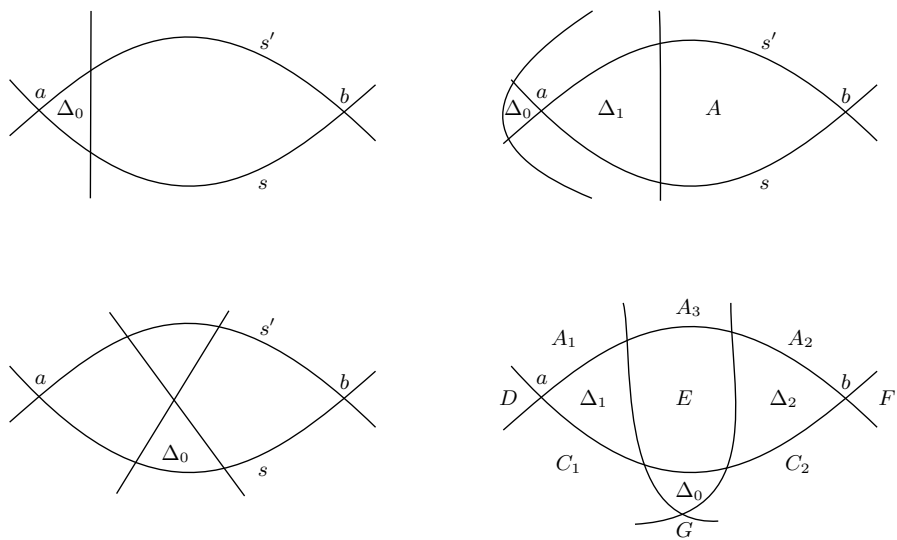}
\caption{Construction of an admissible sequence of RD3 moves}
\label{fig admissible}
\end{center}
\end{figure}

If region $A_1$ is triangular, then there is a different bigon formed by $A_1$ and $\zD_1$, and we may perform the next RD3 move at $A_1$ to obtain a bigon. The same argument applies to the region $A_2$.

If region $A_3$ is triangular, then we have another bigon formed by the regions $A_3,E$ and $\zD_0$, and we may perform the next RD3 move at $A_3$ and then at $\zD_0$ to obtain a bigon. 

Region $C_1$ cannot be triangular, because otherwise it would have been a bigon before the RD3 move at $\zD_0$. 
 The same argument applies to the region $C_2$.
 
 We therefore may assume that none of the regions $A_1,A_2,A_3,C_1,C_2$ is triangular.  The region $E$ has four sides. Thus in total there are at most six triangular regions in the figure: $\zD_0,\zD_1,\zD_2,D,F,G$. By Proposition~\ref{prop 8 triangles}, there exists another triangular region in $K$ that is disjoint from $\zD_0,\zD_1$ and $\zD_2$.  We may perform our next RD3 move at that region and then perform the RD3 moves at $\zD_1$ and $\zD_2$ to obtain an admissible sequence.
\end{proof}

\section{Symmetry of dimension} \label{sect sym dim}
In this section, we prove Theorem \ref{thm dim sym intro} on the dimension symmetry in the representations $T(i)$ of $Q$, where $i$ is a segment of $K$, or, equivalently, a vertex of $Q$.
We start by reviewing the definition of the dimensions of the vector spaces $T(i)_j$. We follow the exposition of \cite{BMS Add}. 
\subsection{Strings and their boundaries} 
\begin{figure}
\begin{center}
 \scalebox{1.5}{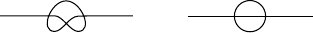}
 \caption{A trefoil string (left) and a Hopf link string (right).}
 \label{fig string}
\end{center}
\end{figure}

Let $K$ be a link diagram with $K_0, K_1,K_2 $ its sets of crossing points, segments and regions, respectively. Let $i\in K_1 $ be a segment. The \emph{$i$-th string} of $K$ is obtained from $K$ by removing an interior point from the segment $i$. The two ends of the segment $i$ that remain are called the start and the terminal of the string.  Figure~\ref{fig string} illustrates a string of the  trefoil knot and  one of the Hopf link.

When we draw the string in the plane, we can indefinitely extend its start segment horizontally to the left and its terminal segment horizontally to the right so that the complement of the string in the plane has two unbounded components. 
A segment of a string is said to be \emph{connecting} if it is incident to both unbounded components.

\begin{figure}
\begin{center}
 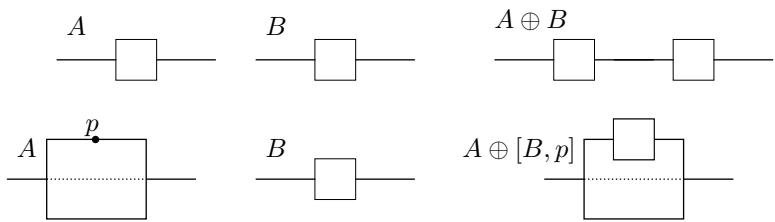
 \caption{Sum of two string $A\oplus B$ (top) and enclosed sum of two strings $A\oplus[B,p]$ (bottom).}
 \label{fig sum}
\end{center}
\end{figure}

In \cite{K}, Kauffman defined the following operations on strings. For illustrations see Figure~\ref{fig sum}.

\begin{definition}(Sum and enclosed sum)

(a)  If $A$ and $B$ are strings, their \emph{sum} $A\oplus B$ is obtained by connecting the terminal of $A$ to the start of $B$.

(b) If $A$ and $B$ are strings, $p$ is a point on a non-connecting edge of $A$ and $p$ is not a crossing point then the \emph{enclosed sum} $A\oplus [B,p]$ of $A$ and $B$ relative to $p$ is obtained by replacing the trivial string at $p$ by $B$. In this situation, the string $A$ is called a \emph{carrier} and the string $B$ is called a \emph{rider} in $A\oplus [B,p]$.

(c) A string is called \emph{atomic} if it is irreducible with respect to both sums.

\end{definition}


\begin{definition} The \emph{boundary} of a string is defined as follows.
 
 (a)  If $A$ is atomic, then $\partial A$ consists of all segments of $A$ that are incident to an unbounded component. 

(b) $\partial(A\oplus B)=\partial A \oplus \partial B$.

(c) $\partial (A\oplus [B,p])=\left\{
\begin{array}
 {ll}
 \partial A \oplus [\partial B,p] &\textup{if $p\in \partial A$};\\
 \partial A &\textup{otherwise.}
\end{array}\right.
$

\end{definition}

\begin{figure}
\begin{center}
\scalebox{0.7}{ 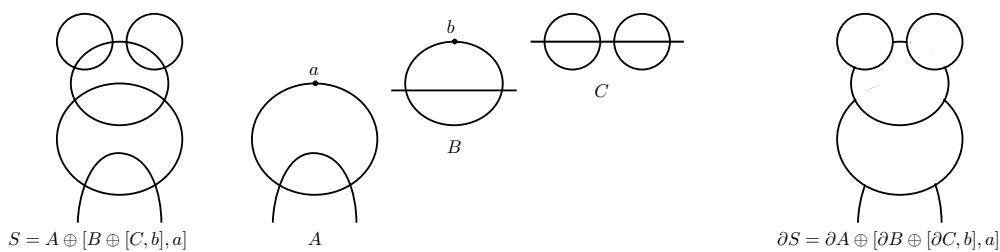}
 \caption{Computation of the boundary. The string $S$ is a double enclosed sum of the strings $A,B,C$. The boundary of $S$ is shown in the rightmost picture.}
 \label{figbdyex}
\end{center}
\end{figure}

We give an example of the computation of the boundary for a more complicated string in Figure~\ref{figbdyex}. The string there is a double enclosed sum $S=A\oplus [B \oplus[C,b],a]$ and the boundary is computed recursively  $\partial S=\partial A \oplus [\partial B\oplus[\partial C,b],a]$.

\subsection{The partition of $K_1$}
Let $i$ be a   segment of $K$. In \cite{BMS,BMS Add}, we define a partition $K_1=\sqcup_{d\ge 0} K(i,d)$ of the set of all segments of $K$ and use it to define the representation $T(i)$. 

\begin{definition}
 For $d=0$, we define
\[K(i,0)=\{j\in K_1\mid \textup{$j$ and $i$ bound the same region of $K$}\}\cup\{i\}\]
Thus $K(i,0)$ is the boundary of the $i$-th string of $K$ together with $i$.

Recursively, $K(i,d)$ is defined as the boundary of the string $K\setminus(\cup_{e<d} K(i,e))$.
\end{definition}

By definition of the representation $T(i)$, the dimension of the vector space $T(i)_j$ of $T(i)$ at the vertex $j$ of $Q$ is equal to $d$ if and only if $j\in K(i,d)$.

We now give an alternative description of $K(i,d)$, which uses the following terminology.
Given two segments $i,j\in K_1$, a curve in $\mathbb{R}^2$ is called a \emph{dimension curve} from  $j$ to  $i$ if it starts at a point on segment $j$, ends at a point on segment $i$ and does not go through a crossing point of $K$. 
Let $\dimprime(i,j)$ be the minimal number of crossings between the segments of $K$ and a dimension curve from segment $j$ to segment $i$.

With this notion, we have the following result.
\begin{prop}\label{prop partition}\cite{BMS Add}
A segment $j\in K_1\setminus(\cup_{e<d} K(e))$ lies in $K(d) $ if and only if one of the following two conditions hold.
\begin{enumerate}
 \item [{\rm (a)}] $\dimprime(i,j)=d$, or
 \item[{\rm (b)}] $\dimprime(i,j)=d+1$ and $K\setminus(\cup_{e<d} K(e))$ contains a string  of the form $ A\oplus[B,p]$ for two strings $A,B$ and $p\in \partial A$ an interior point of a non-separating segment, and $j\in \partial B$ is not incident to an unbounded component of the complement of $ A\oplus[B,p]$ in the plane.
\end{enumerate}
\end{prop}

\subsection{Proof of symmetry} We are now ready to prove the main result of this section. 

\begin{thm}\label{thm dim sym}
 Let $i,j$ be two segments of $K$ and let $T(i),T(j)$ be the corresponding representations of the quiver $Q$. 
 Then
 \[
\begin{array}{rcl}\dim T(i)_j &=&\dim T(j)_i.
 \end{array}\]
\end{thm}
\begin{proof}
 
 It suffices to show that for any two segments $i,j\in K_1$ we have $i\in K(j,d)$ if and only if $j\in K(i,d)$. If both dimensions are realized by a dimension curve then the equation holds because the same curve will realize the minimal number of crossings with $K$ in both directions.

Suppose therefore that $j\in K(i,d)$ and there is no dimension curve from $i$ to $j$ that crosses $K$ exactly $d$ times. Then we are in case (b) of Proposition~\ref{prop partition}. 

\subsubsection{The case $d=1$}
We start with the case $d=1$. Thus $K\setminus K(i,0)$ is an enclosed sum $A\oplus[B,p]$ with $p\in \partial A$ and $j$ a segment in $\partial B$ that is not incident to the unbounded component of the complement of $A\oplus[B,p]$. 
We may assume that $A$ and $B$ are both Hopf strings, since the general case is similar. The situation is illustrated in the  left picture of Figure~\ref{figdimsym1}. The string $K\setminus K(i,0)$ is drawn in red and blue. The string $A$ consists of the segments $a_1,e_1,e_2,c_1,h_2,d_1,b_1,$ and the point $p$ lies on the blue segment $h_2$. The string $B$ consists of the segments $c_1,j,h_2,f_1,h_1,g_1,d_1$. Then $K(i,1)$ is the boundary of $A\oplus[B,p]$ which is equal to $\partial A\oplus[\partial B,p]$. Thus $K(i,1)$ consists of the segment $j$ and all  segments whose labels have the subscript 1. 
\begin{figure}
\begin{center}
\scriptsize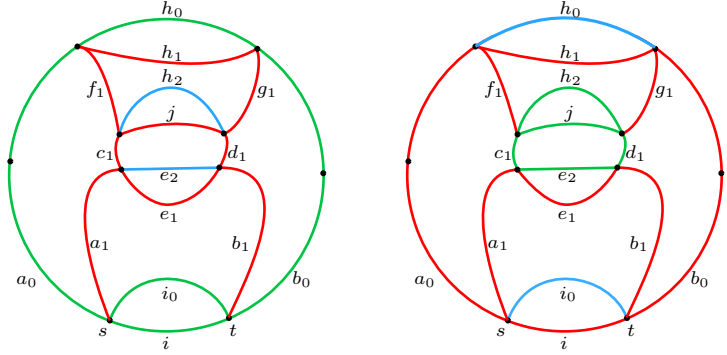
\caption{Proof of Theorem \ref{thm dim sym}. The left picture shows the partition $\sqcup K(i,d)$ relative to the segment $i$ and the right picture shows the partition $\sqcup K(j,d)$ relative to the segment $j.$}
\label{figdimsym1}
\end{center}
\end{figure}

\smallskip
(a)  Assume first that the segment $i$ lies in a bigon $(i,i_0)$ in $K$ and denote the two crossing points by $s$ and $t$.  We denote the four segments at $s$  by $a_0,a_1,i_0,i$ in clockwise order, and those at $t$ by $b_0,b_1,i_0,i$ in counterclockwise order, as shown in  the figure. For simplicity, we assume here without loss of generality that the path from $s$ to $A$ is given by a single segment $a_1$ and the path from $A$ to $t$ by a single segment $b_1$. If this is not the case, we would need to introduce two different letters for the first and the last segment in each of these paths. 
 This situation is illustrated in the  left picture of Figure~\ref{figdimsym1}, where $K(i,0)$ is drawn in green, $K(i,1) $ in red, and $K(i,2)$ in blue.

The right picture  of  the figure shows the sets $K(j,0)$ in green,  $K(j,1)$ in red and $K(j,2)$ in blue.  Note that $c_1,d_1,e_2,h_2$ lie in $K(j,0)$, because they bound a region with $j$. The segments $a_1,b_1,e_1,f_1,g_1,h_1,a_0,b_0$ lie in $K(j,1)$, because they see $c_1,d_1,e_2$ or $h_2$, but not $j$.

In this situation, the string $K\setminus K(j,0)$ is the enclosed sum $C\oplus[D,q]$ of the Hopf string $C$ consisting of the segments $f_1,h_1,h_0,a_0,i_0,b_0,g_1,$ with the point $p$ on the segment $i_0$, and the Hopf string $D$ consisting of the segments $a_0,i,i_0,a_1,e_1,b_1,b_0$. By definition, $K(j,1)$ is the boundary $\partial C\oplus[\partial D,q]$ consisting of the segments $f_1,h_1,a_0,b_0,g_1,i,a_1,e_1,b_1$ colored red in the right picture of the figure. In particular, $i\in K(j,1)$ as desired.

\smallskip (b) Now assume that $i$ is not contained in a bigon. This implies that there exists a crossing point $y\in K(i,0)$ on the on the strand $i_0$ such that $y\ne s,t$. This situation is illustrated in the left picture in Figure~\ref{figdimsym2}. Let $R$ denote the region incident to the segment $e_1$ but not to the segment $e_2$. 
Then $y$ is incident to two segments $i_0',i_0''$ in $K(i,0)$ and to two segments $\za_1,\zb_1$ in $K(i,1)$ as shown in the figure. However, the segments $\za_1,\zb_1$ cannot block the view from the segment $e_1$ to $K(i,0)$, since otherwise $e_1$ would not be in $K(i,1)$. 
Therefore we can draw a closed curve $\zg$ starting near $e_1$ inside the region $R$, following $b_1$, then crossing over the segment $i_0''$ to go around $y$, crossing back over $i_0'$ to enter back in $R$ and following $a_1$ to close near $e_1$. This curve $\zg$ is shown in dashed black ink in the center picture of the figure. Then $\zg$ crosses $K$ exactly twice, and thus $K$ is a connected sum  of the two links obtained by cutting along $\zg$ and connecting the loose ends on either side. This is a contradiction to our assumption that $K$ is a prime link.

\begin{figure}
\begin{center}
\scriptsize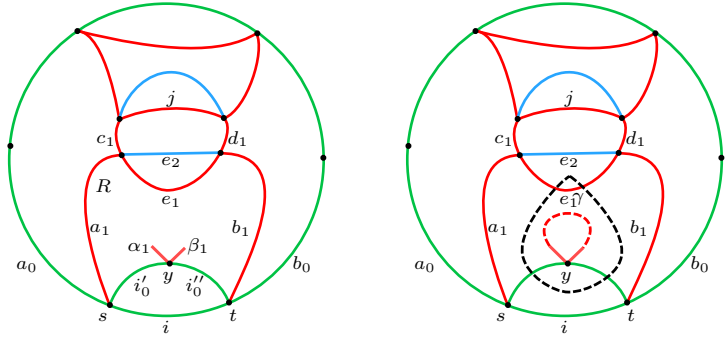
\caption{The case where $i$ is not contained in a bigon. 
Both  pictures show $K(i,r)$, where $r=0,1,2$. As in the previous figure, the colors are green for $r=0$, red for $r=1$ and blue for $r=2$. }
\label{figdimsym2}
\end{center}
\end{figure}

\subsubsection{The case $d>1$}
Assume now that $d>1$ and $x$ is an interior point in $K(i,r)$ with $r\le d$ that causes the dimension drop.

Thus $K\setminus \cup_{e<d} K(i,e)$ is an enclosed sum $A\oplus[B,p]$ with $p\in \partial A$ and $j$ a segment in $\partial B$ that is not incident to the unbounded component of the complement of $A\oplus[B,p]$. This situation is illustrated in the picture on the left in Figure~\ref{figdimsym3}. In this picture, we have drawn $K(i,r)$ in blue for $r=d+1$, red for $r=d$,  green for $r= d-1$ and orange for $r\le d-2$.
\begin{figure}
\begin{center}
\scriptsize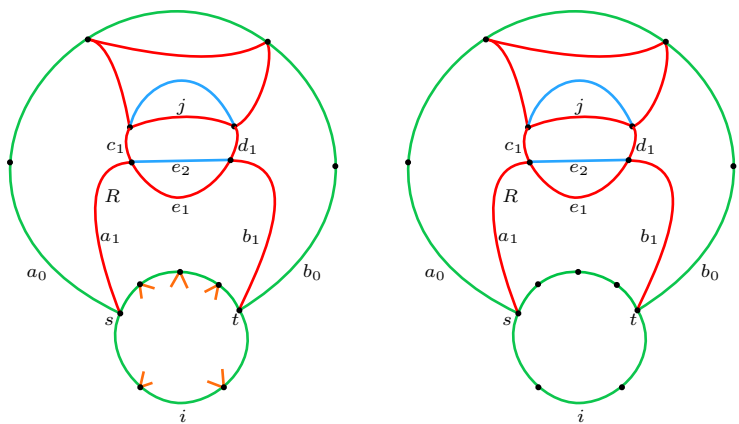
\caption{The case where $x$ is an interior point in $K(j,d)$. Both pictures show $K(i,r)$, where $r$ ranges over $d-2,d-1,d,d+1$ on the left, and over $d-1,d,d+1$ on the right. }
\label{figdimsym3}
\end{center}
\end{figure}

Let $K'$ be the diagram obtained from $K$ after smoothing away $K(i,0),K(i,1),\dots,K(i,d-1)$, as shown in the right picture of the figure. Then $K'$ will not contain any of the crossing points shared by $K(i,d)$ and $K(i,d-1)$, and thus in $K'$ the segment $i$ can be thought of lying in a bigon with crossing points $s$ and $t$. Now the result follows by the same argument as in the case $d=1$.
This completes the proof.
\end{proof}

\section{Applications}\label{sect appl}
In this section, we prove the results stated in section~\ref{sect intro appl}. 
\subsection{Newton polytopes and denominator vectors}
We start by proving Theorems~\ref{thm NewtonPolytop} and \ref{thm denominator}. 
Let $u$ be a cluster variable and denote its $F$-polynomial by $F_u=\sum c_\mathbf{e} \yy^{\mathbf{e}}$. The Newton polytope of $u$ is the convex hull of the lattice points $\mathbf{e}$ in $\RR^N$ for which  $c_\mathbf{e} $ is nonzero.

We use the following two results by Li and Pan. 
\begin{enumerate}
    \item \label{LP1} The Newton polytope is saturated, that is, the lattice points in the Newton polytope are precisely the points $\mathbf{e}$ for which the coefficient $c_\mathbf{e}$ is nonzero. Moreover, the coefficient $c_\mathbf{e}$ is equal to 1 if  and only if the point $\mathbf{e}$ is a vertex of the Newton polytope  \cite[Theorem 4.4]{LP}.
\item  \label{LP2} The $k$-th component of the denominator vector of the Laurent expansion of the cluster variable $u$ is equal to the maximum length of an edge of the Newton polytope that is parallel to the $k$-th standard basis vector of $\RR^n$  \cite[Theorem 5.4]{LP}.
\end{enumerate}

Now consider the cluster variables $x_{i;t}$ in the knot cluster $\xx_t$ of a prime link diagram. From \cite[Corollary 6.8]{BMS} we know that the $F$-polynomial is of the  form
$F_{i;t}=\sum_ {L\subset T(\zs(i))} \yy^{\underline{\dim} L}$, where the sum is over all submodules $L$ of $T(\zs(i))$. In particular, all coefficients in the $F$-polynomial are 0 or 1, and by statement (\ref{LP1}) above this means that every lattice point in the Newton polytope of $F_{i;t}$ is a vertex of the polytope. This proves Theorem~\ref{thm NewtonPolytop}.

Now let $e$ be an edge in the Newton polytope that is parallel to a standard basis vector. If the edge length of $e$ is greater than 1, then it contains a lattice point $p$ that is not an endpoint of $e$.  By (\ref{LP1}) the point $p$ corresponds to a monomial with coefficient 1 in the $F$-polynomial, but the fact that it is an interior point of the edge $e$ implies that it is not a vertex of the polytope. This is a contradiction. Therefore the length of every edge $e$
that is parallel to a standard basis vector is at most 1. Now (\ref{LP2}) implies that every entry in  the denominator vector is 0 or 1. This proves the first statement of Theorem~\ref{thm denominator}. 

Now let $j$ be a segment of $K$ that does not bound a common region with the segment $\zs(i)$. Then the representation $T(\zs(i))$ is nonzero at vertex $j$. It suffices to show that the exponent of $x_j$ in the denominator of $x_{i;t}$ is equal to 1. 
Consider the lattice of submodules of $T(\zs(i))$. By \cite[Theorem 6.9]{BMS} it is isomorphic to the lattice of Kauffman states of the diagram $K$ relative to the segment $i$. The covering relations in that lattice are given by transpositions at the segments  of $K$. Given two states $\cals,\cals'$ such that  $\cals'$ covers $\cals$ via a transposition at a segment $\ell$ then the dimension vectors of the corresponding submodules $L,L'$ of $T(i)$ satisfy $\underline{\dim}\, L' -\underline{\dim}\, L = e_\ell$, where $e_\ell$ is the standard basis vector of $\RR^{2n}$. Recall that the dimension vectors of submodules are precisely the lattice points in the Newton polytope. Thus the edge of the Newton polytope that joins  $\underline{\dim}\, L $ and $\underline{\dim}\, L'$ is parallel to the standard basis vector $e_\ell$. 
Since $T(\zs(i))$ is nonzero at vertex $j$, there exists such an edge for $\ell=j$ and therefore the exponent of $x_j$ in the denominator of $x_{i;t}$ is equal to 1. 
This completes the proof  of Theorem~\ref{thm denominator}. 

\subsection{The third Reidemeister move}
Now we prove Theorem \ref{thm Reidemeister3}. Let $K$ be a link diagram and let $K'$ be the diagram obtained from $K$ by a Reidemeister 3 move that moves a segment $i$ across (over or under) a crossing point formed by the segments $j$ and $k$. For an illustration, see the right picture in Figure~\ref{fig Reidemeister}. This Reidemeister 3 move is a special case of our RD3 moves introduced in section~\ref{sect mutation}. In particular, the quiver of the diagram $K'$ is obtained from the quiver of the diagram $K$ by the mutation sequence $\nu$ at vertices $i,j,k,i$. 
This shows that the cluster algebras $\cala(K)$ and $\cala(K')$ are the same, and that the initial seeds are related by the mutation sequence $\nu$.

\subsection{A new proof of a theorem of Murasugi}
As an application, we obtain a new proof of the following result.
\begin{thm}
 \label{thm Murasugi}
\cite{Murasugi} 
The Alexander polynomial of an alternating knot is an alternating sum \[\sum_{j=m}^M a_j t^j,\] with $a_j\ne 0$ for all $j=m,m+1,\ldots,M$.
\end{thm}
\begin{proof}
 Let $K$ be an oriented alternating diagram for the knot and choose a segment $i$ of $K$. Let $T(i)$ be the corresponding representation and $F_{T(i)}=\sum_{L\subset T(i)} \mathbf{y}^{\underline{\dim}\,L}$ the $F$-polynomial, where the sum is over all submodules $L$ of $T(i)$ and $\mathbf{y}^{\underline{\dim}\,L}=\prod_j y_j^{\dim L_j}$. Recall that the submodules of $T(i)$ form a poset that is isomorphic to the poset of Kauffman states. The unique minimal element is the zero module, and the unique maximal element is $T(i)$. 
 Furthermore, if $L\subset L'$ is a covering relation in the poset then $\dim\, L'=\dim \,L+1$, where $\dim L =\sum_j \dim L_j$ is the total dimension of the representation $L$. More precisely, there exists a unique $j$ such that 
 $\dim L'_j=\dim L_j+1$ and, for all $k\ne j$, $\dim L'_k=\dim L_k$. In particular
 
\begin{equation}
 \label{eq murasugi} 
\frac{\mathbf{y}^{\underline{\dim}\,L'}}{\mathbf{y}^{\underline{\dim}\,L}} = y_j.
\end{equation}

Since $K$ is alternating, the Alexander specialization sends $y_j$ to $-t$ or $-t^{-1}$. For $L=0$, we have ${\mathbf{y}^{\underline{\dim}\,L}}=1$, and if 
$0=L^0\subset L^1\subset L^2\subset L^3\subset \cdots $ is a chain in which each $L^k\subset L^{k+1}$ is a covering relation, then its contribution to the Alexander polynomial is the alternating sum
 $1-t^{\za_1}+t^{\za_2}-t^{\za_3}+\cdots$,
 where $\za_1=\pm 1$ and $|\za_k-\za_{k+1}|=1$ for all $k$. 
 In particular, the sign of $\pm t^{\za_k}$ is positive if and only if the power $\za_k$ is even.
 
 Together with equation (\ref{eq murasugi}) this implies that the contribution of an arbitrary submodule $L\subset T(i)$  to the Alexander polynomial is $(-1)^{\dim\,L}\,t^{\za(L)}$, with $\za(L)$ even if and only if $\dim \,L$ is even.
 In particular, all contributions of a fixed degree $j$ have the same sign. Thus, if there exists a submodule $L\subset T(i)$ whose contribution is of degree $j$, then the degree $j$ term in the Alexander polynomial has a nonzero coefficient.
 
 To show that the coefficients of the Alexander polynomial are nonzero, let $L^0$ and $L^\infty$ be submodules of $T(i)$ whose contributions to the Alexander polynomial have minimal and maximal degree, respectively. Since the poset is connected, there exists a path 
 $L^0 - L^1-L^2-\cdots-L^t=L^\infty$ from $L^0$ to $L^\infty$ in the Hasse diagram. From equation~\ref{eq murasugi} we see that the degrees of the contribution of $L^k$ and $L^{k+1}$ differ exactly by 1. Thus every intermediate degree is realized by this path. 
 This completes the proof.
\end{proof}

\end{document}

%% file: knot2112.pdf_tex
\begingroup%
  \makeatletter%
  \providecommand\color[2][]{%
    \errmessage{(Inkscape) Color is used for the text in Inkscape, but the package 'color.sty' is not loaded}%
    \renewcommand\color[2][]{}%
  }%
  \providecommand\transparent[1]{%
    \errmessage{(Inkscape) Transparency is used (non-zero) for the text in Inkscape, but the package 'transparent.sty' is not loaded}%
    \renewcommand\transparent[1]{}%
  }%
  \providecommand\rotatebox[2]{#2}%
  \newcommand*\fsize{\dimexpr\f@size pt\relax}%
  \newcommand*\lineheight[1]{\fontsize{\fsize}{#1\fsize}\selectfont}%
  \ifx\svgwidth\undefined%
    \setlength{\unitlength}{239.22085019bp}%
    \ifx\svgscale\undefined%
      \relax%
    \else%
      \setlength{\unitlength}{\unitlength * \real{\svgscale}}%
    \fi%
  \else%
    \setlength{\unitlength}{\svgwidth}%
  \fi%
  \global\let\svgwidth\undefined%
  \global\let\svgscale\undefined%
  \makeatother%
  \begin{picture}(1,0.61651023)%
    \lineheight{1}%
    \setlength\tabcolsep{0pt}%
    \put(0,0){\includegraphics[width=\unitlength,page=1]{knot2112.pdf}}%
    \put(0.12177319,0.48324992){\makebox(0,0)[lt]{\lineheight{1.25}\smash{\begin{tabular}[t]{l}1\end{tabular}}}}%
    \put(0.39648907,0.31390353){\makebox(0,0)[lt]{\lineheight{1.25}\smash{\begin{tabular}[t]{l}2\end{tabular}}}}%
    \put(0.67014461,0.13692077){\makebox(0,0)[lt]{\lineheight{1.25}\smash{\begin{tabular}[t]{l}3\end{tabular}}}}%
    \put(0.68614105,0.4833724){\makebox(0,0)[lt]{\lineheight{1.25}\smash{\begin{tabular}[t]{l}4\end{tabular}}}}%
    \put(0.53211912,0.31169124){\makebox(0,0)[lt]{\lineheight{1.25}\smash{\begin{tabular}[t]{l}5\end{tabular}}}}%
    \put(0.32553884,0.13692077){\makebox(0,0)[lt]{\lineheight{1.25}\smash{\begin{tabular}[t]{l}6\end{tabular}}}}%
    \put(0.12301668,0.29560318){\makebox(0,0)[lt]{\lineheight{1.25}\smash{\begin{tabular}[t]{l}7\end{tabular}}}}%
    \put(0.46347271,0.00083277){\makebox(0,0)[lt]{\lineheight{1.25}\smash{\begin{tabular}[t]{l}8\end{tabular}}}}%
    \put(0.68923044,0.30597045){\makebox(0,0)[lt]{\lineheight{1.25}\smash{\begin{tabular}[t]{l}9\end{tabular}}}}%
    \put(0.45471537,0.48263759){\makebox(0,0)[lt]{\lineheight{1.25}\smash{\begin{tabular}[t]{l}10\end{tabular}}}}%
    \put(0.23665467,0.28306246){\makebox(0,0)[lt]{\lineheight{1.25}\smash{\begin{tabular}[t]{l}11\end{tabular}}}}%
    \put(0.12292512,0.13616146){\makebox(0,0)[lt]{\lineheight{1.25}\smash{\begin{tabular}[t]{l}12\end{tabular}}}}%
    \put(0,0){\includegraphics[width=\unitlength,page=2]{knot2112.pdf}}%
  \end{picture}%
\endgroup%

%% file: figReidemeister.pdf_tex
\begingroup%
  \makeatletter%
  \providecommand\color[2][]{%
    \errmessage{(Inkscape) Color is used for the text in Inkscape, but the package 'color.sty' is not loaded}%
    \renewcommand\color[2][]{}%
  }%
  \providecommand\transparent[1]{%
    \errmessage{(Inkscape) Transparency is used (non-zero) for the text in Inkscape, but the package 'transparent.sty' is not loaded}%
    \renewcommand\transparent[1]{}%
  }%
  \providecommand\rotatebox[2]{#2}%
  \newcommand*\fsize{\dimexpr\f@size pt\relax}%
  \newcommand*\lineheight[1]{\fontsize{\fsize}{#1\fsize}\selectfont}%
  \ifx\svgwidth\undefined%
    \setlength{\unitlength}{329.53154701bp}%
    \ifx\svgscale\undefined%
      \relax%
    \else%
      \setlength{\unitlength}{\unitlength * \real{\svgscale}}%
    \fi%
  \else%
    \setlength{\unitlength}{\svgwidth}%
  \fi%
  \global\let\svgwidth\undefined%
  \global\let\svgscale\undefined%
  \makeatother%
  \begin{picture}(1,0.08173924)%
    \lineheight{1}%
    \setlength\tabcolsep{0pt}%
    \put(0,0){\includegraphics[width=\unitlength,page=1]{figReidemeister.pdf}}%
  \end{picture}%
\endgroup%

%% file: figure8knot.pdf_tex
\begingroup%
  \makeatletter%
  \providecommand\color[2][]{%
    \errmessage{(Inkscape) Color is used for the text in Inkscape, but the package 'color.sty' is not loaded}%
    \renewcommand\color[2][]{}%
  }%
  \providecommand\transparent[1]{%
    \errmessage{(Inkscape) Transparency is used (non-zero) for the text in Inkscape, but the package 'transparent.sty' is not loaded}%
    \renewcommand\transparent[1]{}%
  }%
  \providecommand\rotatebox[2]{#2}%
  \newcommand*\fsize{\dimexpr\f@size pt\relax}%
  \newcommand*\lineheight[1]{\fontsize{\fsize}{#1\fsize}\selectfont}%
  \ifx\svgwidth\undefined%
    \setlength{\unitlength}{408.36822556bp}%
    \ifx\svgscale\undefined%
      \relax%
    \else%
      \setlength{\unitlength}{\unitlength * \real{\svgscale}}%
    \fi%
  \else%
    \setlength{\unitlength}{\svgwidth}%
  \fi%
  \global\let\svgwidth\undefined%
  \global\let\svgscale\undefined%
  \makeatother%
  \begin{picture}(1,0.17979329)%
    \lineheight{1}%
    \setlength\tabcolsep{0pt}%
    \put(0,0){\includegraphics[width=\unitlength,page=1]{figure8knot.pdf}}%
    \put(0.14401821,0){\makebox(0,0)[lt]{\lineheight{1.25}\smash{\begin{tabular}[t]{l}1\end{tabular}}}}%
    \put(0.19175726,0.11203126){\makebox(0,0)[lt]{\lineheight{1.25}\smash{\begin{tabular}[t]{l}2\end{tabular}}}}%
    \put(0.07422826,0.16161887){\makebox(0,0)[lt]{\lineheight{1.25}\smash{\begin{tabular}[t]{l}3\end{tabular}}}}%
    \put(0.03868401,0.05509735){\makebox(0,0)[lt]{\lineheight{1.25}\smash{\begin{tabular}[t]{l}4\end{tabular}}}}%
    \put(0.12932557,0.10101179){\makebox(0,0)[lt]{\lineheight{1.25}\smash{\begin{tabular}[t]{l}5\end{tabular}}}}%
    \put(0.19176923,0.16161887){\makebox(0,0)[lt]{\lineheight{1.25}\smash{\begin{tabular}[t]{l}6\end{tabular}}}}%
    \put(0.17340346,0.05326078){\makebox(0,0)[lt]{\lineheight{1.25}\smash{\begin{tabular}[t]{l}7\end{tabular}}}}%
    \put(0.03876067,0.10042426){\makebox(0,0)[lt]{\lineheight{1.25}\smash{\begin{tabular}[t]{l}8\end{tabular}}}}%
    \put(0,0){\includegraphics[width=\unitlength,page=2]{figure8knot.pdf}}%
    \put(0.52969949,0){\makebox(0,0)[lt]{\lineheight{1.25}\smash{\begin{tabular}[t]{l}1\end{tabular}}}}%
    \put(0.57743843,0.11203126){\makebox(0,0)[lt]{\lineheight{1.25}\smash{\begin{tabular}[t]{l}2\end{tabular}}}}%
    \put(0.45990957,0.16161887){\makebox(0,0)[lt]{\lineheight{1.25}\smash{\begin{tabular}[t]{l}3\end{tabular}}}}%
    \put(0.4929679,0.10468495){\makebox(0,0)[lt]{\lineheight{1.25}\smash{\begin{tabular}[t]{l}5\end{tabular}}}}%
    \put(0.57745045,0.16161887){\makebox(0,0)[lt]{\lineheight{1.25}\smash{\begin{tabular}[t]{l}6\end{tabular}}}}%
    \put(0.55908471,0.05326078){\makebox(0,0)[lt]{\lineheight{1.25}\smash{\begin{tabular}[t]{l}7\end{tabular}}}}%
    \put(0,0){\includegraphics[width=\unitlength,page=3]{figure8knot.pdf}}%
    \put(0.3258394,0.08999233){\makebox(0,0)[lt]{\lineheight{1.25}\smash{\begin{tabular}[t]{l}$(4,8)$\end{tabular}}}}%
    \put(0.90436123,0.00734632){\makebox(0,0)[lt]{\lineheight{1.25}\smash{\begin{tabular}[t]{l}1\end{tabular}}}}%
    \put(0.83457126,0.1432531){\makebox(0,0)[lt]{\lineheight{1.25}\smash{\begin{tabular}[t]{l}3\end{tabular}}}}%
    \put(0.86762964,0.10468495){\makebox(0,0)[lt]{\lineheight{1.25}\smash{\begin{tabular}[t]{l}5\end{tabular}}}}%
    \put(0.9300734,0.07162657){\makebox(0,0)[lt]{\lineheight{1.25}\smash{\begin{tabular}[t]{l}7\end{tabular}}}}%
    \put(0,0){\includegraphics[width=\unitlength,page=4]{figure8knot.pdf}}%
    \put(0.70050104,0.08999233){\makebox(0,0)[lt]{\lineheight{1.25}\smash{\begin{tabular}[t]{l}$(2,6)$\end{tabular}}}}%
    \put(0,0){\includegraphics[width=\unitlength,page=5]{figure8knot.pdf}}%
  \end{picture}%
\endgroup%

%% file: BorromeanRings.pdf_tex
\begingroup%
  \makeatletter%
  \providecommand\color[2][]{%
    \errmessage{(Inkscape) Color is used for the text in Inkscape, but the package 'color.sty' is not loaded}%
    \renewcommand\color[2][]{}%
  }%
  \providecommand\transparent[1]{%
    \errmessage{(Inkscape) Transparency is used (non-zero) for the text in Inkscape, but the package 'transparent.sty' is not loaded}%
    \renewcommand\transparent[1]{}%
  }%
  \providecommand\rotatebox[2]{#2}%
  \newcommand*\fsize{\dimexpr\f@size pt\relax}%
  \newcommand*\lineheight[1]{\fontsize{\fsize}{#1\fsize}\selectfont}%
  \ifx\svgwidth\undefined%
    \setlength{\unitlength}{466.73061017bp}%
    \ifx\svgscale\undefined%
      \relax%
    \else%
      \setlength{\unitlength}{\unitlength * \real{\svgscale}}%
    \fi%
  \else%
    \setlength{\unitlength}{\svgwidth}%
  \fi%
  \global\let\svgwidth\undefined%
  \global\let\svgscale\undefined%
  \makeatother%
  \begin{picture}(1,0.2663121)%
    \lineheight{1}%
    \setlength\tabcolsep{0pt}%
    \put(0,0){\includegraphics[width=\unitlength,page=1]{BorromeanRings.pdf}}%
    \put(0.10655565,0.1241756){\makebox(0,0)[lt]{\lineheight{1.25}\smash{\begin{tabular}[t]{l}3\end{tabular}}}}%
    \put(0.08341764,0.20075264){\makebox(0,0)[lt]{\lineheight{1.25}\smash{\begin{tabular}[t]{l}5\end{tabular}}}}%
    \put(0.16429554,0.07118147){\makebox(0,0)[lt]{\lineheight{1.25}\smash{\begin{tabular}[t]{l}8\end{tabular}}}}%
    \put(0.2671029,0.20582727){\makebox(0,0)[lt]{\lineheight{1.25}\smash{\begin{tabular}[t]{l}10\end{tabular}}}}%
    \put(0,0){\includegraphics[width=\unitlength,page=2]{BorromeanRings.pdf}}%
    \put(0.16432291,0.12401343){\makebox(0,0)[lt]{\lineheight{1.25}\smash{\begin{tabular}[t]{l}2\end{tabular}}}}%
    \put(0.18774468,0.20075265){\makebox(0,0)[lt]{\lineheight{1.25}\smash{\begin{tabular}[t]{l}4\end{tabular}}}}%
    \put(0.10216053,0.06866341){\makebox(0,0)[lt]{\lineheight{1.25}\smash{\begin{tabular}[t]{l}7\end{tabular}}}}%
    \put(-0.00290836,0.2057959){\makebox(0,0)[lt]{\lineheight{1.25}\smash{\begin{tabular}[t]{l}11\end{tabular}}}}%
    \put(0.13396055,0.17980275){\makebox(0,0)[lt]{\lineheight{1.25}\smash{\begin{tabular}[t]{l}1\end{tabular}}}}%
    \put(0,0){\includegraphics[width=\unitlength,page=3]{BorromeanRings.pdf}}%
    \put(0.06123186,0.15923752){\makebox(0,0)[lt]{\lineheight{1.25}\smash{\begin{tabular}[t]{l}6\end{tabular}}}}%
    \put(0.21111176,0.16509769){\makebox(0,0)[lt]{\lineheight{1.25}\smash{\begin{tabular}[t]{l}9\end{tabular}}}}%
    \put(0.13262374,0.00000001){\makebox(0,0)[lt]{\lineheight{1.25}\smash{\begin{tabular}[t]{l}12\end{tabular}}}}%
    \put(0,0){\includegraphics[width=\unitlength,page=4]{BorromeanRings.pdf}}%
    \put(0.47614791,0.0754814){\makebox(0,0)[lt]{\lineheight{1.25}\smash{\begin{tabular}[t]{l}2\end{tabular}}}}%
    \put(0.43694066,0.20075264){\makebox(0,0)[lt]{\lineheight{1.25}\smash{\begin{tabular}[t]{l}5\end{tabular}}}}%
    \put(0.54352928,0.07118147){\makebox(0,0)[lt]{\lineheight{1.25}\smash{\begin{tabular}[t]{l}8\end{tabular}}}}%
    \put(0.62062592,0.20582727){\makebox(0,0)[lt]{\lineheight{1.25}\smash{\begin{tabular}[t]{l}10\end{tabular}}}}%
    \put(0,0){\includegraphics[width=\unitlength,page=5]{BorromeanRings.pdf}}%
    \put(0.50177668,0.0758057){\makebox(0,0)[lt]{\lineheight{1.25}\smash{\begin{tabular}[t]{l}3\end{tabular}}}}%
    \put(0.54126769,0.20075265){\makebox(0,0)[lt]{\lineheight{1.25}\smash{\begin{tabular}[t]{l}4\end{tabular}}}}%
    \put(0.42997273,0.07187726){\makebox(0,0)[lt]{\lineheight{1.25}\smash{\begin{tabular}[t]{l}7\end{tabular}}}}%
    \put(0.35061465,0.20579591){\makebox(0,0)[lt]{\lineheight{1.25}\smash{\begin{tabular}[t]{l}11\end{tabular}}}}%
    \put(0,0){\includegraphics[width=\unitlength,page=6]{BorromeanRings.pdf}}%
    \put(0.48748356,0.04160726){\makebox(0,0)[lt]{\lineheight{1.25}\smash{\begin{tabular}[t]{l}1\end{tabular}}}}%
    \put(0.43082412,0.10593574){\makebox(0,0)[lt]{\lineheight{1.25}\smash{\begin{tabular}[t]{l}6\end{tabular}}}}%
    \put(0.54366528,0.10591482){\makebox(0,0)[lt]{\lineheight{1.25}\smash{\begin{tabular}[t]{l}9\end{tabular}}}}%
    \put(0.48293291,0){\makebox(0,0)[lt]{\lineheight{1.25}\smash{\begin{tabular}[t]{l}12\end{tabular}}}}%
    \put(0.81681571,0.137031){\makebox(0,0)[lt]{\lineheight{1.25}\smash{\begin{tabular}[t]{l}2\end{tabular}}}}%
    \put(0.97414911,0.20582727){\makebox(0,0)[lt]{\lineheight{1.25}\smash{\begin{tabular}[t]{l}10\end{tabular}}}}%
    \put(0.864941,0.13686882){\makebox(0,0)[lt]{\lineheight{1.25}\smash{\begin{tabular}[t]{l}3\end{tabular}}}}%
    \put(0.70413746,0.20579591){\makebox(0,0)[lt]{\lineheight{1.25}\smash{\begin{tabular}[t]{l}11\end{tabular}}}}%
    \put(0.84100666,0.04160726){\makebox(0,0)[lt]{\lineheight{1.25}\smash{\begin{tabular}[t]{l}1\end{tabular}}}}%
    \put(0.83645596,0){\makebox(0,0)[lt]{\lineheight{1.25}\smash{\begin{tabular}[t]{l}12\end{tabular}}}}%
    \put(0,0){\includegraphics[width=\unitlength,page=7]{BorromeanRings.pdf}}%
  \end{picture}%
\endgroup%

%% file: figR3.pdf_tex
\begingroup%
  \makeatletter%
  \providecommand\color[2][]{%
    \errmessage{(Inkscape) Color is used for the text in Inkscape, but the package 'color.sty' is not loaded}%
    \renewcommand\color[2][]{}%
  }%
  \providecommand\transparent[1]{%
    \errmessage{(Inkscape) Transparency is used (non-zero) for the text in Inkscape, but the package 'transparent.sty' is not loaded}%
    \renewcommand\transparent[1]{}%
  }%
  \providecommand\rotatebox[2]{#2}%
  \newcommand*\fsize{\dimexpr\f@size pt\relax}%
  \newcommand*\lineheight[1]{\fontsize{\fsize}{#1\fsize}\selectfont}%
  \ifx\svgwidth\undefined%
    \setlength{\unitlength}{348.83784737bp}%
    \ifx\svgscale\undefined%
      \relax%
    \else%
      \setlength{\unitlength}{\unitlength * \real{\svgscale}}%
    \fi%
  \else%
    \setlength{\unitlength}{\svgwidth}%
  \fi%
  \global\let\svgwidth\undefined%
  \global\let\svgscale\undefined%
  \makeatother%
  \begin{picture}(1,0.30734226)%
    \lineheight{1}%
    \setlength\tabcolsep{0pt}%
    \put(0,0){\includegraphics[width=\unitlength,page=1]{figR3.pdf}}%
    \put(0.19349964,0.2221319){\makebox(0,0)[lt]{\lineheight{1.25}\smash{\begin{tabular}[t]{l}1\end{tabular}}}}%
    \put(0.23219957,0.16408204){\makebox(0,0)[lt]{\lineheight{1.25}\smash{\begin{tabular}[t]{l}2\end{tabular}}}}%
    \put(0.15049974,0.16408204){\makebox(0,0)[lt]{\lineheight{1.25}\smash{\begin{tabular}[t]{l}3\end{tabular}}}}%
    \put(0.29024943,0.28233179){\makebox(0,0)[lt]{\lineheight{1.25}\smash{\begin{tabular}[t]{l}4\end{tabular}}}}%
    \put(0.09674988,0.28018182){\makebox(0,0)[lt]{\lineheight{1.25}\smash{\begin{tabular}[t]{l}5\end{tabular}}}}%
    \put(0.00215005,0.18988202){\makebox(0,0)[lt]{\lineheight{1.25}\smash{\begin{tabular}[t]{l}6\end{tabular}}}}%
    \put(0.09459988,0.0006825){\makebox(0,0)[lt]{\lineheight{1.25}\smash{\begin{tabular}[t]{l}7\end{tabular}}}}%
    \put(0.28164944,0.00283246){\makebox(0,0)[lt]{\lineheight{1.25}\smash{\begin{tabular}[t]{l}8\end{tabular}}}}%
    \put(0.37839926,0.18988196){\makebox(0,0)[lt]{\lineheight{1.25}\smash{\begin{tabular}[t]{l}9\end{tabular}}}}%
    \put(0.79549834,0.05873229){\makebox(0,0)[lt]{\lineheight{1.25}\smash{\begin{tabular}[t]{l}1\end{tabular}}}}%
    \put(0.83419825,0.12108214){\makebox(0,0)[lt]{\lineheight{1.25}\smash{\begin{tabular}[t]{l}3\end{tabular}}}}%
    \put(0.7567985,0.12108214){\makebox(0,0)[lt]{\lineheight{1.25}\smash{\begin{tabular}[t]{l}2\end{tabular}}}}%
    \put(0.89224829,0.28233179){\makebox(0,0)[lt]{\lineheight{1.25}\smash{\begin{tabular}[t]{l}4\end{tabular}}}}%
    \put(0.69874864,0.28018182){\makebox(0,0)[lt]{\lineheight{1.25}\smash{\begin{tabular}[t]{l}5\end{tabular}}}}%
    \put(0.60414891,0.06088232){\makebox(0,0)[lt]{\lineheight{1.25}\smash{\begin{tabular}[t]{l}6\end{tabular}}}}%
    \put(0.69659855,0.0006825){\makebox(0,0)[lt]{\lineheight{1.25}\smash{\begin{tabular}[t]{l}7\end{tabular}}}}%
    \put(0.88364817,0.00283246){\makebox(0,0)[lt]{\lineheight{1.25}\smash{\begin{tabular}[t]{l}8\end{tabular}}}}%
    \put(0.98039812,0.06088226){\makebox(0,0)[lt]{\lineheight{1.25}\smash{\begin{tabular}[t]{l}9\end{tabular}}}}%
  \end{picture}%
\endgroup%

%% file: figtransposition.pdf_tex
\begingroup%
  \makeatletter%
  \providecommand\color[2][]{%
    \errmessage{(Inkscape) Color is used for the text in Inkscape, but the package 'color.sty' is not loaded}%
    \renewcommand\color[2][]{}%
  }%
  \providecommand\transparent[1]{%
    \errmessage{(Inkscape) Transparency is used (non-zero) for the text in Inkscape, but the package 'transparent.sty' is not loaded}%
    \renewcommand\transparent[1]{}%
  }%
  \providecommand\rotatebox[2]{#2}%
  \newcommand*\fsize{\dimexpr\f@size pt\relax}%
  \newcommand*\lineheight[1]{\fontsize{\fsize}{#1\fsize}\selectfont}%
  \ifx\svgwidth\undefined%
    \setlength{\unitlength}{300.00007559bp}%
    \ifx\svgscale\undefined%
      \relax%
    \else%
      \setlength{\unitlength}{\unitlength * \real{\svgscale}}%
    \fi%
  \else%
    \setlength{\unitlength}{\svgwidth}%
  \fi%
  \global\let\svgwidth\undefined%
  \global\let\svgscale\undefined%
  \makeatother%
  \begin{picture}(1,0.19999994)%
    \lineheight{1}%
    \setlength\tabcolsep{0pt}%
    \put(0,0){\includegraphics[width=\unitlength,page=1]{figtransposition.pdf}}%
    \put(0.19499992,0.06999994){\makebox(0,0)[lt]{\lineheight{1.25}\smash{\begin{tabular}[t]{l}$j$\end{tabular}}}}%
    \put(0,0){\includegraphics[width=\unitlength,page=2]{figtransposition.pdf}}%
    \put(0.79499963,0.06999994){\makebox(0,0)[lt]{\lineheight{1.25}\smash{\begin{tabular}[t]{l}$j$\end{tabular}}}}%
    \put(0,0){\includegraphics[width=\unitlength,page=3]{figtransposition.pdf}}%
  \end{picture}%
\endgroup%

%% file: figextranspositions.pdf_tex
\begingroup%
  \makeatletter%
  \providecommand\color[2][]{%
    \errmessage{(Inkscape) Color is used for the text in Inkscape, but the package 'color.sty' is not loaded}%
    \renewcommand\color[2][]{}%
  }%
  \providecommand\transparent[1]{%
    \errmessage{(Inkscape) Transparency is used (non-zero) for the text in Inkscape, but the package 'transparent.sty' is not loaded}%
    \renewcommand\transparent[1]{}%
  }%
  \providecommand\rotatebox[2]{#2}%
  \newcommand*\fsize{\dimexpr\f@size pt\relax}%
  \newcommand*\lineheight[1]{\fontsize{\fsize}{#1\fsize}\selectfont}%
  \ifx\svgwidth\undefined%
    \setlength{\unitlength}{407.7649564bp}%
    \ifx\svgscale\undefined%
      \relax%
    \else%
      \setlength{\unitlength}{\unitlength * \real{\svgscale}}%
    \fi%
  \else%
    \setlength{\unitlength}{\svgwidth}%
  \fi%
  \global\let\svgwidth\undefined%
  \global\let\svgscale\undefined%
  \makeatother%
  \begin{picture}(1,0.24878015)%
    \lineheight{1}%
    \setlength\tabcolsep{0pt}%
    \put(0,0){\includegraphics[width=\unitlength,page=1]{figextranspositions.pdf}}%
    \put(0.07101462,0.06242942){\makebox(0,0)[lt]{\lineheight{1.25}\smash{\begin{tabular}[t]{l}1\end{tabular}}}}%
    \put(0.09308613,0.03575966){\makebox(0,0)[lt]{\lineheight{1.25}\smash{\begin{tabular}[t]{l}2\end{tabular}}}}%
    \put(0.04894304,0.03575966){\makebox(0,0)[lt]{\lineheight{1.25}\smash{\begin{tabular}[t]{l}3\end{tabular}}}}%
    \put(0.11975575,0.08817957){\makebox(0,0)[lt]{\lineheight{1.25}\smash{\begin{tabular}[t]{l}4\end{tabular}}}}%
    \put(0.01859472,0.08909921){\makebox(0,0)[lt]{\lineheight{1.25}\smash{\begin{tabular}[t]{l}5\end{tabular}}}}%
    \put(-0.00098431,0.05415263){\makebox(0,0)[lt]{\lineheight{1.25}\smash{\begin{tabular}[t]{l}6\end{tabular}}}}%
    \put(0.04407844,0.00024428){\makebox(0,0)[lt]{\lineheight{1.25}\smash{\begin{tabular}[t]{l}7\end{tabular}}}}%
    \put(0.10445844,0.00024428){\makebox(0,0)[lt]{\lineheight{1.25}\smash{\begin{tabular}[t]{l}8\end{tabular}}}}%
    \put(0.14002205,0.05413825){\makebox(0,0)[lt]{\lineheight{1.25}\smash{\begin{tabular}[t]{l}9\end{tabular}}}}%
    \put(0,0){\includegraphics[width=\unitlength,page=2]{figextranspositions.pdf}}%
    \put(0.16114008,0.21141231){\makebox(0,0)[lt]{\lineheight{1.25}\smash{\begin{tabular}[t]{l}1\end{tabular}}}}%
    \put(0.18321158,0.18474255){\makebox(0,0)[lt]{\lineheight{1.25}\smash{\begin{tabular}[t]{l}2\end{tabular}}}}%
    \put(0.13906849,0.18474255){\makebox(0,0)[lt]{\lineheight{1.25}\smash{\begin{tabular}[t]{l}3\end{tabular}}}}%
    \put(0.20988121,0.23716246){\makebox(0,0)[lt]{\lineheight{1.25}\smash{\begin{tabular}[t]{l}4\end{tabular}}}}%
    \put(0.10872016,0.2380821){\makebox(0,0)[lt]{\lineheight{1.25}\smash{\begin{tabular}[t]{l}5\end{tabular}}}}%
    \put(0.08914113,0.20313551){\makebox(0,0)[lt]{\lineheight{1.25}\smash{\begin{tabular}[t]{l}6\end{tabular}}}}%
    \put(0.13420388,0.14922717){\makebox(0,0)[lt]{\lineheight{1.25}\smash{\begin{tabular}[t]{l}7\end{tabular}}}}%
    \put(0.19458391,0.14922717){\makebox(0,0)[lt]{\lineheight{1.25}\smash{\begin{tabular}[t]{l}8\end{tabular}}}}%
    \put(0.23014753,0.20312114){\makebox(0,0)[lt]{\lineheight{1.25}\smash{\begin{tabular}[t]{l}9\end{tabular}}}}%
    \put(0,0){\includegraphics[width=\unitlength,page=3]{figextranspositions.pdf}}%
    \put(0.0859024,0.11600251){\makebox(0,0)[lt]{\lineheight{1.25}\smash{\begin{tabular}[t]{l}3\end{tabular}}}}%
    \put(0,0){\includegraphics[width=\unitlength,page=4]{figextranspositions.pdf}}%
    \put(0.38185544,0.21141231){\makebox(0,0)[lt]{\lineheight{1.25}\smash{\begin{tabular}[t]{l}1\end{tabular}}}}%
    \put(0.40392694,0.18474255){\makebox(0,0)[lt]{\lineheight{1.25}\smash{\begin{tabular}[t]{l}2\end{tabular}}}}%
    \put(0.35978388,0.18474255){\makebox(0,0)[lt]{\lineheight{1.25}\smash{\begin{tabular}[t]{l}3\end{tabular}}}}%
    \put(0.43059655,0.23716246){\makebox(0,0)[lt]{\lineheight{1.25}\smash{\begin{tabular}[t]{l}4\end{tabular}}}}%
    \put(0.32943558,0.2380821){\makebox(0,0)[lt]{\lineheight{1.25}\smash{\begin{tabular}[t]{l}5\end{tabular}}}}%
    \put(0.30985654,0.20313551){\makebox(0,0)[lt]{\lineheight{1.25}\smash{\begin{tabular}[t]{l}6\end{tabular}}}}%
    \put(0.35491927,0.14922717){\makebox(0,0)[lt]{\lineheight{1.25}\smash{\begin{tabular}[t]{l}7\end{tabular}}}}%
    \put(0.41529924,0.14922717){\makebox(0,0)[lt]{\lineheight{1.25}\smash{\begin{tabular}[t]{l}8\end{tabular}}}}%
    \put(0.45086283,0.20312114){\makebox(0,0)[lt]{\lineheight{1.25}\smash{\begin{tabular}[t]{l}9\end{tabular}}}}%
    \put(0,0){\includegraphics[width=\unitlength,page=5]{figextranspositions.pdf}}%
    \put(0.2735105,0.18221708){\makebox(0,0)[lt]{\lineheight{1.25}\smash{\begin{tabular}[t]{l}1\end{tabular}}}}%
    \put(0,0){\includegraphics[width=\unitlength,page=6]{figextranspositions.pdf}}%
    \put(0.38185544,0.06426879){\makebox(0,0)[lt]{\lineheight{1.25}\smash{\begin{tabular}[t]{l}1\end{tabular}}}}%
    \put(0.40392694,0.03759902){\makebox(0,0)[lt]{\lineheight{1.25}\smash{\begin{tabular}[t]{l}2\end{tabular}}}}%
    \put(0.35978388,0.03759902){\makebox(0,0)[lt]{\lineheight{1.25}\smash{\begin{tabular}[t]{l}3\end{tabular}}}}%
    \put(0.43059655,0.09001894){\makebox(0,0)[lt]{\lineheight{1.25}\smash{\begin{tabular}[t]{l}4\end{tabular}}}}%
    \put(0.32943558,0.09093855){\makebox(0,0)[lt]{\lineheight{1.25}\smash{\begin{tabular}[t]{l}5\end{tabular}}}}%
    \put(0.30985654,0.05599199){\makebox(0,0)[lt]{\lineheight{1.25}\smash{\begin{tabular}[t]{l}6\end{tabular}}}}%
    \put(0.35491927,0.00208365){\makebox(0,0)[lt]{\lineheight{1.25}\smash{\begin{tabular}[t]{l}7\end{tabular}}}}%
    \put(0.41529924,0.00208365){\makebox(0,0)[lt]{\lineheight{1.25}\smash{\begin{tabular}[t]{l}8\end{tabular}}}}%
    \put(0.45086283,0.05597762){\makebox(0,0)[lt]{\lineheight{1.25}\smash{\begin{tabular}[t]{l}9\end{tabular}}}}%
    \put(0,0){\includegraphics[width=\unitlength,page=7]{figextranspositions.pdf}}%
    \put(0.2735105,0.11232402){\makebox(0,0)[lt]{\lineheight{1.25}\smash{\begin{tabular}[t]{l}7\end{tabular}}}}%
    \put(0,0){\includegraphics[width=\unitlength,page=8]{figextranspositions.pdf}}%
    \put(0.60257072,0.21141231){\makebox(0,0)[lt]{\lineheight{1.25}\smash{\begin{tabular}[t]{l}1\end{tabular}}}}%
    \put(0.62464228,0.18474255){\makebox(0,0)[lt]{\lineheight{1.25}\smash{\begin{tabular}[t]{l}2\end{tabular}}}}%
    \put(0.58049916,0.18474255){\makebox(0,0)[lt]{\lineheight{1.25}\smash{\begin{tabular}[t]{l}3\end{tabular}}}}%
    \put(0.65131188,0.23716246){\makebox(0,0)[lt]{\lineheight{1.25}\smash{\begin{tabular}[t]{l}4\end{tabular}}}}%
    \put(0.55015086,0.2380821){\makebox(0,0)[lt]{\lineheight{1.25}\smash{\begin{tabular}[t]{l}5\end{tabular}}}}%
    \put(0.53057182,0.20313551){\makebox(0,0)[lt]{\lineheight{1.25}\smash{\begin{tabular}[t]{l}6\end{tabular}}}}%
    \put(0.57563455,0.14922717){\makebox(0,0)[lt]{\lineheight{1.25}\smash{\begin{tabular}[t]{l}7\end{tabular}}}}%
    \put(0.63601458,0.14922717){\makebox(0,0)[lt]{\lineheight{1.25}\smash{\begin{tabular}[t]{l}8\end{tabular}}}}%
    \put(0.67157811,0.20312114){\makebox(0,0)[lt]{\lineheight{1.25}\smash{\begin{tabular}[t]{l}9\end{tabular}}}}%
    \put(0,0){\includegraphics[width=\unitlength,page=9]{figextranspositions.pdf}}%
    \put(0.49422578,0.18221708){\makebox(0,0)[lt]{\lineheight{1.25}\smash{\begin{tabular}[t]{l}7\end{tabular}}}}%
    \put(0,0){\includegraphics[width=\unitlength,page=10]{figextranspositions.pdf}}%
    \put(0.60257072,0.06426879){\makebox(0,0)[lt]{\lineheight{1.25}\smash{\begin{tabular}[t]{l}1\end{tabular}}}}%
    \put(0.62464228,0.03759902){\makebox(0,0)[lt]{\lineheight{1.25}\smash{\begin{tabular}[t]{l}2\end{tabular}}}}%
    \put(0.58049916,0.03759902){\makebox(0,0)[lt]{\lineheight{1.25}\smash{\begin{tabular}[t]{l}3\end{tabular}}}}%
    \put(0.65131188,0.09001894){\makebox(0,0)[lt]{\lineheight{1.25}\smash{\begin{tabular}[t]{l}4\end{tabular}}}}%
    \put(0.55015086,0.09093855){\makebox(0,0)[lt]{\lineheight{1.25}\smash{\begin{tabular}[t]{l}5\end{tabular}}}}%
    \put(0.53057182,0.05599199){\makebox(0,0)[lt]{\lineheight{1.25}\smash{\begin{tabular}[t]{l}6\end{tabular}}}}%
    \put(0.57563455,0.00208365){\makebox(0,0)[lt]{\lineheight{1.25}\smash{\begin{tabular}[t]{l}7\end{tabular}}}}%
    \put(0.63601458,0.00208365){\makebox(0,0)[lt]{\lineheight{1.25}\smash{\begin{tabular}[t]{l}8\end{tabular}}}}%
    \put(0.67157811,0.05597762){\makebox(0,0)[lt]{\lineheight{1.25}\smash{\begin{tabular}[t]{l}9\end{tabular}}}}%
    \put(0,0){\includegraphics[width=\unitlength,page=11]{figextranspositions.pdf}}%
    \put(0.49422578,0.03507356){\makebox(0,0)[lt]{\lineheight{1.25}\smash{\begin{tabular}[t]{l}8\end{tabular}}}}%
    \put(0,0){\includegraphics[width=\unitlength,page=12]{figextranspositions.pdf}}%
    \put(0.49422578,0.11232381){\makebox(0,0)[lt]{\lineheight{1.25}\smash{\begin{tabular}[t]{l}1\end{tabular}}}}%
    \put(0,0){\includegraphics[width=\unitlength,page=13]{figextranspositions.pdf}}%
    \put(0.82328616,0.21141231){\makebox(0,0)[lt]{\lineheight{1.25}\smash{\begin{tabular}[t]{l}1\end{tabular}}}}%
    \put(0.84535772,0.18474255){\makebox(0,0)[lt]{\lineheight{1.25}\smash{\begin{tabular}[t]{l}2\end{tabular}}}}%
    \put(0.8012146,0.18474255){\makebox(0,0)[lt]{\lineheight{1.25}\smash{\begin{tabular}[t]{l}3\end{tabular}}}}%
    \put(0.87202737,0.23716246){\makebox(0,0)[lt]{\lineheight{1.25}\smash{\begin{tabular}[t]{l}4\end{tabular}}}}%
    \put(0.7708663,0.2380821){\makebox(0,0)[lt]{\lineheight{1.25}\smash{\begin{tabular}[t]{l}5\end{tabular}}}}%
    \put(0.75128726,0.20313551){\makebox(0,0)[lt]{\lineheight{1.25}\smash{\begin{tabular}[t]{l}6\end{tabular}}}}%
    \put(0.79635004,0.14922717){\makebox(0,0)[lt]{\lineheight{1.25}\smash{\begin{tabular}[t]{l}7\end{tabular}}}}%
    \put(0.85673012,0.14922717){\makebox(0,0)[lt]{\lineheight{1.25}\smash{\begin{tabular}[t]{l}8\end{tabular}}}}%
    \put(0.89229366,0.20312114){\makebox(0,0)[lt]{\lineheight{1.25}\smash{\begin{tabular}[t]{l}9\end{tabular}}}}%
    \put(0,0){\includegraphics[width=\unitlength,page=14]{figextranspositions.pdf}}%
    \put(0.71494101,0.18221708){\makebox(0,0)[lt]{\lineheight{1.25}\smash{\begin{tabular}[t]{l}8\end{tabular}}}}%
    \put(0,0){\includegraphics[width=\unitlength,page=15]{figextranspositions.pdf}}%
    \put(0.71494101,0.11232381){\makebox(0,0)[lt]{\lineheight{1.25}\smash{\begin{tabular}[t]{l}1\end{tabular}}}}%
    \put(0,0){\includegraphics[width=\unitlength,page=16]{figextranspositions.pdf}}%
    \put(0.92260793,0.06426871){\makebox(0,0)[lt]{\lineheight{1.25}\smash{\begin{tabular}[t]{l}1\end{tabular}}}}%
    \put(0.94467949,0.03759895){\makebox(0,0)[lt]{\lineheight{1.25}\smash{\begin{tabular}[t]{l}2\end{tabular}}}}%
    \put(0.90053637,0.03759895){\makebox(0,0)[lt]{\lineheight{1.25}\smash{\begin{tabular}[t]{l}3\end{tabular}}}}%
    \put(0.97134914,0.09001887){\makebox(0,0)[lt]{\lineheight{1.25}\smash{\begin{tabular}[t]{l}4\end{tabular}}}}%
    \put(0.87018806,0.0909385){\makebox(0,0)[lt]{\lineheight{1.25}\smash{\begin{tabular}[t]{l}5\end{tabular}}}}%
    \put(0.85060903,0.05599192){\makebox(0,0)[lt]{\lineheight{1.25}\smash{\begin{tabular}[t]{l}6\end{tabular}}}}%
    \put(0.89567181,0.00208357){\makebox(0,0)[lt]{\lineheight{1.25}\smash{\begin{tabular}[t]{l}7\end{tabular}}}}%
    \put(0.95605189,0.00208357){\makebox(0,0)[lt]{\lineheight{1.25}\smash{\begin{tabular}[t]{l}8\end{tabular}}}}%
    \put(0.99161543,0.05597755){\makebox(0,0)[lt]{\lineheight{1.25}\smash{\begin{tabular}[t]{l}9\end{tabular}}}}%
    \put(0,0){\includegraphics[width=\unitlength,page=17]{figextranspositions.pdf}}%
    \put(0.89151412,0.12703818){\makebox(0,0)[lt]{\lineheight{1.25}\smash{\begin{tabular}[t]{l}2\end{tabular}}}}%
  \end{picture}%
\endgroup%

%% file: staircase.pdf_tex
\begingroup%
  \makeatletter%
  \providecommand\color[2][]{%
    \errmessage{(Inkscape) Color is used for the text in Inkscape, but the package 'color.sty' is not loaded}%
    \renewcommand\color[2][]{}%
  }%
  \providecommand\transparent[1]{%
    \errmessage{(Inkscape) Transparency is used (non-zero) for the text in Inkscape, but the package 'transparent.sty' is not loaded}%
    \renewcommand\transparent[1]{}%
  }%
  \providecommand\rotatebox[2]{#2}%
  \newcommand*\fsize{\dimexpr\f@size pt\relax}%
  \newcommand*\lineheight[1]{\fontsize{\fsize}{#1\fsize}\selectfont}%
  \ifx\svgwidth\undefined%
    \setlength{\unitlength}{420.44769085bp}%
    \ifx\svgscale\undefined%
      \relax%
    \else%
      \setlength{\unitlength}{\unitlength * \real{\svgscale}}%
    \fi%
  \else%
    \setlength{\unitlength}{\svgwidth}%
  \fi%
  \global\let\svgwidth\undefined%
  \global\let\svgscale\undefined%
  \makeatother%
  \begin{picture}(1,1.30656841)%
    \lineheight{1}%
    \setlength\tabcolsep{0pt}%
    \put(0,0){\includegraphics[width=\unitlength,page=1]{staircase.pdf}}%
    \put(0.12539909,1.18194748){\makebox(0,0)[lt]{\lineheight{1.25}\smash{\begin{tabular}[t]{l}2\end{tabular}}}}%
    \put(0.1252133,1.11416258){\makebox(0,0)[lt]{\lineheight{1.25}\smash{\begin{tabular}[t]{l}1\end{tabular}}}}%
    \put(0.09329043,1.07135108){\makebox(0,0)[lt]{\lineheight{1.25}\smash{\begin{tabular}[t]{l}4\end{tabular}}}}%
    \put(0.02550553,1.03567482){\makebox(0,0)[lt]{\lineheight{1.25}\smash{\begin{tabular}[t]{l}9\end{tabular}}}}%
    \put(0.39297092,1.11773021){\makebox(0,0)[lt]{\lineheight{1.25}\smash{\begin{tabular}[t]{l}6\end{tabular}}}}%
    \put(0.32518633,1.08205397){\makebox(0,0)[lt]{\lineheight{1.25}\smash{\begin{tabular}[t]{l}5\end{tabular}}}}%
    \put(0.29166525,1.03567482){\makebox(0,0)[lt]{\lineheight{1.25}\smash{\begin{tabular}[t]{l}1\end{tabular}}}}%
    \put(0.29193813,0.96788994){\makebox(0,0)[lt]{\lineheight{1.25}\smash{\begin{tabular}[t]{l}3\end{tabular}}}}%
    \put(0.32161844,1.21405611){\makebox(0,0)[lt]{\lineheight{1.25}\smash{\begin{tabular}[t]{l}8\end{tabular}}}}%
    \put(0.39297092,1.17837985){\makebox(0,0)[lt]{\lineheight{1.25}\smash{\begin{tabular}[t]{l}7\end{tabular}}}}%
    \put(0.12539909,0.89653738){\makebox(0,0)[lt]{\lineheight{1.25}\smash{\begin{tabular}[t]{l}3\end{tabular}}}}%
    \put(0.1252133,0.82875253){\makebox(0,0)[lt]{\lineheight{1.25}\smash{\begin{tabular}[t]{l}2\end{tabular}}}}%
    \put(0.09329043,0.78594108){\makebox(0,0)[lt]{\lineheight{1.25}\smash{\begin{tabular}[t]{l}8\end{tabular}}}}%
    \put(0.02550553,0.75026484){\makebox(0,0)[lt]{\lineheight{1.25}\smash{\begin{tabular}[t]{l}7\end{tabular}}}}%
    \put(0.39297092,0.83232016){\makebox(0,0)[lt]{\lineheight{1.25}\smash{\begin{tabular}[t]{l}4\end{tabular}}}}%
    \put(0.32518633,0.79664392){\makebox(0,0)[lt]{\lineheight{1.25}\smash{\begin{tabular}[t]{l}9\end{tabular}}}}%
    \put(0.29166525,0.75026484){\makebox(0,0)[lt]{\lineheight{1.25}\smash{\begin{tabular}[t]{l}2\end{tabular}}}}%
    \put(0.29193813,0.68247999){\makebox(0,0)[lt]{\lineheight{1.25}\smash{\begin{tabular}[t]{l}1\end{tabular}}}}%
    \put(0.32161844,0.92864599){\makebox(0,0)[lt]{\lineheight{1.25}\smash{\begin{tabular}[t]{l}6\end{tabular}}}}%
    \put(0.39297092,0.89296975){\makebox(0,0)[lt]{\lineheight{1.25}\smash{\begin{tabular}[t]{l}5\end{tabular}}}}%
    \put(0.12539909,0.61112751){\makebox(0,0)[lt]{\lineheight{1.25}\smash{\begin{tabular}[t]{l}1\end{tabular}}}}%
    \put(0.1252133,0.54334256){\makebox(0,0)[lt]{\lineheight{1.25}\smash{\begin{tabular}[t]{l}3\end{tabular}}}}%
    \put(0.09329043,0.50053111){\makebox(0,0)[lt]{\lineheight{1.25}\smash{\begin{tabular}[t]{l}6\end{tabular}}}}%
    \put(0.02550553,0.46485482){\makebox(0,0)[lt]{\lineheight{1.25}\smash{\begin{tabular}[t]{l}5\end{tabular}}}}%
    \put(0.39297092,0.54691024){\makebox(0,0)[lt]{\lineheight{1.25}\smash{\begin{tabular}[t]{l}8\end{tabular}}}}%
    \put(0.32518633,0.511234){\makebox(0,0)[lt]{\lineheight{1.25}\smash{\begin{tabular}[t]{l}7\end{tabular}}}}%
    \put(0.29166525,0.46485482){\makebox(0,0)[lt]{\lineheight{1.25}\smash{\begin{tabular}[t]{l}3\end{tabular}}}}%
    \put(0.29193813,0.39706981){\makebox(0,0)[lt]{\lineheight{1.25}\smash{\begin{tabular}[t]{l}2\end{tabular}}}}%
    \put(0.32161844,0.64323607){\makebox(0,0)[lt]{\lineheight{1.25}\smash{\begin{tabular}[t]{l}4\end{tabular}}}}%
    \put(0.39297092,0.60755983){\makebox(0,0)[lt]{\lineheight{1.25}\smash{\begin{tabular}[t]{l}9\end{tabular}}}}%
    \put(0.12539909,0.32571723){\makebox(0,0)[lt]{\lineheight{1.25}\smash{\begin{tabular}[t]{l}2\end{tabular}}}}%
    \put(0.1252133,0.25793233){\makebox(0,0)[lt]{\lineheight{1.25}\smash{\begin{tabular}[t]{l}1\end{tabular}}}}%
    \put(0.09329043,0.21512078){\makebox(0,0)[lt]{\lineheight{1.25}\smash{\begin{tabular}[t]{l}4\end{tabular}}}}%
    \put(0.02550553,0.17944449){\makebox(0,0)[lt]{\lineheight{1.25}\smash{\begin{tabular}[t]{l}9\end{tabular}}}}%
    \put(0.39297092,0.26150001){\makebox(0,0)[lt]{\lineheight{1.25}\smash{\begin{tabular}[t]{l}6\end{tabular}}}}%
    \put(0.32518633,0.22582372){\makebox(0,0)[lt]{\lineheight{1.25}\smash{\begin{tabular}[t]{l}5\end{tabular}}}}%
    \put(0.29166525,0.17944449){\makebox(0,0)[lt]{\lineheight{1.25}\smash{\begin{tabular}[t]{l}1\end{tabular}}}}%
    \put(0.32161844,0.35782594){\makebox(0,0)[lt]{\lineheight{1.25}\smash{\begin{tabular}[t]{l}8\end{tabular}}}}%
    \put(0.39297092,0.32214965){\makebox(0,0)[lt]{\lineheight{1.25}\smash{\begin{tabular}[t]{l}7\end{tabular}}}}%
    \put(0.02550556,0.68604757){\makebox(0,0)[lt]{\lineheight{1.25}\smash{\begin{tabular}[t]{l}6\end{tabular}}}}%
    \put(0.09685806,0.65037138){\makebox(0,0)[lt]{\lineheight{1.25}\smash{\begin{tabular}[t]{l}5\end{tabular}}}}%
    \put(0.02550556,0.97145752){\makebox(0,0)[lt]{\lineheight{1.25}\smash{\begin{tabular}[t]{l}8\end{tabular}}}}%
    \put(0.09685806,0.93578125){\makebox(0,0)[lt]{\lineheight{1.25}\smash{\begin{tabular}[t]{l}7\end{tabular}}}}%
    \put(0.02550556,0.40063749){\makebox(0,0)[lt]{\lineheight{1.25}\smash{\begin{tabular}[t]{l}4\end{tabular}}}}%
    \put(0.09685806,0.3649612){\makebox(0,0)[lt]{\lineheight{1.25}\smash{\begin{tabular}[t]{l}9\end{tabular}}}}%
    \put(0.02550556,0.11522716){\makebox(0,0)[lt]{\lineheight{1.25}\smash{\begin{tabular}[t]{l}8\end{tabular}}}}%
    \put(0.09685806,0.07955087){\makebox(0,0)[lt]{\lineheight{1.25}\smash{\begin{tabular}[t]{l}7\end{tabular}}}}%
    \put(0,0){\includegraphics[width=\unitlength,page=2]{staircase.pdf}}%
    \put(0.23599544,1.28540862){\makebox(0,0)[lt]{\lineheight{1.25}\smash{\begin{tabular}[t]{l}6\end{tabular}}}}%
    \put(0.1682109,1.24973237){\makebox(0,0)[lt]{\lineheight{1.25}\smash{\begin{tabular}[t]{l}5\end{tabular}}}}%
    \put(0,0){\includegraphics[width=\unitlength,page=3]{staircase.pdf}}%
    \put(0.29193813,0.11165948){\makebox(0,0)[lt]{\lineheight{1.25}\smash{\begin{tabular}[t]{l}3\end{tabular}}}}%
    \put(0.18248101,0.01176628){\color[rgb]{0,0,0}\makebox(0,0)[lt]{\lineheight{1.25}\smash{\begin{tabular}[t]{l}9\end{tabular}}}}%
    \put(0,0){\includegraphics[width=\unitlength,page=4]{staircase.pdf}}%
    \put(0.25026591,0.04744257){\makebox(0,0)[lt]{\lineheight{1.25}\smash{\begin{tabular}[t]{l}4\end{tabular}}}}%
    \put(0,0){\includegraphics[width=\unitlength,page=5]{staircase.pdf}}%
    \put(0.69621925,1.03924242){\makebox(0,0)[lt]{\lineheight{1.25}\smash{\begin{tabular}[t]{l}1\end{tabular}}}}%
    \put(0.69603345,0.97145752){\makebox(0,0)[lt]{\lineheight{1.25}\smash{\begin{tabular}[t]{l}3\end{tabular}}}}%
    \put(0.66411063,0.92864604){\makebox(0,0)[lt]{\lineheight{1.25}\smash{\begin{tabular}[t]{l}9\end{tabular}}}}%
    \put(0.59632563,0.8929698){\makebox(0,0)[lt]{\lineheight{1.25}\smash{\begin{tabular}[t]{l}8\end{tabular}}}}%
    \put(0.96379138,0.97502515){\makebox(0,0)[lt]{\lineheight{1.25}\smash{\begin{tabular}[t]{l}5\end{tabular}}}}%
    \put(0.89600679,0.93934891){\makebox(0,0)[lt]{\lineheight{1.25}\smash{\begin{tabular}[t]{l}4\end{tabular}}}}%
    \put(0.86248571,0.8929698){\makebox(0,0)[lt]{\lineheight{1.25}\smash{\begin{tabular}[t]{l}3\end{tabular}}}}%
    \put(0.86275864,0.82518495){\makebox(0,0)[lt]{\lineheight{1.25}\smash{\begin{tabular}[t]{l}2\end{tabular}}}}%
    \put(0.8924388,1.07135106){\makebox(0,0)[lt]{\lineheight{1.25}\smash{\begin{tabular}[t]{l}7\end{tabular}}}}%
    \put(0.96379138,1.03567479){\makebox(0,0)[lt]{\lineheight{1.25}\smash{\begin{tabular}[t]{l}6\end{tabular}}}}%
    \put(0.69621925,0.75383242){\makebox(0,0)[lt]{\lineheight{1.25}\smash{\begin{tabular}[t]{l}2\end{tabular}}}}%
    \put(0.69603345,0.68604757){\makebox(0,0)[lt]{\lineheight{1.25}\smash{\begin{tabular}[t]{l}1\end{tabular}}}}%
    \put(0.66411063,0.64323617){\makebox(0,0)[lt]{\lineheight{1.25}\smash{\begin{tabular}[t]{l}7\end{tabular}}}}%
    \put(0.59632563,0.60755993){\makebox(0,0)[lt]{\lineheight{1.25}\smash{\begin{tabular}[t]{l}6\end{tabular}}}}%
    \put(0.96379138,0.6896152){\makebox(0,0)[lt]{\lineheight{1.25}\smash{\begin{tabular}[t]{l}9\end{tabular}}}}%
    \put(0.89600679,0.65393896){\makebox(0,0)[lt]{\lineheight{1.25}\smash{\begin{tabular}[t]{l}8\end{tabular}}}}%
    \put(0.86248571,0.60755993){\makebox(0,0)[lt]{\lineheight{1.25}\smash{\begin{tabular}[t]{l}1\end{tabular}}}}%
    \put(0.86275864,0.53977498){\makebox(0,0)[lt]{\lineheight{1.25}\smash{\begin{tabular}[t]{l}3\end{tabular}}}}%
    \put(0.8924388,0.78594103){\makebox(0,0)[lt]{\lineheight{1.25}\smash{\begin{tabular}[t]{l}5\end{tabular}}}}%
    \put(0.96379138,0.75026479){\makebox(0,0)[lt]{\lineheight{1.25}\smash{\begin{tabular}[t]{l}4\end{tabular}}}}%
    \put(0.69621925,0.4684224){\makebox(0,0)[lt]{\lineheight{1.25}\smash{\begin{tabular}[t]{l}3\end{tabular}}}}%
    \put(0.69603345,0.40063749){\makebox(0,0)[lt]{\lineheight{1.25}\smash{\begin{tabular}[t]{l}2\end{tabular}}}}%
    \put(0.66411063,0.35782594){\makebox(0,0)[lt]{\lineheight{1.25}\smash{\begin{tabular}[t]{l}5\end{tabular}}}}%
    \put(0.59632563,0.32214965){\makebox(0,0)[lt]{\lineheight{1.25}\smash{\begin{tabular}[t]{l}4\end{tabular}}}}%
    \put(0.96379138,0.40420518){\makebox(0,0)[lt]{\lineheight{1.25}\smash{\begin{tabular}[t]{l}7\end{tabular}}}}%
    \put(0.89600679,0.36852888){\makebox(0,0)[lt]{\lineheight{1.25}\smash{\begin{tabular}[t]{l}6\end{tabular}}}}%
    \put(0.86248571,0.32214965){\makebox(0,0)[lt]{\lineheight{1.25}\smash{\begin{tabular}[t]{l}2\end{tabular}}}}%
    \put(0.86275864,0.25436465){\makebox(0,0)[lt]{\lineheight{1.25}\smash{\begin{tabular}[t]{l}1\end{tabular}}}}%
    \put(0.8924388,0.50053111){\makebox(0,0)[lt]{\lineheight{1.25}\smash{\begin{tabular}[t]{l}9\end{tabular}}}}%
    \put(0.96379138,0.46485482){\makebox(0,0)[lt]{\lineheight{1.25}\smash{\begin{tabular}[t]{l}8\end{tabular}}}}%
    \put(0.69621925,0.18301206){\makebox(0,0)[lt]{\lineheight{1.25}\smash{\begin{tabular}[t]{l}1\end{tabular}}}}%
    \put(0.69603345,0.11522716){\makebox(0,0)[lt]{\lineheight{1.25}\smash{\begin{tabular}[t]{l}3\end{tabular}}}}%
    \put(0.66411063,0.07241561){\makebox(0,0)[lt]{\lineheight{1.25}\smash{\begin{tabular}[t]{l}9\end{tabular}}}}%
    \put(0.8924388,0.15447113){\makebox(0,0)[lt]{\lineheight{1.25}\smash{\begin{tabular}[t]{l}5\end{tabular}}}}%
    \put(0.8246542,0.11879484){\makebox(0,0)[lt]{\lineheight{1.25}\smash{\begin{tabular}[t]{l}4\end{tabular}}}}%
    \put(0.8924388,0.21512078){\makebox(0,0)[lt]{\lineheight{1.25}\smash{\begin{tabular}[t]{l}7\end{tabular}}}}%
    \put(0.59632563,0.54334256){\makebox(0,0)[lt]{\lineheight{1.25}\smash{\begin{tabular}[t]{l}5\end{tabular}}}}%
    \put(0.66767821,0.50766652){\makebox(0,0)[lt]{\lineheight{1.25}\smash{\begin{tabular}[t]{l}4\end{tabular}}}}%
    \put(0.59632563,0.82875253){\makebox(0,0)[lt]{\lineheight{1.25}\smash{\begin{tabular}[t]{l}7\end{tabular}}}}%
    \put(0.66767821,0.79307629){\makebox(0,0)[lt]{\lineheight{1.25}\smash{\begin{tabular}[t]{l}6\end{tabular}}}}%
    \put(0.59632563,0.25793233){\makebox(0,0)[lt]{\lineheight{1.25}\smash{\begin{tabular}[t]{l}9\end{tabular}}}}%
    \put(0.66767821,0.22225604){\makebox(0,0)[lt]{\lineheight{1.25}\smash{\begin{tabular}[t]{l}8\end{tabular}}}}%
    \put(0,0){\includegraphics[width=\unitlength,page=6]{staircase.pdf}}%
    \put(0.73546322,1.17837989){\makebox(0,0)[lt]{\lineheight{1.25}\smash{\begin{tabular}[t]{l}5\end{tabular}}}}%
    \put(0.66767821,1.14270364){\makebox(0,0)[lt]{\lineheight{1.25}\smash{\begin{tabular}[t]{l}4\end{tabular}}}}%
    \put(0,0){\includegraphics[width=\unitlength,page=7]{staircase.pdf}}%
    \put(0.86248571,1.17837989){\makebox(0,0)[lt]{\lineheight{1.25}\smash{\begin{tabular}[t]{l}2\end{tabular}}}}%
    \put(0.86248571,1.10702737){\makebox(0,0)[lt]{\lineheight{1.25}\smash{\begin{tabular}[t]{l}1\end{tabular}}}}%
    \put(0.66054295,1.08205397){\makebox(0,0)[lt]{\lineheight{1.25}\smash{\begin{tabular}[t]{l}8\end{tabular}}}}%
    \put(0.79140605,0.07598319){\makebox(0,0)[lt]{\lineheight{1.25}\smash{\begin{tabular}[t]{l}3\end{tabular}}}}%
    \put(0.74973363,0.01176597){\makebox(0,0)[lt]{\lineheight{1.25}\smash{\begin{tabular}[t]{l}9\end{tabular}}}}%
    \put(0.66767821,0.00819829){\makebox(0,0)[lt]{\lineheight{1.25}\smash{\begin{tabular}[t]{l}7\end{tabular}}}}%
    \put(0,0){\includegraphics[width=\unitlength,page=8]{staircase.pdf}}%
    \put(0.8068158,1.28540866){\makebox(0,0)[lt]{\lineheight{1.25}\smash{\begin{tabular}[t]{l}6\end{tabular}}}}%
    \put(0.76757183,1.21762372){\makebox(0,0)[lt]{\lineheight{1.25}\smash{\begin{tabular}[t]{l}2\end{tabular}}}}%
    \put(0.89600679,1.22475905){\makebox(0,0)[lt]{\lineheight{1.25}\smash{\begin{tabular}[t]{l}6\end{tabular}}}}%
    \put(0.8924388,1.28540871){\makebox(0,0)[lt]{\lineheight{1.25}\smash{\begin{tabular}[t]{l}8\end{tabular}}}}%
    \put(0,0){\includegraphics[width=\unitlength,page=9]{staircase.pdf}}%
  \end{picture}%
\endgroup%

%% file: figlemintervals.pdf_tex
\begingroup%
  \makeatletter%
  \providecommand\color[2][]{%
    \errmessage{(Inkscape) Color is used for the text in Inkscape, but the package 'color.sty' is not loaded}%
    \renewcommand\color[2][]{}%
  }%
  \providecommand\transparent[1]{%
    \errmessage{(Inkscape) Transparency is used (non-zero) for the text in Inkscape, but the package 'transparent.sty' is not loaded}%
    \renewcommand\transparent[1]{}%
  }%
  \providecommand\rotatebox[2]{#2}%
  \newcommand*\fsize{\dimexpr\f@size pt\relax}%
  \newcommand*\lineheight[1]{\fontsize{\fsize}{#1\fsize}\selectfont}%
  \ifx\svgwidth\undefined%
    \setlength{\unitlength}{377.46329648bp}%
    \ifx\svgscale\undefined%
      \relax%
    \else%
      \setlength{\unitlength}{\unitlength * \real{\svgscale}}%
    \fi%
  \else%
    \setlength{\unitlength}{\svgwidth}%
  \fi%
  \global\let\svgwidth\undefined%
  \global\let\svgscale\undefined%
  \makeatother%
  \begin{picture}(1,1.79690424)%
    \lineheight{1}%
    \setlength\tabcolsep{0pt}%
    \put(0,0){\includegraphics[width=\unitlength,page=1]{figlemintervals.pdf}}%
    \put(0.90567166,1.63000047){\makebox(0,0)[lt]{\lineheight{1.25}\smash{\begin{tabular}[t]{l}3\end{tabular}}}}%
    \put(0.9387316,1.58628776){\makebox(0,0)[lt]{\lineheight{1.25}\smash{\begin{tabular}[t]{l}9\end{tabular}}}}%
    \put(0.85925314,1.7015311){\makebox(0,0)[lt]{\lineheight{1.25}\smash{\begin{tabular}[t]{l}4\end{tabular}}}}%
    \put(0,0){\includegraphics[width=\unitlength,page=2]{figlemintervals.pdf}}%
    \put(0.16382164,1.74921769){\makebox(0,0)[lt]{\lineheight{1.25}\smash{\begin{tabular}[t]{l}1\end{tabular}}}}%
    \put(0.16361475,1.6737135){\makebox(0,0)[lt]{\lineheight{1.25}\smash{\begin{tabular}[t]{l}3\end{tabular}}}}%
    \put(0.34902158,1.58628776){\makebox(0,0)[lt]{\lineheight{1.25}\smash{\begin{tabular}[t]{l}3\end{tabular}}}}%
    \put(0,0){\includegraphics[width=\unitlength,page=3]{figlemintervals.pdf}}%
    \put(0.30290823,1.66576569){\makebox(0,0)[lt]{\lineheight{1.25}\smash{\begin{tabular}[t]{l}9\end{tabular}}}}%
    \put(0.22342983,1.7015311){\makebox(0,0)[lt]{\lineheight{1.25}\smash{\begin{tabular}[t]{l}4\end{tabular}}}}%
    \put(0.34902158,1.50680981){\makebox(0,0)[lt]{\lineheight{1.25}\smash{\begin{tabular}[t]{l}2\end{tabular}}}}%
    \put(0,0){\includegraphics[width=\unitlength,page=4]{figlemintervals.pdf}}%
    \put(0.58386253,1.65888632){\makebox(0,0)[lt]{\lineheight{1.25}\smash{\begin{tabular}[t]{l}\LARGE$\varphi$\end{tabular}}}}%
    \put(-0.00090549,1.65735652){\makebox(0,0)[lt]{\lineheight{1.25}\smash{\begin{tabular}[t]{l}\huge$\frac{1}{F_1}$\end{tabular}}}}%
    \put(0,0){\includegraphics[width=\unitlength,page=5]{figlemintervals.pdf}}%
    \put(0.98484554,1.34785395){\makebox(0,0)[lt]{\lineheight{1.25}\smash{\begin{tabular}[t]{l}1\end{tabular}}}}%
    \put(0.98514966,1.27234979){\makebox(0,0)[lt]{\lineheight{1.25}\smash{\begin{tabular}[t]{l}3\end{tabular}}}}%
    \put(0,0){\includegraphics[width=\unitlength,page=6]{figlemintervals.pdf}}%
    \put(0.94270544,1.2008197){\makebox(0,0)[lt]{\lineheight{1.25}\smash{\begin{tabular}[t]{l}7\end{tabular}}}}%
    \put(0.87912292,1.23261085){\makebox(0,0)[lt]{\lineheight{1.25}\smash{\begin{tabular}[t]{l}3\end{tabular}}}}%
    \put(0.84733161,1.27632421){\makebox(0,0)[lt]{\lineheight{1.25}\smash{\begin{tabular}[t]{l}4\end{tabular}}}}%
    \put(0.87891597,1.39554072){\makebox(0,0)[lt]{\lineheight{1.25}\smash{\begin{tabular}[t]{l}1\end{tabular}}}}%
    \put(0.9387316,1.42733187){\makebox(0,0)[lt]{\lineheight{1.25}\smash{\begin{tabular}[t]{l}4\end{tabular}}}}%
    \put(0.84335765,1.34785395){\makebox(0,0)[lt]{\lineheight{1.25}\smash{\begin{tabular}[t]{l}7\end{tabular}}}}%
    \put(0,0){\includegraphics[width=\unitlength,page=7]{figlemintervals.pdf}}%
    \put(0.58386247,1.34097458){\makebox(0,0)[lt]{\lineheight{1.25}\smash{\begin{tabular}[t]{l}\LARGE$\varphi$\end{tabular}}}}%
    \put(0,0){\includegraphics[width=\unitlength,page=8]{figlemintervals.pdf}}%
    \put(0.16382164,1.35182802){\makebox(0,0)[lt]{\lineheight{1.25}\smash{\begin{tabular}[t]{l}1\end{tabular}}}}%
    \put(0.16361475,1.27632386){\makebox(0,0)[lt]{\lineheight{1.25}\smash{\begin{tabular}[t]{l}3\end{tabular}}}}%
    \put(0,0){\includegraphics[width=\unitlength,page=9]{figlemintervals.pdf}}%
    \put(0.20753435,1.42733187){\makebox(0,0)[lt]{\lineheight{1.25}\smash{\begin{tabular}[t]{l}7\end{tabular}}}}%
    \put(0.98484554,1.34785395){\makebox(0,0)[lt]{\lineheight{1.25}\smash{\begin{tabular}[t]{l}1\end{tabular}}}}%
    \put(0,0){\includegraphics[width=\unitlength,page=10]{figlemintervals.pdf}}%
    \put(0.26984833,1.23261085){\makebox(0,0)[lt]{\lineheight{1.25}\smash{\begin{tabular}[t]{l}3\end{tabular}}}}%
    \put(0.26954427,1.38759289){\makebox(0,0)[lt]{\lineheight{1.25}\smash{\begin{tabular}[t]{l}1\end{tabular}}}}%
    \put(0.30688208,1.28029759){\makebox(0,0)[lt]{\lineheight{1.25}\smash{\begin{tabular}[t]{l}7\end{tabular}}}}%
    \put(0.21150842,1.19684632){\makebox(0,0)[lt]{\lineheight{1.25}\smash{\begin{tabular}[t]{l}4\end{tabular}}}}%
    \put(0.30290823,1.34785395){\makebox(0,0)[lt]{\lineheight{1.25}\smash{\begin{tabular}[t]{l}4\end{tabular}}}}%
    \put(-0.00090551,1.31957537){\makebox(0,0)[lt]{\lineheight{1.25}\smash{\begin{tabular}[t]{l}\huge$\frac{F_3}{y_1F_1}$\end{tabular}}}}%
    \put(0,0){\includegraphics[width=\unitlength,page=11]{figlemintervals.pdf}}%
    \put(0.98484554,1.06968134){\makebox(0,0)[lt]{\lineheight{1.25}\smash{\begin{tabular}[t]{l}1\end{tabular}}}}%
    \put(0.98514966,0.99417713){\makebox(0,0)[lt]{\lineheight{1.25}\smash{\begin{tabular}[t]{l}3\end{tabular}}}}%
    \put(0.79964492,0.91469924){\makebox(0,0)[lt]{\lineheight{1.25}\smash{\begin{tabular}[t]{l}3\end{tabular}}}}%
    \put(0.79943797,0.83919531){\makebox(0,0)[lt]{\lineheight{1.25}\smash{\begin{tabular}[t]{l}2\end{tabular}}}}%
    \put(0,0){\includegraphics[width=\unitlength,page=12]{figlemintervals.pdf}}%
    \put(0.92283565,1.02994234){\makebox(0,0)[lt]{\lineheight{1.25}\smash{\begin{tabular}[t]{l}7\end{tabular}}}}%
    \put(0.84733138,0.9902034){\makebox(0,0)[lt]{\lineheight{1.25}\smash{\begin{tabular}[t]{l}6\end{tabular}}}}%
    \put(0.92680949,0.91469947){\makebox(0,0)[lt]{\lineheight{1.25}\smash{\begin{tabular}[t]{l}7\end{tabular}}}}%
    \put(0.85130522,0.87496047){\makebox(0,0)[lt]{\lineheight{1.25}\smash{\begin{tabular}[t]{l}6\end{tabular}}}}%
    \put(0,0){\includegraphics[width=\unitlength,page=13]{figlemintervals.pdf}}%
    \put(0.58386253,0.98332408){\makebox(0,0)[lt]{\lineheight{1.25}\smash{\begin{tabular}[t]{l}\LARGE$\varphi$\end{tabular}}}}%
    \put(0,0){\includegraphics[width=\unitlength,page=14]{figlemintervals.pdf}}%
    \put(0.28701229,1.02994234){\makebox(0,0)[lt]{\lineheight{1.25}\smash{\begin{tabular}[t]{l}7\end{tabular}}}}%
    \put(0.21548192,0.87496047){\makebox(0,0)[lt]{\lineheight{1.25}\smash{\begin{tabular}[t]{l}6\end{tabular}}}}%
    \put(0.26984833,0.95443818){\makebox(0,0)[lt]{\lineheight{1.25}\smash{\begin{tabular}[t]{l}3\end{tabular}}}}%
    \put(-0.00090551,0.96192482){\makebox(0,0)[lt]{\lineheight{1.25}\smash{\begin{tabular}[t]{l}\huge$\frac{F_2}{y_1y_3}$\end{tabular}}}}%
    \put(0,0){\includegraphics[width=\unitlength,page=15]{figlemintervals.pdf}}%
    \put(0.87912292,0.63652634){\makebox(0,0)[lt]{\lineheight{1.25}\smash{\begin{tabular}[t]{l}2\end{tabular}}}}%
    \put(0.87891597,0.5610223){\makebox(0,0)[lt]{\lineheight{1.25}\smash{\begin{tabular}[t]{l}1\end{tabular}}}}%
    \put(0.84335765,0.51333561){\makebox(0,0)[lt]{\lineheight{1.25}\smash{\begin{tabular}[t]{l}7\end{tabular}}}}%
    \put(0.85130545,0.71203004){\makebox(0,0)[lt]{\lineheight{1.25}\smash{\begin{tabular}[t]{l}6\end{tabular}}}}%
    \put(0,0){\includegraphics[width=\unitlength,page=16]{figlemintervals.pdf}}%
    \put(0.9387316,0.5928135){\makebox(0,0)[lt]{\lineheight{1.25}\smash{\begin{tabular}[t]{l}5\end{tabular}}}}%
    \put(0.85527906,0.40206698){\makebox(0,0)[lt]{\lineheight{1.25}\smash{\begin{tabular}[t]{l}6\end{tabular}}}}%
    \put(0.83540985,0.44577925){\makebox(0,0)[lt]{\lineheight{1.25}\smash{\begin{tabular}[t]{l}5\end{tabular}}}}%
    \put(0.94270544,0.44577914){\makebox(0,0)[lt]{\lineheight{1.25}\smash{\begin{tabular}[t]{l}7\end{tabular}}}}%
    \put(0.98484554,0.5133355){\makebox(0,0)[lt]{\lineheight{1.25}\smash{\begin{tabular}[t]{l}1\end{tabular}}}}%
    \put(0,0){\includegraphics[width=\unitlength,page=17]{figlemintervals.pdf}}%
    \put(0.24329959,0.55704868){\makebox(0,0)[lt]{\lineheight{1.25}\smash{\begin{tabular}[t]{l}1\end{tabular}}}}%
    \put(0.2430927,0.48154441){\makebox(0,0)[lt]{\lineheight{1.25}\smash{\begin{tabular}[t]{l}3\end{tabular}}}}%
    \put(0.20753435,0.43385749){\makebox(0,0)[lt]{\lineheight{1.25}\smash{\begin{tabular}[t]{l}6\end{tabular}}}}%
    \put(0.21150842,0.60076141){\makebox(0,0)[lt]{\lineheight{1.25}\smash{\begin{tabular}[t]{l}5\end{tabular}}}}%
    \put(0,0){\includegraphics[width=\unitlength,page=18]{figlemintervals.pdf}}%
    \put(0.30290823,0.6722915){\makebox(0,0)[lt]{\lineheight{1.25}\smash{\begin{tabular}[t]{l}5\end{tabular}}}}%
    \put(0.30688208,0.52525714){\makebox(0,0)[lt]{\lineheight{1.25}\smash{\begin{tabular}[t]{l}7\end{tabular}}}}%
    \put(0.34902221,0.5928135){\makebox(0,0)[lt]{\lineheight{1.25}\smash{\begin{tabular}[t]{l}1\end{tabular}}}}%
    \put(0.20753435,0.66831777){\makebox(0,0)[lt]{\lineheight{1.25}\smash{\begin{tabular}[t]{l}7\end{tabular}}}}%
    \put(0.21548215,0.71203004){\makebox(0,0)[lt]{\lineheight{1.25}\smash{\begin{tabular}[t]{l}6\end{tabular}}}}%
    \put(0,0){\includegraphics[width=\unitlength,page=19]{figlemintervals.pdf}}%
    \put(0.58386253,0.58593413){\makebox(0,0)[lt]{\lineheight{1.25}\smash{\begin{tabular}[t]{l}\LARGE$\varphi$\end{tabular}}}}%
    \put(0,0){\includegraphics[width=\unitlength,page=20]{figlemintervals.pdf}}%
    \put(0.94270601,0.20734513){\makebox(0,0)[lt]{\lineheight{1.25}\smash{\begin{tabular}[t]{l}8\end{tabular}}}}%
    \put(0.90536754,0.15568471){\makebox(0,0)[lt]{\lineheight{1.25}\smash{\begin{tabular}[t]{l}1\end{tabular}}}}%
    \put(0.90567166,0.08018044){\makebox(0,0)[lt]{\lineheight{1.25}\smash{\begin{tabular}[t]{l}3\end{tabular}}}}%
    \put(0.9387316,0.03646748){\makebox(0,0)[lt]{\lineheight{1.25}\smash{\begin{tabular}[t]{l}9\end{tabular}}}}%
    \put(0,0){\includegraphics[width=\unitlength,page=21]{figlemintervals.pdf}}%
    \put(0.84733161,0.12389363){\makebox(0,0)[lt]{\lineheight{1.25}\smash{\begin{tabular}[t]{l}4\end{tabular}}}}%
    \put(0.79943797,0.2033714){\makebox(0,0)[lt]{\lineheight{1.25}\smash{\begin{tabular}[t]{l}1\end{tabular}}}}%
    \put(0,0){\includegraphics[width=\unitlength,page=22]{figlemintervals.pdf}}%
    \put(0.92283565,0.31464049){\makebox(0,0)[lt]{\lineheight{1.25}\smash{\begin{tabular}[t]{l}9\end{tabular}}}}%
    \put(0.84733138,0.27490149){\makebox(0,0)[lt]{\lineheight{1.25}\smash{\begin{tabular}[t]{l}8\end{tabular}}}}%
    \put(0.9387316,0.27092765){\makebox(0,0)[lt]{\lineheight{1.25}\smash{\begin{tabular}[t]{l}4\end{tabular}}}}%
    \put(0,0){\includegraphics[width=\unitlength,page=23]{figlemintervals.pdf}}%
    \put(0.26954364,0.23516249){\makebox(0,0)[lt]{\lineheight{1.25}\smash{\begin{tabular}[t]{l}2\end{tabular}}}}%
    \put(0.26984758,0.15965867){\makebox(0,0)[lt]{\lineheight{1.25}\smash{\begin{tabular}[t]{l}1\end{tabular}}}}%
    \put(0.21150807,0.04441517){\makebox(0,0)[lt]{\lineheight{1.25}\smash{\begin{tabular}[t]{l}4\end{tabular}}}}%
    \put(0.29098602,0.00467639){\makebox(0,0)[lt]{\lineheight{1.25}\smash{\begin{tabular}[t]{l}9\end{tabular}}}}%
    \put(0,0){\includegraphics[width=\unitlength,page=24]{figlemintervals.pdf}}%
    \put(0.30688265,0.04838912){\makebox(0,0)[lt]{\lineheight{1.25}\smash{\begin{tabular}[t]{l}8\end{tabular}}}}%
    \put(0.16361475,0.1238934){\makebox(0,0)[lt]{\lineheight{1.25}\smash{\begin{tabular}[t]{l}1\end{tabular}}}}%
    \put(0.28701229,0.31464049){\makebox(0,0)[lt]{\lineheight{1.25}\smash{\begin{tabular}[t]{l}9\end{tabular}}}}%
    \put(0.21150807,0.19542349){\makebox(0,0)[lt]{\lineheight{1.25}\smash{\begin{tabular}[t]{l}8\end{tabular}}}}%
    \put(0.30290823,0.11197176){\makebox(0,0)[lt]{\lineheight{1.25}\smash{\begin{tabular}[t]{l}4\end{tabular}}}}%
    \put(0,0){\includegraphics[width=\unitlength,page=25]{figlemintervals.pdf}}%
    \put(0.58386253,0.18854412){\makebox(0,0)[lt]{\lineheight{1.25}\smash{\begin{tabular}[t]{l}\LARGE$\varphi$\end{tabular}}}}%
  \end{picture}%
\endgroup%

%% file: figminimalstateNEW.pdf_tex
\begingroup%
  \makeatletter%
  \providecommand\color[2][]{%
    \errmessage{(Inkscape) Color is used for the text in Inkscape, but the package 'color.sty' is not loaded}%
    \renewcommand\color[2][]{}%
  }%
  \providecommand\transparent[1]{%
    \errmessage{(Inkscape) Transparency is used (non-zero) for the text in Inkscape, but the package 'transparent.sty' is not loaded}%
    \renewcommand\transparent[1]{}%
  }%
  \providecommand\rotatebox[2]{#2}%
  \newcommand*\fsize{\dimexpr\f@size pt\relax}%
  \newcommand*\lineheight[1]{\fontsize{\fsize}{#1\fsize}\selectfont}%
  \ifx\svgwidth\undefined%
    \setlength{\unitlength}{743.86598171bp}%
    \ifx\svgscale\undefined%
      \relax%
    \else%
      \setlength{\unitlength}{\unitlength * \real{\svgscale}}%
    \fi%
  \else%
    \setlength{\unitlength}{\svgwidth}%
  \fi%
  \global\let\svgwidth\undefined%
  \global\let\svgscale\undefined%
  \makeatother%
  \begin{picture}(1,0.69620149)%
    \lineheight{1}%
    \setlength\tabcolsep{0pt}%
    \put(0,0){\includegraphics[width=\unitlength,page=1]{figminimalstateNEW.pdf}}%
    \put(0.12719131,0.6028037){\makebox(0,0)[lt]{\lineheight{1.25}\smash{\begin{tabular}[t]{l}1\end{tabular}}}}%
    \put(0.14332325,0.63002634){\makebox(0,0)[lt]{\lineheight{1.25}\smash{\begin{tabular}[t]{l}3\end{tabular}}}}%
    \put(0.10904292,0.63002634){\makebox(0,0)[lt]{\lineheight{1.25}\smash{\begin{tabular}[t]{l}2\end{tabular}}}}%
    \put(0.06568837,0.6017955){\makebox(0,0)[lt]{\lineheight{1.25}\smash{\begin{tabular}[t]{l}6\end{tabular}}}}%
    \put(0.08887798,0.58566357){\makebox(0,0)[lt]{\lineheight{1.25}\smash{\begin{tabular}[t]{l}7\end{tabular}}}}%
    \put(0.16247989,0.58667184){\makebox(0,0)[lt]{\lineheight{1.25}\smash{\begin{tabular}[t]{l}8\end{tabular}}}}%
    \put(0.18566962,0.60381196){\makebox(0,0)[lt]{\lineheight{1.25}\smash{\begin{tabular}[t]{l}9\end{tabular}}}}%
    \put(0,0){\includegraphics[width=\unitlength,page=2]{figminimalstateNEW.pdf}}%
    \put(0.14332325,0.63002634){\makebox(0,0)[lt]{\lineheight{1.25}\smash{\begin{tabular}[t]{l}3\end{tabular}}}}%
    \put(0.10904292,0.63002634){\makebox(0,0)[lt]{\lineheight{1.25}\smash{\begin{tabular}[t]{l}2\end{tabular}}}}%
    \put(0.15239752,0.66531489){\makebox(0,0)[lt]{\lineheight{1.25}\smash{\begin{tabular}[t]{l}4\end{tabular}}}}%
    \put(0.09391926,0.66632315){\makebox(0,0)[lt]{\lineheight{1.25}\smash{\begin{tabular}[t]{l}5\end{tabular}}}}%
    \put(0,0){\includegraphics[width=\unitlength,page=3]{figminimalstateNEW.pdf}}%
    \put(0.12719132,0.44148432){\makebox(0,0)[lt]{\lineheight{1.25}\smash{\begin{tabular}[t]{l}1\end{tabular}}}}%
    \put(0.14332324,0.46870696){\makebox(0,0)[lt]{\lineheight{1.25}\smash{\begin{tabular}[t]{l}3\end{tabular}}}}%
    \put(0.10904293,0.46870696){\makebox(0,0)[lt]{\lineheight{1.25}\smash{\begin{tabular}[t]{l}2\end{tabular}}}}%
    \put(0.06568838,0.44047611){\makebox(0,0)[lt]{\lineheight{1.25}\smash{\begin{tabular}[t]{l}6\end{tabular}}}}%
    \put(0.08887799,0.42434419){\makebox(0,0)[lt]{\lineheight{1.25}\smash{\begin{tabular}[t]{l}7\end{tabular}}}}%
    \put(0.16247987,0.42535245){\makebox(0,0)[lt]{\lineheight{1.25}\smash{\begin{tabular}[t]{l}8\end{tabular}}}}%
    \put(0.18566959,0.44249258){\makebox(0,0)[lt]{\lineheight{1.25}\smash{\begin{tabular}[t]{l}9\end{tabular}}}}%
    \put(0,0){\includegraphics[width=\unitlength,page=4]{figminimalstateNEW.pdf}}%
    \put(0.14332324,0.46870696){\makebox(0,0)[lt]{\lineheight{1.25}\smash{\begin{tabular}[t]{l}3\end{tabular}}}}%
    \put(0.10904293,0.46870696){\makebox(0,0)[lt]{\lineheight{1.25}\smash{\begin{tabular}[t]{l}2\end{tabular}}}}%
    \put(0.15239751,0.50399551){\makebox(0,0)[lt]{\lineheight{1.25}\smash{\begin{tabular}[t]{l}4\end{tabular}}}}%
    \put(0.09391927,0.50500377){\makebox(0,0)[lt]{\lineheight{1.25}\smash{\begin{tabular}[t]{l}5\end{tabular}}}}%
    \put(0,0){\includegraphics[width=\unitlength,page=5]{figminimalstateNEW.pdf}}%
    \put(0.12719123,0.13901043){\makebox(0,0)[lt]{\lineheight{1.25}\smash{\begin{tabular}[t]{l}1\end{tabular}}}}%
    \put(0.14332316,0.16623307){\makebox(0,0)[lt]{\lineheight{1.25}\smash{\begin{tabular}[t]{l}3\end{tabular}}}}%
    \put(0.10904284,0.16623307){\makebox(0,0)[lt]{\lineheight{1.25}\smash{\begin{tabular}[t]{l}2\end{tabular}}}}%
    \put(0.06568829,0.13800223){\makebox(0,0)[lt]{\lineheight{1.25}\smash{\begin{tabular}[t]{l}6\end{tabular}}}}%
    \put(0.08887793,0.1218703){\makebox(0,0)[lt]{\lineheight{1.25}\smash{\begin{tabular}[t]{l}7\end{tabular}}}}%
    \put(0.16247981,0.12287856){\makebox(0,0)[lt]{\lineheight{1.25}\smash{\begin{tabular}[t]{l}8\end{tabular}}}}%
    \put(0.18566953,0.14001869){\makebox(0,0)[lt]{\lineheight{1.25}\smash{\begin{tabular}[t]{l}9\end{tabular}}}}%
    \put(0,0){\includegraphics[width=\unitlength,page=6]{figminimalstateNEW.pdf}}%
    \put(0.14332316,0.16623307){\makebox(0,0)[lt]{\lineheight{1.25}\smash{\begin{tabular}[t]{l}3\end{tabular}}}}%
    \put(0.10904284,0.16623307){\makebox(0,0)[lt]{\lineheight{1.25}\smash{\begin{tabular}[t]{l}2\end{tabular}}}}%
    \put(0.15239745,0.20152162){\makebox(0,0)[lt]{\lineheight{1.25}\smash{\begin{tabular}[t]{l}4\end{tabular}}}}%
    \put(0.09391915,0.20252988){\makebox(0,0)[lt]{\lineheight{1.25}\smash{\begin{tabular}[t]{l}5\end{tabular}}}}%
    \put(0,0){\includegraphics[width=\unitlength,page=7]{figminimalstateNEW.pdf}}%
    \put(0.12719132,0.28016499){\makebox(0,0)[lt]{\lineheight{1.25}\smash{\begin{tabular}[t]{l}1\end{tabular}}}}%
    \put(0.14332324,0.30738762){\makebox(0,0)[lt]{\lineheight{1.25}\smash{\begin{tabular}[t]{l}3\end{tabular}}}}%
    \put(0.10904293,0.30738762){\makebox(0,0)[lt]{\lineheight{1.25}\smash{\begin{tabular}[t]{l}2\end{tabular}}}}%
    \put(0.06568838,0.27915679){\makebox(0,0)[lt]{\lineheight{1.25}\smash{\begin{tabular}[t]{l}6\end{tabular}}}}%
    \put(0.08887799,0.26302486){\makebox(0,0)[lt]{\lineheight{1.25}\smash{\begin{tabular}[t]{l}7\end{tabular}}}}%
    \put(0.16247987,0.26403315){\makebox(0,0)[lt]{\lineheight{1.25}\smash{\begin{tabular}[t]{l}8\end{tabular}}}}%
    \put(0.18566959,0.28117325){\makebox(0,0)[lt]{\lineheight{1.25}\smash{\begin{tabular}[t]{l}9\end{tabular}}}}%
    \put(0,0){\includegraphics[width=\unitlength,page=8]{figminimalstateNEW.pdf}}%
    \put(0.14332324,0.30738762){\makebox(0,0)[lt]{\lineheight{1.25}\smash{\begin{tabular}[t]{l}3\end{tabular}}}}%
    \put(0.10904293,0.30738762){\makebox(0,0)[lt]{\lineheight{1.25}\smash{\begin{tabular}[t]{l}2\end{tabular}}}}%
    \put(0.15239751,0.34267614){\makebox(0,0)[lt]{\lineheight{1.25}\smash{\begin{tabular}[t]{l}4\end{tabular}}}}%
    \put(0.09391927,0.3436844){\makebox(0,0)[lt]{\lineheight{1.25}\smash{\begin{tabular}[t]{l}5\end{tabular}}}}%
    \put(0,0){\includegraphics[width=\unitlength,page=9]{figminimalstateNEW.pdf}}%
    \put(0.12719123,0.01802095){\makebox(0,0)[lt]{\lineheight{1.25}\smash{\begin{tabular}[t]{l}1\end{tabular}}}}%
    \put(0.14332316,0.04524358){\makebox(0,0)[lt]{\lineheight{1.25}\smash{\begin{tabular}[t]{l}3\end{tabular}}}}%
    \put(0.10904284,0.04524358){\makebox(0,0)[lt]{\lineheight{1.25}\smash{\begin{tabular}[t]{l}2\end{tabular}}}}%
    \put(0.06568829,0.01701275){\makebox(0,0)[lt]{\lineheight{1.25}\smash{\begin{tabular}[t]{l}6\end{tabular}}}}%
    \put(0.08887793,0.00088082){\makebox(0,0)[lt]{\lineheight{1.25}\smash{\begin{tabular}[t]{l}7\end{tabular}}}}%
    \put(0.16247981,0.00188911){\makebox(0,0)[lt]{\lineheight{1.25}\smash{\begin{tabular}[t]{l}8\end{tabular}}}}%
    \put(0.18566953,0.01902921){\makebox(0,0)[lt]{\lineheight{1.25}\smash{\begin{tabular}[t]{l}9\end{tabular}}}}%
    \put(0,0){\includegraphics[width=\unitlength,page=10]{figminimalstateNEW.pdf}}%
    \put(0.14332316,0.04524358){\makebox(0,0)[lt]{\lineheight{1.25}\smash{\begin{tabular}[t]{l}3\end{tabular}}}}%
    \put(0.10904284,0.04524358){\makebox(0,0)[lt]{\lineheight{1.25}\smash{\begin{tabular}[t]{l}2\end{tabular}}}}%
    \put(0.15239745,0.0805321){\makebox(0,0)[lt]{\lineheight{1.25}\smash{\begin{tabular}[t]{l}4\end{tabular}}}}%
    \put(0.09391915,0.08154036){\makebox(0,0)[lt]{\lineheight{1.25}\smash{\begin{tabular}[t]{l}5\end{tabular}}}}%
    \put(0,0){\includegraphics[width=\unitlength,page=11]{figminimalstateNEW.pdf}}%
    \put(0.3062594,0.67200364){\makebox(0,0)[lt]{\lineheight{1.25}\smash{\begin{tabular}[t]{l}3\end{tabular}}}}%
    \put(0.30615438,0.6336901){\makebox(0,0)[lt]{\lineheight{1.25}\smash{\begin{tabular}[t]{l}2\end{tabular}}}}%
    \put(0.28811086,0.609492){\makebox(0,0)[lt]{\lineheight{1.25}\smash{\begin{tabular}[t]{l}5\end{tabular}}}}%
    \put(0,0){\includegraphics[width=\unitlength,page=12]{figminimalstateNEW.pdf}}%
    \put(0.3365069,0.649822){\makebox(0,0)[lt]{\lineheight{1.25}\smash{\begin{tabular}[t]{l}9\end{tabular}}}}%
    \put(0.2941602,0.55303051){\makebox(0,0)[lt]{\lineheight{1.25}\smash{\begin{tabular}[t]{l}4\end{tabular}}}}%
    \put(0.28407787,0.57521163){\makebox(0,0)[lt]{\lineheight{1.25}\smash{\begin{tabular}[t]{l}9\end{tabular}}}}%
    \put(0.33852333,0.57521157){\makebox(0,0)[lt]{\lineheight{1.25}\smash{\begin{tabular}[t]{l}5\end{tabular}}}}%
    \put(0.35990668,0.609492){\makebox(0,0)[lt]{\lineheight{1.25}\smash{\begin{tabular}[t]{l}2\end{tabular}}}}%
    \put(0,0){\includegraphics[width=\unitlength,page=13]{figminimalstateNEW.pdf}}%
    \put(0.30625928,0.51068425){\makebox(0,0)[lt]{\lineheight{1.25}\smash{\begin{tabular}[t]{l}3\end{tabular}}}}%
    \put(0.3061543,0.47237071){\makebox(0,0)[lt]{\lineheight{1.25}\smash{\begin{tabular}[t]{l}2\end{tabular}}}}%
    \put(0.28811077,0.44817261){\makebox(0,0)[lt]{\lineheight{1.25}\smash{\begin{tabular}[t]{l}5\end{tabular}}}}%
    \put(0,0){\includegraphics[width=\unitlength,page=14]{figminimalstateNEW.pdf}}%
    \put(0.33650678,0.48850261){\makebox(0,0)[lt]{\lineheight{1.25}\smash{\begin{tabular}[t]{l}9\end{tabular}}}}%
    \put(0.28407779,0.41389224){\makebox(0,0)[lt]{\lineheight{1.25}\smash{\begin{tabular}[t]{l}9\end{tabular}}}}%
    \put(0.33852322,0.41389218){\makebox(0,0)[lt]{\lineheight{1.25}\smash{\begin{tabular}[t]{l}5\end{tabular}}}}%
    \put(0.35990662,0.44817261){\makebox(0,0)[lt]{\lineheight{1.25}\smash{\begin{tabular}[t]{l}2\end{tabular}}}}%
    \put(0,0){\includegraphics[width=\unitlength,page=15]{figminimalstateNEW.pdf}}%
    \put(0.30625934,0.34936486){\makebox(0,0)[lt]{\lineheight{1.25}\smash{\begin{tabular}[t]{l}3\end{tabular}}}}%
    \put(0.30615435,0.31105133){\makebox(0,0)[lt]{\lineheight{1.25}\smash{\begin{tabular}[t]{l}2\end{tabular}}}}%
    \put(0.28811083,0.28685316){\makebox(0,0)[lt]{\lineheight{1.25}\smash{\begin{tabular}[t]{l}5\end{tabular}}}}%
    \put(0,0){\includegraphics[width=\unitlength,page=16]{figminimalstateNEW.pdf}}%
    \put(0.30625931,0.20821044){\makebox(0,0)[lt]{\lineheight{1.25}\smash{\begin{tabular}[t]{l}3\end{tabular}}}}%
    \put(0.30615435,0.16989691){\makebox(0,0)[lt]{\lineheight{1.25}\smash{\begin{tabular}[t]{l}2\end{tabular}}}}%
    \put(0,0){\includegraphics[width=\unitlength,page=17]{figminimalstateNEW.pdf}}%
    \put(0.33650687,0.1860288){\makebox(0,0)[lt]{\lineheight{1.25}\smash{\begin{tabular}[t]{l}9\end{tabular}}}}%
    \put(0.35990674,0.1456988){\makebox(0,0)[lt]{\lineheight{1.25}\smash{\begin{tabular}[t]{l}2\end{tabular}}}}%
    \put(0,0){\includegraphics[width=\unitlength,page=18]{figminimalstateNEW.pdf}}%
    \put(0.30625934,0.06705581){\makebox(0,0)[lt]{\lineheight{1.25}\smash{\begin{tabular}[t]{l}3\end{tabular}}}}%
    \put(0.30615438,0.02874227){\makebox(0,0)[lt]{\lineheight{1.25}\smash{\begin{tabular}[t]{l}2\end{tabular}}}}%
    \put(0,0){\includegraphics[width=\unitlength,page=19]{figminimalstateNEW.pdf}}%
    \put(0.46240115,0.06060041){\makebox(0,0)[lt]{\lineheight{1.25}\smash{\begin{tabular}[t]{l}$\begin{array}{c}3 \\ 2\end{array}$\end{tabular}}}}%
    \put(0.46240118,0.18159007){\makebox(0,0)[lt]{\lineheight{1.25}\smash{\begin{tabular}[t]{l}$\begin{array}{c}3 \\ 2\ \ 9\end{array}$\end{tabular}}}}%
    \put(0.46240118,0.32274466){\makebox(0,0)[lt]{\lineheight{1.25}\smash{\begin{tabular}[t]{l}$\begin{array}{c}5\\ 3 \\ 2\end{array}$\end{tabular}}}}%
    \put(0.46240118,0.48406399){\makebox(0,0)[lt]{\lineheight{1.25}\smash{\begin{tabular}[t]{l}$\begin{array}{r} 3\hspace{6pt}  \\ 9\ \ 2\\5\end{array}$\end{tabular}}}}%
    \put(0.46240118,0.64538334){\makebox(0,0)[lt]{\lineheight{1.25}\smash{\begin{tabular}[t]{l}$\begin{array}{r} 3\hspace{6pt}   \\ 9\ \ 2\\5\\ 4\hspace{6pt} \end{array}$\end{tabular}}}}%
    \put(0,0){\includegraphics[width=\unitlength,page=20]{figminimalstateNEW.pdf}}%
  \end{picture}%
\endgroup%

%% file: bigon.pdf_tex
\begingroup%
  \makeatletter%
  \providecommand\color[2][]{%
    \errmessage{(Inkscape) Color is used for the text in Inkscape, but the package 'color.sty' is not loaded}%
    \renewcommand\color[2][]{}%
  }%
  \providecommand\transparent[1]{%
    \errmessage{(Inkscape) Transparency is used (non-zero) for the text in Inkscape, but the package 'transparent.sty' is not loaded}%
    \renewcommand\transparent[1]{}%
  }%
  \providecommand\rotatebox[2]{#2}%
  \newcommand*\fsize{\dimexpr\f@size pt\relax}%
  \newcommand*\lineheight[1]{\fontsize{\fsize}{#1\fsize}\selectfont}%
  \ifx\svgwidth\undefined%
    \setlength{\unitlength}{372.48495601bp}%
    \ifx\svgscale\undefined%
      \relax%
    \else%
      \setlength{\unitlength}{\unitlength * \real{\svgscale}}%
    \fi%
  \else%
    \setlength{\unitlength}{\svgwidth}%
  \fi%
  \global\let\svgwidth\undefined%
  \global\let\svgscale\undefined%
  \makeatother%
  \begin{picture}(1,0.17376394)%
    \lineheight{1}%
    \setlength\tabcolsep{0pt}%
    \put(0,0){\includegraphics[width=\unitlength,page=1]{bigon.pdf}}%
    \put(0.19696166,0.15417951){\makebox(0,0)[lt]{\lineheight{1.25}\smash{\begin{tabular}[t]{l}1\end{tabular}}}}%
    \put(0.1971714,0){\makebox(0,0)[lt]{\lineheight{1.25}\smash{\begin{tabular}[t]{l}2\end{tabular}}}}%
    \put(0.01567685,0.01960099){\makebox(0,0)[lt]{\lineheight{1.25}\smash{\begin{tabular}[t]{l}3\end{tabular}}}}%
    \put(0.01495585,0.12978273){\makebox(0,0)[lt]{\lineheight{1.25}\smash{\begin{tabular}[t]{l}4\end{tabular}}}}%
    \put(0.37968071,0.1301891){\makebox(0,0)[lt]{\lineheight{1.25}\smash{\begin{tabular}[t]{l}5\end{tabular}}}}%
    \put(0.38033615,0.01960099){\makebox(0,0)[lt]{\lineheight{1.25}\smash{\begin{tabular}[t]{l}6\end{tabular}}}}%
    \put(0.80101287,0.14186377){\makebox(0,0)[lt]{\lineheight{1.25}\smash{\begin{tabular}[t]{l}1\end{tabular}}}}%
    \put(0.80122247,0.01960099){\makebox(0,0)[lt]{\lineheight{1.25}\smash{\begin{tabular}[t]{l}2\end{tabular}}}}%
    \put(0.61972794,0.01960099){\makebox(0,0)[lt]{\lineheight{1.25}\smash{\begin{tabular}[t]{l}3\end{tabular}}}}%
    \put(0.61900689,0.14186377){\makebox(0,0)[lt]{\lineheight{1.25}\smash{\begin{tabular}[t]{l}4\end{tabular}}}}%
    \put(0.98373213,0.14227013){\makebox(0,0)[lt]{\lineheight{1.25}\smash{\begin{tabular}[t]{l}5\end{tabular}}}}%
    \put(0.98438751,0.01960099){\makebox(0,0)[lt]{\lineheight{1.25}\smash{\begin{tabular}[t]{l}6\end{tabular}}}}%
    \put(0,0){\includegraphics[width=\unitlength,page=2]{bigon.pdf}}%
  \end{picture}%
\endgroup%

%% file: bigonposet.pdf_tex
\begingroup%
  \makeatletter%
  \providecommand\color[2][]{%
    \errmessage{(Inkscape) Color is used for the text in Inkscape, but the package 'color.sty' is not loaded}%
    \renewcommand\color[2][]{}%
  }%
  \providecommand\transparent[1]{%
    \errmessage{(Inkscape) Transparency is used (non-zero) for the text in Inkscape, but the package 'transparent.sty' is not loaded}%
    \renewcommand\transparent[1]{}%
  }%
  \providecommand\rotatebox[2]{#2}%
  \newcommand*\fsize{\dimexpr\f@size pt\relax}%
  \newcommand*\lineheight[1]{\fontsize{\fsize}{#1\fsize}\selectfont}%
  \ifx\svgwidth\undefined%
    \setlength{\unitlength}{342.43490688bp}%
    \ifx\svgscale\undefined%
      \relax%
    \else%
      \setlength{\unitlength}{\unitlength * \real{\svgscale}}%
    \fi%
  \else%
    \setlength{\unitlength}{\svgwidth}%
  \fi%
  \global\let\svgwidth\undefined%
  \global\let\svgscale\undefined%
  \makeatother%
  \begin{picture}(1,1.57694247)%
    \lineheight{1}%
    \setlength\tabcolsep{0pt}%
    \put(0,0){\includegraphics[width=\unitlength,page=1]{bigonposet.pdf}}%
    \put(0.25455837,0.22778082){\makebox(0,0)[lt]{\lineheight{1.25}\smash{\begin{tabular}[t]{l}6\end{tabular}}}}%
    \put(0.25455837,0.4468005){\makebox(0,0)[lt]{\lineheight{1.25}\smash{\begin{tabular}[t]{l}5\end{tabular}}}}%
    \put(0.25455837,0.66143966){\makebox(0,0)[lt]{\lineheight{1.25}\smash{\begin{tabular}[t]{l}1\end{tabular}}}}%
    \put(0.25455837,0.88483975){\makebox(0,0)[lt]{\lineheight{1.25}\smash{\begin{tabular}[t]{l}4\end{tabular}}}}%
    \put(0.25455837,1.1038595){\makebox(0,0)[lt]{\lineheight{1.25}\smash{\begin{tabular}[t]{l}3\end{tabular}}}}%
    \put(0.25455837,1.31849883){\makebox(0,0)[lt]{\lineheight{1.25}\smash{\begin{tabular}[t]{l}2\end{tabular}}}}%
    \put(0.25455837,1.53751857){\makebox(0,0)[lt]{\lineheight{1.25}\smash{\begin{tabular}[t]{l}6\end{tabular}}}}%
    \put(0,0){\includegraphics[width=\unitlength,page=2]{bigonposet.pdf}}%
    \put(0.95542138,0.22778082){\makebox(0,0)[lt]{\lineheight{1.25}\smash{\begin{tabular}[t]{l}6\end{tabular}}}}%
    \put(0.95542138,0.4468005){\makebox(0,0)[lt]{\lineheight{1.25}\smash{\begin{tabular}[t]{l}5\end{tabular}}}}%
    \put(0.95542138,0.66582019){\makebox(0,0)[lt]{\lineheight{1.25}\smash{\begin{tabular}[t]{l}4\end{tabular}}}}%
    \put(0.95542138,0.88483981){\makebox(0,0)[lt]{\lineheight{1.25}\smash{\begin{tabular}[t]{l}3\end{tabular}}}}%
    \put(0.95542138,1.10385981){\makebox(0,0)[lt]{\lineheight{1.25}\smash{\begin{tabular}[t]{l}6\end{tabular}}}}%
  \end{picture}%
\endgroup%

%% file: figposetbij.pdf_tex
\begingroup%
  \makeatletter%
  \providecommand\color[2][]{%
    \errmessage{(Inkscape) Color is used for the text in Inkscape, but the package 'color.sty' is not loaded}%
    \renewcommand\color[2][]{}%
  }%
  \providecommand\transparent[1]{%
    \errmessage{(Inkscape) Transparency is used (non-zero) for the text in Inkscape, but the package 'transparent.sty' is not loaded}%
    \renewcommand\transparent[1]{}%
  }%
  \providecommand\rotatebox[2]{#2}%
  \newcommand*\fsize{\dimexpr\f@size pt\relax}%
  \newcommand*\lineheight[1]{\fontsize{\fsize}{#1\fsize}\selectfont}%
  \ifx\svgwidth\undefined%
    \setlength{\unitlength}{611.37034038bp}%
    \ifx\svgscale\undefined%
      \relax%
    \else%
      \setlength{\unitlength}{\unitlength * \real{\svgscale}}%
    \fi%
  \else%
    \setlength{\unitlength}{\svgwidth}%
  \fi%
  \global\let\svgwidth\undefined%
  \global\let\svgscale\undefined%
  \makeatother%
  \begin{picture}(1,0.7391094)%
    \lineheight{1}%
    \setlength\tabcolsep{0pt}%
    \put(0,0){\includegraphics[width=\unitlength,page=1]{figposetbij.pdf}}%
    \put(0.15898572,0.66108962){\makebox(0,0)[lt]{\lineheight{1.25}\smash{\begin{tabular}[t]{l}, $y_1\Big)$\end{tabular}}}}%
    \put(0,0){\includegraphics[width=\unitlength,page=2]{figposetbij.pdf}}%
    \put(-0.00123793,0.66116856){\makebox(0,0)[lt]{\lineheight{1.25}\smash{\begin{tabular}[t]{l}$C_y=\Big($\end{tabular}}}}%
    \put(-0.00123793,0.52867929){\makebox(0,0)[lt]{\lineheight{1.25}\smash{\begin{tabular}[t]{l}$B_y=\Big($\end{tabular}}}}%
    \put(-0.00123793,0.39619002){\makebox(0,0)[lt]{\lineheight{1.25}\smash{\begin{tabular}[t]{l}$A_y=\Big($\end{tabular}}}}%
    \put(-0.00123793,0.59492391){\makebox(0,0)[lt]{\lineheight{1.25}\smash{\begin{tabular}[t]{l}$C_1=\Big($\end{tabular}}}}%
    \put(-0.00123793,0.46243464){\makebox(0,0)[lt]{\lineheight{1.25}\smash{\begin{tabular}[t]{l}$B_1=\Big($\end{tabular}}}}%
    \put(-0.00123793,0.32994536){\makebox(0,0)[lt]{\lineheight{1.25}\smash{\begin{tabular}[t]{l}$A_1=\Big($\end{tabular}}}}%
    \put(0.15898572,0.52860035){\makebox(0,0)[lt]{\lineheight{1.25}\smash{\begin{tabular}[t]{l}, $y_1\Big)$\end{tabular}}}}%
    \put(0.15898572,0.39611106){\makebox(0,0)[lt]{\lineheight{1.25}\smash{\begin{tabular}[t]{l}, $y_1\Big)$\end{tabular}}}}%
    \put(0.15898572,0.32986644){\makebox(0,0)[lt]{\lineheight{1.25}\smash{\begin{tabular}[t]{l}, $1\Big)$\end{tabular}}}}%
    \put(0.15898572,0.46235572){\makebox(0,0)[lt]{\lineheight{1.25}\smash{\begin{tabular}[t]{l}, $1\Big)$\end{tabular}}}}%
    \put(0.15898572,0.59484499){\makebox(0,0)[lt]{\lineheight{1.25}\smash{\begin{tabular}[t]{l}, $1\Big)$\end{tabular}}}}%
    \put(0.01093703,0.72652244){\makebox(0,0)[lt]{\lineheight{1.25}\smash{\begin{tabular}[t]{l}$\underline{\cals(T(2))\times \{1,y_1\}}$\end{tabular}}}}%
    \put(0,0){\includegraphics[width=\unitlength,page=3]{figposetbij.pdf}}%
    \put(0.44284643,0.59492391){\makebox(0,0)[lt]{\lineheight{1.25}\smash{\begin{tabular}[t]{l}$\chi(C_1)=$\end{tabular}}}}%
    \put(0,0){\includegraphics[width=\unitlength,page=4]{figposetbij.pdf}}%
    \put(0.46728885,0.72652244){\makebox(0,0)[lt]{\lineheight{1.25}\smash{\begin{tabular}[t]{l}$\underline{\cals(T(3))}$\end{tabular}}}}%
    \put(0.44284643,0.66362203){\makebox(0,0)[lt]{\lineheight{1.25}\smash{\begin{tabular}[t]{l}$\chi(C_y)=$\end{tabular}}}}%
    \put(0.44284643,0.46243464){\makebox(0,0)[lt]{\lineheight{1.25}\smash{\begin{tabular}[t]{l}$\chi(B_1)=$\end{tabular}}}}%
    \put(0.78624526,0.72652244){\makebox(0,0)[lt]{\lineheight{1.25}\smash{\begin{tabular}[t]{l}$\underline{\cals(T(6))}$\end{tabular}}}}%
    \put(0.76425535,0.52867931){\makebox(0,0)[lt]{\lineheight{1.25}\smash{\begin{tabular}[t]{l}$\chi(B_y)=$\end{tabular}}}}%
    \put(0.76425535,0.39619002){\makebox(0,0)[lt]{\lineheight{1.25}\smash{\begin{tabular}[t]{l}$\chi(A_y)=$\end{tabular}}}}%
    \put(0.76425535,0.3299454){\makebox(0,0)[lt]{\lineheight{1.25}\smash{\begin{tabular}[t]{l}$\chi(A_1)=$\end{tabular}}}}%
    \put(0.31811436,0.50128262){\makebox(0,0)[lt]{\lineheight{1.25}\smash{\begin{tabular}[t]{l}\Huge$\mapsto$\end{tabular}}}}%
    \put(0.05160786,0.0032905){\makebox(0,0)[lt]{\lineheight{1.25}\smash{\begin{tabular}[t]{l}$A_1$\end{tabular}}}}%
    \put(0.05160786,0.07689571){\makebox(0,0)[lt]{\lineheight{1.25}\smash{\begin{tabular}[t]{l}$B_1$\end{tabular}}}}%
    \put(0.05160786,0.15050082){\makebox(0,0)[lt]{\lineheight{1.25}\smash{\begin{tabular}[t]{l}$C_1$\end{tabular}}}}%
    \put(0.12521299,0.17503585){\makebox(0,0)[lt]{\lineheight{1.25}\smash{\begin{tabular}[t]{l}$C_y$\end{tabular}}}}%
    \put(0.12521299,0.10143078){\makebox(0,0)[lt]{\lineheight{1.25}\smash{\begin{tabular}[t]{l}$B_y$\end{tabular}}}}%
    \put(0.12521299,0.02782557){\makebox(0,0)[lt]{\lineheight{1.25}\smash{\begin{tabular}[t]{l}$A_y$\end{tabular}}}}%
    \put(0,0){\includegraphics[width=\unitlength,page=5]{figposetbij.pdf}}%
    \put(0.0907457,0.17648884){\makebox(0,0)[lt]{\lineheight{1.25}\smash{\begin{tabular}[t]{l}$_1$\end{tabular}}}}%
    \put(0,0){\includegraphics[width=\unitlength,page=6]{figposetbij.pdf}}%
    \put(0.0907457,0.1028837){\makebox(0,0)[lt]{\lineheight{1.25}\smash{\begin{tabular}[t]{l}$_1$\end{tabular}}}}%
    \put(0,0){\includegraphics[width=\unitlength,page=7]{figposetbij.pdf}}%
    \put(0.0907457,0.02927855){\makebox(0,0)[lt]{\lineheight{1.25}\smash{\begin{tabular}[t]{l}$_1$\end{tabular}}}}%
    \put(0.14219574,0.13826805){\makebox(0,0)[lt]{\lineheight{1.25}\smash{\begin{tabular}[t]{l}$_5$\end{tabular}}}}%
    \put(0,0){\includegraphics[width=\unitlength,page=8]{figposetbij.pdf}}%
    \put(0.05141607,0.11373306){\makebox(0,0)[lt]{\lineheight{1.25}\smash{\begin{tabular}[t]{l}$_5$\end{tabular}}}}%
    \put(0,0){\includegraphics[width=\unitlength,page=9]{figposetbij.pdf}}%
    \put(0.05141607,0.04012788){\makebox(0,0)[lt]{\lineheight{1.25}\smash{\begin{tabular}[t]{l}$_4$\end{tabular}}}}%
    \put(0,0){\includegraphics[width=\unitlength,page=10]{figposetbij.pdf}}%
    \put(0.14219574,0.06466295){\makebox(0,0)[lt]{\lineheight{1.25}\smash{\begin{tabular}[t]{l}$_4$\end{tabular}}}}%
    \put(0.44416856,0.07689571){\makebox(0,0)[lt]{\lineheight{1.25}\smash{\begin{tabular}[t]{l}$\chi(B_1)$\end{tabular}}}}%
    \put(0.44416856,0.15050075){\makebox(0,0)[lt]{\lineheight{1.25}\smash{\begin{tabular}[t]{l}$\chi(C_1)$\end{tabular}}}}%
    \put(0,0){\includegraphics[width=\unitlength,page=11]{figposetbij.pdf}}%
    \put(0.46851184,0.11373306){\makebox(0,0)[lt]{\lineheight{1.25}\smash{\begin{tabular}[t]{l}$_5$\end{tabular}}}}%
    \put(0.54230867,0.17503575){\makebox(0,0)[lt]{\lineheight{1.25}\smash{\begin{tabular}[t]{l}$\chi(C_y)$\end{tabular}}}}%
    \put(0,0){\includegraphics[width=\unitlength,page=12]{figposetbij.pdf}}%
    \put(0.50784151,0.17648884){\makebox(0,0)[lt]{\lineheight{1.25}\smash{\begin{tabular}[t]{l}$_1$\end{tabular}}}}%
    \put(0.75085652,0.0032905){\makebox(0,0)[lt]{\lineheight{1.25}\smash{\begin{tabular}[t]{l}$\chi(A_1)$\end{tabular}}}}%
    \put(0.86126452,0.10143078){\makebox(0,0)[lt]{\lineheight{1.25}\smash{\begin{tabular}[t]{l}$\chi(B_y)$\end{tabular}}}}%
    \put(0.86126452,0.02782557){\makebox(0,0)[lt]{\lineheight{1.25}\smash{\begin{tabular}[t]{l}$\chi(A_y)$\end{tabular}}}}%
    \put(0,0){\includegraphics[width=\unitlength,page=13]{figposetbij.pdf}}%
    \put(0.82679729,0.02927855){\makebox(0,0)[lt]{\lineheight{1.25}\smash{\begin{tabular}[t]{l}$_1$\end{tabular}}}}%
    \put(0,0){\includegraphics[width=\unitlength,page=14]{figposetbij.pdf}}%
    \put(0.8782473,0.06466295){\makebox(0,0)[lt]{\lineheight{1.25}\smash{\begin{tabular}[t]{l}$_4$\end{tabular}}}}%
  \end{picture}%
\endgroup%

%% file: figgenbigon1.pdf_tex
\begingroup%
  \makeatletter%
  \providecommand\color[2][]{%
    \errmessage{(Inkscape) Color is used for the text in Inkscape, but the package 'color.sty' is not loaded}%
    \renewcommand\color[2][]{}%
  }%
  \providecommand\transparent[1]{%
    \errmessage{(Inkscape) Transparency is used (non-zero) for the text in Inkscape, but the package 'transparent.sty' is not loaded}%
    \renewcommand\transparent[1]{}%
  }%
  \providecommand\rotatebox[2]{#2}%
  \newcommand*\fsize{\dimexpr\f@size pt\relax}%
  \newcommand*\lineheight[1]{\fontsize{\fsize}{#1\fsize}\selectfont}%
  \ifx\svgwidth\undefined%
    \setlength{\unitlength}{298.21899029bp}%
    \ifx\svgscale\undefined%
      \relax%
    \else%
      \setlength{\unitlength}{\unitlength * \real{\svgscale}}%
    \fi%
  \else%
    \setlength{\unitlength}{\svgwidth}%
  \fi%
  \global\let\svgwidth\undefined%
  \global\let\svgscale\undefined%
  \makeatother%
  \begin{picture}(1,0.26567221)%
    \lineheight{1}%
    \setlength\tabcolsep{0pt}%
    \put(0,0){\includegraphics[width=\unitlength,page=1]{figgenbigon1.pdf}}%
    \put(0.5105079,0.15124837){\makebox(0,0)[lt]{\lineheight{1.25}\smash{\begin{tabular}[t]{l}$a$\end{tabular}}}}%
    \put(0.9387642,0.14406287){\makebox(0,0)[lt]{\lineheight{1.25}\smash{\begin{tabular}[t]{l}$b$\end{tabular}}}}%
    \put(0.76343779,0.01208748){\makebox(0,0)[lt]{\lineheight{1.25}\smash{\begin{tabular}[t]{l}$s$\end{tabular}}}}%
    \put(0.75050375,0.23986803){\makebox(0,0)[lt]{\lineheight{1.25}\smash{\begin{tabular}[t]{l}$s'$\end{tabular}}}}%
    \put(0,0){\includegraphics[width=\unitlength,page=2]{figgenbigon1.pdf}}%
    \put(0.21015341,0.01208748){\makebox(0,0)[lt]{\lineheight{1.25}\smash{\begin{tabular}[t]{l}$s$\end{tabular}}}}%
    \put(-0.00253785,0.15124837){\makebox(0,0)[lt]{\lineheight{1.25}\smash{\begin{tabular}[t]{l}$a$\end{tabular}}}}%
  \end{picture}%
\endgroup%

%% file: figgenbigon.pdf_tex
\begingroup%
  \makeatletter%
  \providecommand\color[2][]{%
    \errmessage{(Inkscape) Color is used for the text in Inkscape, but the package 'color.sty' is not loaded}%
    \renewcommand\color[2][]{}%
  }%
  \providecommand\transparent[1]{%
    \errmessage{(Inkscape) Transparency is used (non-zero) for the text in Inkscape, but the package 'transparent.sty' is not loaded}%
    \renewcommand\transparent[1]{}%
  }%
  \providecommand\rotatebox[2]{#2}%
  \newcommand*\fsize{\dimexpr\f@size pt\relax}%
  \newcommand*\lineheight[1]{\fontsize{\fsize}{#1\fsize}\selectfont}%
  \ifx\svgwidth\undefined%
    \setlength{\unitlength}{156.77378386bp}%
    \ifx\svgscale\undefined%
      \relax%
    \else%
      \setlength{\unitlength}{\unitlength * \real{\svgscale}}%
    \fi%
  \else%
    \setlength{\unitlength}{\svgwidth}%
  \fi%
  \global\let\svgwidth\undefined%
  \global\let\svgscale\undefined%
  \makeatother%
  \begin{picture}(1,0.50782938)%
    \lineheight{1}%
    \setlength\tabcolsep{0pt}%
    \put(0,0){\includegraphics[width=\unitlength,page=1]{figgenbigon.pdf}}%
    \put(0.06887533,0.2901695){\makebox(0,0)[lt]{\lineheight{1.25}\smash{\begin{tabular}[t]{l}$a$\end{tabular}}}}%
    \put(0.88351574,0.27650106){\makebox(0,0)[lt]{\lineheight{1.25}\smash{\begin{tabular}[t]{l}$b$\end{tabular}}}}%
    \put(0.55000527,0.02545418){\makebox(0,0)[lt]{\lineheight{1.25}\smash{\begin{tabular}[t]{l}$s$\end{tabular}}}}%
    \put(0.52540197,0.45874403){\makebox(0,0)[lt]{\lineheight{1.25}\smash{\begin{tabular}[t]{l}$s'$\end{tabular}}}}%
    \put(0,0){\includegraphics[width=\unitlength,page=2]{figgenbigon.pdf}}%
    \put(0.29874878,0.32785786){\makebox(0,0)[lt]{\lineheight{1.25}\smash{\begin{tabular}[t]{l}$y$\end{tabular}}}}%
    \put(0.21601709,0.2155017){\makebox(0,0)[lt]{\lineheight{1.25}\smash{\begin{tabular}[t]{l}$x_2$\end{tabular}}}}%
    \put(0.28733958,0.12923881){\makebox(0,0)[lt]{\lineheight{1.25}\smash{\begin{tabular}[t]{l}$x_1$\end{tabular}}}}%
    \put(0.66360041,0.26021796){\makebox(0,0)[lt]{\lineheight{1.25}\smash{\begin{tabular}[t]{l}$x_3$\end{tabular}}}}%
    \put(0,0){\includegraphics[width=\unitlength,page=3]{figgenbigon.pdf}}%
    \put(0.31490753,0.45874403){\makebox(0,0)[lt]{\lineheight{1.25}\smash{\begin{tabular}[t]{l}$s_2$\end{tabular}}}}%
    \put(0.1713886,0.45874403){\makebox(0,0)[lt]{\lineheight{1.25}\smash{\begin{tabular}[t]{l}$s_1$\end{tabular}}}}%
  \end{picture}%
\endgroup%

%% file: figadmissible.pdf_tex
\begingroup%
  \makeatletter%
  \providecommand\color[2][]{%
    \errmessage{(Inkscape) Color is used for the text in Inkscape, but the package 'color.sty' is not loaded}%
    \renewcommand\color[2][]{}%
  }%
  \providecommand\transparent[1]{%
    \errmessage{(Inkscape) Transparency is used (non-zero) for the text in Inkscape, but the package 'transparent.sty' is not loaded}%
    \renewcommand\transparent[1]{}%
  }%
  \providecommand\rotatebox[2]{#2}%
  \newcommand*\fsize{\dimexpr\f@size pt\relax}%
  \newcommand*\lineheight[1]{\fontsize{\fsize}{#1\fsize}\selectfont}%
  \ifx\svgwidth\undefined%
    \setlength{\unitlength}{381.77542516bp}%
    \ifx\svgscale\undefined%
      \relax%
    \else%
      \setlength{\unitlength}{\unitlength * \real{\svgscale}}%
    \fi%
  \else%
    \setlength{\unitlength}{\svgwidth}%
  \fi%
  \global\let\svgwidth\undefined%
  \global\let\svgscale\undefined%
  \makeatother%
  \begin{picture}(1,0.58280197)%
    \lineheight{1}%
    \setlength\tabcolsep{0pt}%
    \put(0,0){\includegraphics[width=\unitlength,page=1]{figadmissible.pdf}}%
    \put(0.02805595,0.4866703){\makebox(0,0)[lt]{\lineheight{1.25}\smash{\begin{tabular}[t]{l}$a$\end{tabular}}}}%
    \put(0.36258314,0.48105741){\makebox(0,0)[lt]{\lineheight{1.25}\smash{\begin{tabular}[t]{l}$b$\end{tabular}}}}%
    \put(0.27277725,0.38582438){\makebox(0,0)[lt]{\lineheight{1.25}\smash{\begin{tabular}[t]{l}$s$\end{tabular}}}}%
    \put(0.27446107,0.54017839){\makebox(0,0)[lt]{\lineheight{1.25}\smash{\begin{tabular}[t]{l}$s'$\end{tabular}}}}%
    \put(0,0){\includegraphics[width=\unitlength,page=2]{figadmissible.pdf}}%
    \put(0.05196643,0.46744918){\makebox(0,0)[lt]{\lineheight{1.25}\smash{\begin{tabular}[t]{l}$\zD_0$\end{tabular}}}}%
    \put(0,0){\includegraphics[width=\unitlength,page=3]{figadmissible.pdf}}%
    \put(0.57811748,0.4866703){\makebox(0,0)[lt]{\lineheight{1.25}\smash{\begin{tabular}[t]{l}$a$\end{tabular}}}}%
    \put(0.91264481,0.48105741){\makebox(0,0)[lt]{\lineheight{1.25}\smash{\begin{tabular}[t]{l}$b$\end{tabular}}}}%
    \put(0.8245228,0.54017839){\makebox(0,0)[lt]{\lineheight{1.25}\smash{\begin{tabular}[t]{l}$s'$\end{tabular}}}}%
    \put(0.54309278,0.46744918){\makebox(0,0)[lt]{\lineheight{1.25}\smash{\begin{tabular}[t]{l}$\zD_0$\end{tabular}}}}%
    \put(0,0){\includegraphics[width=\unitlength,page=4]{figadmissible.pdf}}%
    \put(0.64524717,0.46744918){\makebox(0,0)[lt]{\lineheight{1.25}\smash{\begin{tabular}[t]{l}$\zD_1$\end{tabular}}}}%
    \put(0,0){\includegraphics[width=\unitlength,page=5]{figadmissible.pdf}}%
    \put(0.76340262,0.46705137){\makebox(0,0)[lt]{\lineheight{1.25}\smash{\begin{tabular}[t]{l}$A$\end{tabular}}}}%
    \put(0.82283855,0.38582438){\makebox(0,0)[lt]{\lineheight{1.25}\smash{\begin{tabular}[t]{l}$s$\end{tabular}}}}%
    \put(0,0){\includegraphics[width=\unitlength,page=6]{figadmissible.pdf}}%
    \put(0.02805595,0.17234943){\makebox(0,0)[lt]{\lineheight{1.25}\smash{\begin{tabular}[t]{l}$a$\end{tabular}}}}%
    \put(0.36258314,0.16673656){\makebox(0,0)[lt]{\lineheight{1.25}\smash{\begin{tabular}[t]{l}$b$\end{tabular}}}}%
    \put(0.27277725,0.07150359){\makebox(0,0)[lt]{\lineheight{1.25}\smash{\begin{tabular}[t]{l}$s$\end{tabular}}}}%
    \put(0.27446107,0.22585748){\makebox(0,0)[lt]{\lineheight{1.25}\smash{\begin{tabular}[t]{l}$s'$\end{tabular}}}}%
    \put(0.16712334,0.08938353){\makebox(0,0)[lt]{\lineheight{1.25}\smash{\begin{tabular}[t]{l}$\zD_0$\end{tabular}}}}%
    \put(0,0){\includegraphics[width=\unitlength,page=7]{figadmissible.pdf}}%
    \put(0.57811748,0.17234943){\makebox(0,0)[lt]{\lineheight{1.25}\smash{\begin{tabular}[t]{l}$a$\end{tabular}}}}%
    \put(0.91264481,0.16673656){\makebox(0,0)[lt]{\lineheight{1.25}\smash{\begin{tabular}[t]{l}$b$\end{tabular}}}}%
    \put(0.73071383,0.0480356){\makebox(0,0)[lt]{\lineheight{1.25}\smash{\begin{tabular}[t]{l}$\zD_0$\end{tabular}}}}%
    \put(0.62560232,0.14919929){\makebox(0,0)[lt]{\lineheight{1.25}\smash{\begin{tabular}[t]{l}$\zD_1$\end{tabular}}}}%
    \put(0.59445519,0.21559465){\makebox(0,0)[lt]{\lineheight{1.25}\smash{\begin{tabular}[t]{l}$A_1$\end{tabular}}}}%
    \put(0,0){\includegraphics[width=\unitlength,page=8]{figadmissible.pdf}}%
    \put(0.83384018,0.14919929){\makebox(0,0)[lt]{\lineheight{1.25}\smash{\begin{tabular}[t]{l}$\zD_2$\end{tabular}}}}%
    \put(0.85376946,0.21559465){\makebox(0,0)[lt]{\lineheight{1.25}\smash{\begin{tabular}[t]{l}$A_2$\end{tabular}}}}%
    \put(0.73071383,0.2492592){\makebox(0,0)[lt]{\lineheight{1.25}\smash{\begin{tabular}[t]{l}$A_3$\end{tabular}}}}%
    \put(0.59950883,0.07773221){\makebox(0,0)[lt]{\lineheight{1.25}\smash{\begin{tabular}[t]{l}$C_1$\end{tabular}}}}%
    \put(0.85096518,0.07773221){\makebox(0,0)[lt]{\lineheight{1.25}\smash{\begin{tabular}[t]{l}$C_2$\end{tabular}}}}%
    \put(0.53738163,0.14919929){\makebox(0,0)[lt]{\lineheight{1.25}\smash{\begin{tabular}[t]{l}$D$\end{tabular}}}}%
    \put(0.73071383,0.14919929){\makebox(0,0)[lt]{\lineheight{1.25}\smash{\begin{tabular}[t]{l}$E$\end{tabular}}}}%
    \put(0.95344541,0.14919929){\makebox(0,0)[lt]{\lineheight{1.25}\smash{\begin{tabular}[t]{l}$F$\end{tabular}}}}%
    \put(0.73120015,0.00461709){\makebox(0,0)[lt]{\lineheight{1.25}\smash{\begin{tabular}[t]{l}$G$\end{tabular}}}}%
  \end{picture}%
\endgroup%

%% file: figstrings.pdf_tex
\begingroup%
  \makeatletter%
  \providecommand\color[2][]{%
    \errmessage{(Inkscape) Color is used for the text in Inkscape, but the package 'color.sty' is not loaded}%
    \renewcommand\color[2][]{}%
  }%
  \providecommand\transparent[1]{%
    \errmessage{(Inkscape) Transparency is used (non-zero) for the text in Inkscape, but the package 'transparent.sty' is not loaded}%
    \renewcommand\transparent[1]{}%
  }%
  \providecommand\rotatebox[2]{#2}%
  \newcommand*\fsize{\dimexpr\f@size pt\relax}%
  \newcommand*\lineheight[1]{\fontsize{\fsize}{#1\fsize}\selectfont}%
  \ifx\svgwidth\undefined%
    \setlength{\unitlength}{150.00000382bp}%
    \ifx\svgscale\undefined%
      \relax%
    \else%
      \setlength{\unitlength}{\unitlength * \real{\svgscale}}%
    \fi%
  \else%
    \setlength{\unitlength}{\svgwidth}%
  \fi%
  \global\let\svgwidth\undefined%
  \global\let\svgscale\undefined%
  \makeatother%
  \begin{picture}(1,0.10333333)%
    \lineheight{1}%
    \setlength\tabcolsep{0pt}%
    \put(0,0){\includegraphics[width=\unitlength,page=1]{figstrings.pdf}}%
  \end{picture}%
\endgroup%

%% file: figsumNEW.pdf_tex
\begingroup%
  \makeatletter%
  \providecommand\color[2][]{%
    \errmessage{(Inkscape) Color is used for the text in Inkscape, but the package 'color.sty' is not loaded}%
    \renewcommand\color[2][]{}%
  }%
  \providecommand\transparent[1]{%
    \errmessage{(Inkscape) Transparency is used (non-zero) for the text in Inkscape, but the package 'transparent.sty' is not loaded}%
    \renewcommand\transparent[1]{}%
  }%
  \providecommand\rotatebox[2]{#2}%
  \newcommand*\fsize{\dimexpr\f@size pt\relax}%
  \newcommand*\lineheight[1]{\fontsize{\fsize}{#1\fsize}\selectfont}%
  \ifx\svgwidth\undefined%
    \setlength{\unitlength}{288.75002884bp}%
    \ifx\svgscale\undefined%
      \relax%
    \else%
      \setlength{\unitlength}{\unitlength * \real{\svgscale}}%
    \fi%
  \else%
    \setlength{\unitlength}{\svgwidth}%
  \fi%
  \global\let\svgwidth\undefined%
  \global\let\svgscale\undefined%
  \makeatother%
  \begin{picture}(1,0.27686686)%
    \lineheight{1}%
    \setlength\tabcolsep{0pt}%
    \put(0,0){\includegraphics[width=\unitlength,page=1]{figsumNEW.pdf}}%
    \put(0.33766217,0.24129802){\makebox(0,0)[lt]{\lineheight{1.25}\smash{\begin{tabular}[t]{l}$B$\end{tabular}}}}%
    \put(0,0){\includegraphics[width=\unitlength,page=2]{figsumNEW.pdf}}%
    \put(0.07792195,0.24129802){\makebox(0,0)[lt]{\lineheight{1.25}\smash{\begin{tabular}[t]{l}$A$\end{tabular}}}}%
    \put(0.63636354,0.25021647){\makebox(0,0)[lt]{\lineheight{1.25}\smash{\begin{tabular}[t]{l}$A\oplus B$\end{tabular}}}}%
    \put(0,0){\includegraphics[width=\unitlength,page=3]{figsumNEW.pdf}}%
    \put(0.10238951,0.11543902){\makebox(0,0)[lt]{\lineheight{1.25}\smash{\begin{tabular}[t]{l}$p$\end{tabular}}}}%
    \put(0.59480551,0.08224771){\makebox(0,0)[lt]{\lineheight{1.25}\smash{\begin{tabular}[t]{l}$A\oplus[B,p]$\end{tabular}}}}%
    \put(0,0){\includegraphics[width=\unitlength,page=4]{figsumNEW.pdf}}%
    \put(0.33766217,0.08138527){\makebox(0,0)[lt]{\lineheight{1.25}\smash{\begin{tabular}[t]{l}$B$\end{tabular}}}}%
    \put(0.01298698,0.08138527){\makebox(0,0)[lt]{\lineheight{1.25}\smash{\begin{tabular}[t]{l}$A$\end{tabular}}}}%
    \put(0,0){\includegraphics[width=\unitlength,page=5]{figsumNEW.pdf}}%
  \end{picture}%
\endgroup%

%% file: figbdyex.pdf_tex
\begingroup%
  \makeatletter%
  \providecommand\color[2][]{%
    \errmessage{(Inkscape) Color is used for the text in Inkscape, but the package 'color.sty' is not loaded}%
    \renewcommand\color[2][]{}%
  }%
  \providecommand\transparent[1]{%
    \errmessage{(Inkscape) Transparency is used (non-zero) for the text in Inkscape, but the package 'transparent.sty' is not loaded}%
    \renewcommand\transparent[1]{}%
  }%
  \providecommand\rotatebox[2]{#2}%
  \newcommand*\fsize{\dimexpr\f@size pt\relax}%
  \newcommand*\lineheight[1]{\fontsize{\fsize}{#1\fsize}\selectfont}%
  \ifx\svgwidth\undefined%
    \setlength{\unitlength}{728.26197111bp}%
    \ifx\svgscale\undefined%
      \relax%
    \else%
      \setlength{\unitlength}{\unitlength * \real{\svgscale}}%
    \fi%
  \else%
    \setlength{\unitlength}{\svgwidth}%
  \fi%
  \global\let\svgwidth\undefined%
  \global\let\svgscale\undefined%
  \makeatother%
  \begin{picture}(1,0.17440389)%
    \lineheight{1}%
    \setlength\tabcolsep{0pt}%
    \put(0,0){\includegraphics[width=\unitlength,page=1]{figbdyex.pdf}}%
    \put(0.22037835,0.00242041){\makebox(0,0)[lt]{\lineheight{1.25}\smash{\begin{tabular}[t]{l}$A$\end{tabular}}}}%
    \put(0.32336325,0.07073241){\makebox(0,0)[lt]{\lineheight{1.25}\smash{\begin{tabular}[t]{l}$B$\end{tabular}}}}%
    \put(0.43252722,0.11192637){\makebox(0,0)[lt]{\lineheight{1.25}\smash{\begin{tabular}[t]{l}$C$\end{tabular}}}}%
    \put(0,0){\includegraphics[width=\unitlength,page=2]{figbdyex.pdf}}%
    \put(0.22140821,0.12840395){\makebox(0,0)[lt]{\lineheight{1.25}\smash{\begin{tabular}[t]{l}$a$\end{tabular}}}}%
    \put(0.32336325,0.15929944){\makebox(0,0)[lt]{\lineheight{1.25}\smash{\begin{tabular}[t]{l}$b$\end{tabular}}}}%
    \put(-0.00103923,0.00276235){\makebox(0,0)[lt]{\lineheight{1.25}\smash{\begin{tabular}[t]{l}$S=A\oplus [B \oplus[C,b],a]$\end{tabular}}}}%
    \put(0,0){\includegraphics[width=\unitlength,page=3]{figbdyex.pdf}}%
    \put(0.56743746,0.00276235){\makebox(0,0)[lt]{\lineheight{1.25}\smash{\begin{tabular}[t]{l}$\partial S=\partial A \oplus [\partial B\oplus[\partial C,b],a]$\end{tabular}}}}%
    \put(0,0){\includegraphics[width=\unitlength,page=4]{figbdyex.pdf}}%
  \end{picture}%
\endgroup%

%% file: figdimsym1NEW.pdf_tex
\begingroup%
  \makeatletter%
  \providecommand\color[2][]{%
    \errmessage{(Inkscape) Color is used for the text in Inkscape, but the package 'color.sty' is not loaded}%
    \renewcommand\color[2][]{}%
  }%
  \providecommand\transparent[1]{%
    \errmessage{(Inkscape) Transparency is used (non-zero) for the text in Inkscape, but the package 'transparent.sty' is not loaded}%
    \renewcommand\transparent[1]{}%
  }%
  \providecommand\rotatebox[2]{#2}%
  \newcommand*\fsize{\dimexpr\f@size pt\relax}%
  \newcommand*\lineheight[1]{\fontsize{\fsize}{#1\fsize}\selectfont}%
  \ifx\svgwidth\undefined%
    \setlength{\unitlength}{288.37779017bp}%
    \ifx\svgscale\undefined%
      \relax%
    \else%
      \setlength{\unitlength}{\unitlength * \real{\svgscale}}%
    \fi%
  \else%
    \setlength{\unitlength}{\svgwidth}%
  \fi%
  \global\let\svgwidth\undefined%
  \global\let\svgscale\undefined%
  \makeatother%
  \begin{picture}(1,0.46987715)%
    \lineheight{1}%
    \setlength\tabcolsep{0pt}%
    \put(0,0){\includegraphics[width=\unitlength,page=1]{figdimsym1NEW.pdf}}%
    \put(0.20241807,0.00611244){\makebox(0,0)[lt]{\lineheight{1.25}\smash{\begin{tabular}[t]{l}$i$\end{tabular}}}}%
    \put(0.201295,0.07150815){\makebox(0,0)[lt]{\lineheight{1.25}\smash{\begin{tabular}[t]{l}$i_{0}$\end{tabular}}}}%
    \put(0.20777843,0.30739025){\makebox(0,0)[lt]{\lineheight{1.25}\smash{\begin{tabular}[t]{l}$j$\end{tabular}}}}%
    \put(0.01134736,0.09000529){\makebox(0,0)[lt]{\lineheight{1.25}\smash{\begin{tabular}[t]{l}$a_0$\end{tabular}}}}%
    \put(0.10732789,0.13810084){\makebox(0,0)[lt]{\lineheight{1.25}\smash{\begin{tabular}[t]{l}$a_1$\end{tabular}}}}%
    \put(0.29302713,0.13666009){\makebox(0,0)[lt]{\lineheight{1.25}\smash{\begin{tabular}[t]{l}$b_1$\end{tabular}}}}%
    \put(0.37236503,0.0900903){\makebox(0,0)[lt]{\lineheight{1.25}\smash{\begin{tabular}[t]{l}$b_0$\end{tabular}}}}%
    \put(0.11565471,0.25599939){\makebox(0,0)[lt]{\lineheight{1.25}\smash{\begin{tabular}[t]{l}$c_1$\end{tabular}}}}%
    \put(0.28703132,0.25512019){\makebox(0,0)[lt]{\lineheight{1.25}\smash{\begin{tabular}[t]{l}$d_1$\end{tabular}}}}%
    \put(0.1181039,0.02222072){\makebox(0,0)[lt]{\lineheight{1.25}\smash{\begin{tabular}[t]{l}$s$\end{tabular}}}}%
    \put(0.28975379,0.02222072){\makebox(0,0)[lt]{\lineheight{1.25}\smash{\begin{tabular}[t]{l}$t$\end{tabular}}}}%
    \put(0,0){\includegraphics[width=\unitlength,page=2]{figdimsym1NEW.pdf}}%
    \put(0.1987158,0.17505691){\makebox(0,0)[lt]{\lineheight{1.25}\smash{\begin{tabular}[t]{l}$e_1$\end{tabular}}}}%
    \put(0.1987158,0.22707203){\makebox(0,0)[lt]{\lineheight{1.25}\smash{\begin{tabular}[t]{l}$e_2$\end{tabular}}}}%
    \put(0.10232683,0.33916221){\makebox(0,0)[lt]{\lineheight{1.25}\smash{\begin{tabular}[t]{l}$f_1$\end{tabular}}}}%
    \put(0.20115562,0.38597577){\makebox(0,0)[lt]{\lineheight{1.25}\smash{\begin{tabular}[t]{l}$h_1$\end{tabular}}}}%
    \put(0.32599184,0.33916221){\makebox(0,0)[lt]{\lineheight{1.25}\smash{\begin{tabular}[t]{l}$g_1$\end{tabular}}}}%
    \put(0.20115562,0.35476668){\makebox(0,0)[lt]{\lineheight{1.25}\smash{\begin{tabular}[t]{l}$h_2$\end{tabular}}}}%
    \put(0,0){\includegraphics[width=\unitlength,page=3]{figdimsym1NEW.pdf}}%
    \put(0.72256902,0.00611244){\makebox(0,0)[lt]{\lineheight{1.25}\smash{\begin{tabular}[t]{l}$i$\end{tabular}}}}%
    \put(0.72144591,0.07150815){\makebox(0,0)[lt]{\lineheight{1.25}\smash{\begin{tabular}[t]{l}$i_{0}$\end{tabular}}}}%
    \put(0.7279293,0.30739025){\makebox(0,0)[lt]{\lineheight{1.25}\smash{\begin{tabular}[t]{l}$j$\end{tabular}}}}%
    \put(0.53149841,0.09000529){\makebox(0,0)[lt]{\lineheight{1.25}\smash{\begin{tabular}[t]{l}$a_0$\end{tabular}}}}%
    \put(0.62747884,0.13810084){\makebox(0,0)[lt]{\lineheight{1.25}\smash{\begin{tabular}[t]{l}$a_1$\end{tabular}}}}%
    \put(0.81317801,0.13666009){\makebox(0,0)[lt]{\lineheight{1.25}\smash{\begin{tabular}[t]{l}$b_1$\end{tabular}}}}%
    \put(0.89251591,0.0900903){\makebox(0,0)[lt]{\lineheight{1.25}\smash{\begin{tabular}[t]{l}$b_0$\end{tabular}}}}%
    \put(0.63580566,0.25599939){\makebox(0,0)[lt]{\lineheight{1.25}\smash{\begin{tabular}[t]{l}$c_1$\end{tabular}}}}%
    \put(0.80718216,0.25512019){\makebox(0,0)[lt]{\lineheight{1.25}\smash{\begin{tabular}[t]{l}$d_1$\end{tabular}}}}%
    \put(0.63825489,0.02222072){\makebox(0,0)[lt]{\lineheight{1.25}\smash{\begin{tabular}[t]{l}$s$\end{tabular}}}}%
    \put(0.80990467,0.02222072){\makebox(0,0)[lt]{\lineheight{1.25}\smash{\begin{tabular}[t]{l}$t$\end{tabular}}}}%
    \put(0,0){\includegraphics[width=\unitlength,page=4]{figdimsym1NEW.pdf}}%
    \put(0.71886672,0.17505691){\makebox(0,0)[lt]{\lineheight{1.25}\smash{\begin{tabular}[t]{l}$e_1$\end{tabular}}}}%
    \put(0.71886672,0.22707203){\makebox(0,0)[lt]{\lineheight{1.25}\smash{\begin{tabular}[t]{l}$e_2$\end{tabular}}}}%
    \put(0.62247779,0.33916221){\makebox(0,0)[lt]{\lineheight{1.25}\smash{\begin{tabular}[t]{l}$f_1$\end{tabular}}}}%
    \put(0.72130657,0.38597577){\makebox(0,0)[lt]{\lineheight{1.25}\smash{\begin{tabular}[t]{l}$h_1$\end{tabular}}}}%
    \put(0.84614256,0.33916221){\makebox(0,0)[lt]{\lineheight{1.25}\smash{\begin{tabular}[t]{l}$g_1$\end{tabular}}}}%
    \put(0.72130657,0.35476668){\makebox(0,0)[lt]{\lineheight{1.25}\smash{\begin{tabular}[t]{l}$h_2$\end{tabular}}}}%
    \put(0,0){\includegraphics[width=\unitlength,page=5]{figdimsym1NEW.pdf}}%
    \put(0.20115562,0.44319238){\makebox(0,0)[lt]{\lineheight{1.25}\smash{\begin{tabular}[t]{l}$h_0$\end{tabular}}}}%
    \put(0.72130657,0.44319238){\makebox(0,0)[lt]{\lineheight{1.25}\smash{\begin{tabular}[t]{l}$h_0$\end{tabular}}}}%
  \end{picture}%
\endgroup%

%% file: figdimsym2NEW.pdf_tex
\begingroup%
  \makeatletter%
  \providecommand\color[2][]{%
    \errmessage{(Inkscape) Color is used for the text in Inkscape, but the package 'color.sty' is not loaded}%
    \renewcommand\color[2][]{}%
  }%
  \providecommand\transparent[1]{%
    \errmessage{(Inkscape) Transparency is used (non-zero) for the text in Inkscape, but the package 'transparent.sty' is not loaded}%
    \renewcommand\transparent[1]{}%
  }%
  \providecommand\rotatebox[2]{#2}%
  \newcommand*\fsize{\dimexpr\f@size pt\relax}%
  \newcommand*\lineheight[1]{\fontsize{\fsize}{#1\fsize}\selectfont}%
  \ifx\svgwidth\undefined%
    \setlength{\unitlength}{288.34863858bp}%
    \ifx\svgscale\undefined%
      \relax%
    \else%
      \setlength{\unitlength}{\unitlength * \real{\svgscale}}%
    \fi%
  \else%
    \setlength{\unitlength}{\svgwidth}%
  \fi%
  \global\let\svgwidth\undefined%
  \global\let\svgscale\undefined%
  \makeatother%
  \begin{picture}(1,0.43880857)%
    \lineheight{1}%
    \setlength\tabcolsep{0pt}%
    \put(0,0){\includegraphics[width=\unitlength,page=1]{figdimsym2NEW.pdf}}%
    \put(0.11609009,0.25461725){\makebox(0,0)[lt]{\lineheight{1.25}\smash{\begin{tabular}[t]{l}$c_1$\end{tabular}}}}%
    \put(0.28748406,0.25373808){\makebox(0,0)[lt]{\lineheight{1.25}\smash{\begin{tabular}[t]{l}$d_1$\end{tabular}}}}%
    \put(0,0){\includegraphics[width=\unitlength,page=2]{figdimsym2NEW.pdf}}%
    \put(0.20080602,0.17450433){\makebox(0,0)[lt]{\lineheight{1.25}\smash{\begin{tabular}[t]{l}$e_1$\end{tabular}}}}%
    \put(0,0){\includegraphics[width=\unitlength,page=3]{figdimsym2NEW.pdf}}%
    \put(0.20080602,0.22652471){\makebox(0,0)[lt]{\lineheight{1.25}\smash{\begin{tabular}[t]{l}$e_2$\end{tabular}}}}%
    \put(0,0){\includegraphics[width=\unitlength,page=4]{figdimsym2NEW.pdf}}%
    \put(0.20243866,0.00611309){\makebox(0,0)[lt]{\lineheight{1.25}\smash{\begin{tabular}[t]{l}$i$\end{tabular}}}}%
    \put(0.16490126,0.06111152){\makebox(0,0)[lt]{\lineheight{1.25}\smash{\begin{tabular}[t]{l}$i_{0}'$\end{tabular}}}}%
    \put(0.20779955,0.30742147){\makebox(0,0)[lt]{\lineheight{1.25}\smash{\begin{tabular}[t]{l}$j$\end{tabular}}}}%
    \put(0.11448059,0.1902203){\makebox(0,0)[lt]{\lineheight{1.25}\smash{\begin{tabular}[t]{l}$R$\end{tabular}}}}%
    \put(0.01134855,0.09001446){\makebox(0,0)[lt]{\lineheight{1.25}\smash{\begin{tabular}[t]{l}$a_0$\end{tabular}}}}%
    \put(0.1073388,0.13811483){\makebox(0,0)[lt]{\lineheight{1.25}\smash{\begin{tabular}[t]{l}$a_1$\end{tabular}}}}%
    \put(0.29305684,0.1366739){\makebox(0,0)[lt]{\lineheight{1.25}\smash{\begin{tabular}[t]{l}$b_1$\end{tabular}}}}%
    \put(0.37240283,0.09009944){\makebox(0,0)[lt]{\lineheight{1.25}\smash{\begin{tabular}[t]{l}$b_0$\end{tabular}}}}%
    \put(0.11811587,0.022223){\makebox(0,0)[lt]{\lineheight{1.25}\smash{\begin{tabular}[t]{l}$s$\end{tabular}}}}%
    \put(0.28978339,0.022223){\makebox(0,0)[lt]{\lineheight{1.25}\smash{\begin{tabular}[t]{l}$t$\end{tabular}}}}%
    \put(0,0){\includegraphics[width=\unitlength,page=5]{figdimsym2NEW.pdf}}%
    \put(0.20262451,0.07563934){\makebox(0,0)[lt]{\lineheight{1.25}\smash{\begin{tabular}[t]{l}$y$\end{tabular}}}}%
    \put(0.23252753,0.06111152){\makebox(0,0)[lt]{\lineheight{1.25}\smash{\begin{tabular}[t]{l}$i_{0}''$\end{tabular}}}}%
    \put(0.23475755,0.11417535){\makebox(0,0)[lt]{\lineheight{1.25}\smash{\begin{tabular}[t]{l}$\zb_1$\end{tabular}}}}%
    \put(0.15672693,0.1141762){\makebox(0,0)[lt]{\lineheight{1.25}\smash{\begin{tabular}[t]{l}$\za_1$\end{tabular}}}}%
    \put(0,0){\includegraphics[width=\unitlength,page=6]{figdimsym2NEW.pdf}}%
    \put(0.63629372,0.25461725){\makebox(0,0)[lt]{\lineheight{1.25}\smash{\begin{tabular}[t]{l}$c_1$\end{tabular}}}}%
    \put(0.80768769,0.25373808){\makebox(0,0)[lt]{\lineheight{1.25}\smash{\begin{tabular}[t]{l}$d_1$\end{tabular}}}}%
    \put(0,0){\includegraphics[width=\unitlength,page=7]{figdimsym2NEW.pdf}}%
    \put(0.72100956,0.17450433){\makebox(0,0)[lt]{\lineheight{1.25}\smash{\begin{tabular}[t]{l}$e_1$\end{tabular}}}}%
    \put(0,0){\includegraphics[width=\unitlength,page=8]{figdimsym2NEW.pdf}}%
    \put(0.72100956,0.22652471){\makebox(0,0)[lt]{\lineheight{1.25}\smash{\begin{tabular}[t]{l}$e_2$\end{tabular}}}}%
    \put(0,0){\includegraphics[width=\unitlength,page=9]{figdimsym2NEW.pdf}}%
    \put(0.72264219,0.00611309){\makebox(0,0)[lt]{\lineheight{1.25}\smash{\begin{tabular}[t]{l}$i$\end{tabular}}}}%
    \put(0.72800324,0.30742147){\makebox(0,0)[lt]{\lineheight{1.25}\smash{\begin{tabular}[t]{l}$j$\end{tabular}}}}%
    \put(0.53155241,0.09001446){\makebox(0,0)[lt]{\lineheight{1.25}\smash{\begin{tabular}[t]{l}$a_0$\end{tabular}}}}%
    \put(0.62754217,0.13811483){\makebox(0,0)[lt]{\lineheight{1.25}\smash{\begin{tabular}[t]{l}$a_1$\end{tabular}}}}%
    \put(0.81326019,0.1366739){\makebox(0,0)[lt]{\lineheight{1.25}\smash{\begin{tabular}[t]{l}$b_1$\end{tabular}}}}%
    \put(0.89260656,0.09009944){\makebox(0,0)[lt]{\lineheight{1.25}\smash{\begin{tabular}[t]{l}$b_0$\end{tabular}}}}%
    \put(0.63831954,0.022223){\makebox(0,0)[lt]{\lineheight{1.25}\smash{\begin{tabular}[t]{l}$s$\end{tabular}}}}%
    \put(0.80998689,0.022223){\makebox(0,0)[lt]{\lineheight{1.25}\smash{\begin{tabular}[t]{l}$t$\end{tabular}}}}%
    \put(0,0){\includegraphics[width=\unitlength,page=10]{figdimsym2NEW.pdf}}%
    \put(0.72282805,0.07563934){\makebox(0,0)[lt]{\lineheight{1.25}\smash{\begin{tabular}[t]{l}$y$\end{tabular}}}}%
    \put(0,0){\includegraphics[width=\unitlength,page=11]{figdimsym2NEW.pdf}}%
    \put(0.73794554,0.17864375){\makebox(0,0)[lt]{\lineheight{1.25}\smash{\begin{tabular}[t]{l}$\zg$\end{tabular}}}}%
    \put(0,0){\includegraphics[width=\unitlength,page=12]{figdimsym2NEW.pdf}}%
  \end{picture}%
\endgroup%

%% file: figdimsym3NEW.pdf_tex
\begingroup%
  \makeatletter%
  \providecommand\color[2][]{%
    \errmessage{(Inkscape) Color is used for the text in Inkscape, but the package 'color.sty' is not loaded}%
    \renewcommand\color[2][]{}%
  }%
  \providecommand\transparent[1]{%
    \errmessage{(Inkscape) Transparency is used (non-zero) for the text in Inkscape, but the package 'transparent.sty' is not loaded}%
    \renewcommand\transparent[1]{}%
  }%
  \providecommand\rotatebox[2]{#2}%
  \newcommand*\fsize{\dimexpr\f@size pt\relax}%
  \newcommand*\lineheight[1]{\fontsize{\fsize}{#1\fsize}\selectfont}%
  \ifx\svgwidth\undefined%
    \setlength{\unitlength}{292.40354523bp}%
    \ifx\svgscale\undefined%
      \relax%
    \else%
      \setlength{\unitlength}{\unitlength * \real{\svgscale}}%
    \fi%
  \else%
    \setlength{\unitlength}{\svgwidth}%
  \fi%
  \global\let\svgwidth\undefined%
  \global\let\svgscale\undefined%
  \makeatother%
  \begin{picture}(1,0.53714809)%
    \lineheight{1}%
    \setlength\tabcolsep{0pt}%
    \put(0,0){\includegraphics[width=\unitlength,page=1]{figdimsym3NEW.pdf}}%
    \put(0.2217664,0.00602836){\makebox(0,0)[lt]{\lineheight{1.25}\smash{\begin{tabular}[t]{l}$i$\end{tabular}}}}%
    \put(0.21878604,0.40667532){\makebox(0,0)[lt]{\lineheight{1.25}\smash{\begin{tabular}[t]{l}$j$\end{tabular}}}}%
    \put(0.02505936,0.19228302){\makebox(0,0)[lt]{\lineheight{1.25}\smash{\begin{tabular}[t]{l}$a_0$\end{tabular}}}}%
    \put(0.11971842,0.23971658){\makebox(0,0)[lt]{\lineheight{1.25}\smash{\begin{tabular}[t]{l}$a_1$\end{tabular}}}}%
    \put(0.3028608,0.23829563){\makebox(0,0)[lt]{\lineheight{1.25}\smash{\begin{tabular}[t]{l}$b_1$\end{tabular}}}}%
    \put(0.38110647,0.19236704){\makebox(0,0)[lt]{\lineheight{1.25}\smash{\begin{tabular}[t]{l}$b_0$\end{tabular}}}}%
    \put(0.12521598,0.1305617){\makebox(0,0)[lt]{\lineheight{1.25}\smash{\begin{tabular}[t]{l}$s$\end{tabular}}}}%
    \put(0.2893728,0.1305617){\makebox(0,0)[lt]{\lineheight{1.25}\smash{\begin{tabular}[t]{l}$t$\end{tabular}}}}%
    \put(0,0){\includegraphics[width=\unitlength,page=2]{figdimsym3NEW.pdf}}%
    \put(0.12782098,0.35602116){\makebox(0,0)[lt]{\lineheight{1.25}\smash{\begin{tabular}[t]{l}$c_1$\end{tabular}}}}%
    \put(0.29683815,0.35515419){\makebox(0,0)[lt]{\lineheight{1.25}\smash{\begin{tabular}[t]{l}$d_1$\end{tabular}}}}%
    \put(0,0){\includegraphics[width=\unitlength,page=3]{figdimsym3NEW.pdf}}%
    \put(0.21136204,0.3283182){\makebox(0,0)[lt]{\lineheight{1.25}\smash{\begin{tabular}[t]{l}$e_2$\end{tabular}}}}%
    \put(0,0){\includegraphics[width=\unitlength,page=4]{figdimsym3NEW.pdf}}%
    \put(0.21136204,0.27701925){\makebox(0,0)[lt]{\lineheight{1.25}\smash{\begin{tabular}[t]{l}$e_1$\end{tabular}}}}%
    \put(0,0){\includegraphics[width=\unitlength,page=5]{figdimsym3NEW.pdf}}%
    \put(0.12532515,0.28988221){\makebox(0,0)[lt]{\lineheight{1.25}\smash{\begin{tabular}[t]{l}$R$\end{tabular}}}}%
    \put(0,0){\includegraphics[width=\unitlength,page=6]{figdimsym3NEW.pdf}}%
    \put(0.73475594,0.00602836){\makebox(0,0)[lt]{\lineheight{1.25}\smash{\begin{tabular}[t]{l}$i$\end{tabular}}}}%
    \put(0.73177558,0.40667532){\makebox(0,0)[lt]{\lineheight{1.25}\smash{\begin{tabular}[t]{l}$j$\end{tabular}}}}%
    \put(0.53804904,0.19228302){\makebox(0,0)[lt]{\lineheight{1.25}\smash{\begin{tabular}[t]{l}$a_0$\end{tabular}}}}%
    \put(0.63270803,0.23971658){\makebox(0,0)[lt]{\lineheight{1.25}\smash{\begin{tabular}[t]{l}$a_1$\end{tabular}}}}%
    \put(0.81585031,0.23829563){\makebox(0,0)[lt]{\lineheight{1.25}\smash{\begin{tabular}[t]{l}$b_1$\end{tabular}}}}%
    \put(0.8940959,0.19236704){\makebox(0,0)[lt]{\lineheight{1.25}\smash{\begin{tabular}[t]{l}$b_0$\end{tabular}}}}%
    \put(0.6382056,0.1305617){\makebox(0,0)[lt]{\lineheight{1.25}\smash{\begin{tabular}[t]{l}$s$\end{tabular}}}}%
    \put(0.8023623,0.1305617){\makebox(0,0)[lt]{\lineheight{1.25}\smash{\begin{tabular}[t]{l}$t$\end{tabular}}}}%
    \put(0,0){\includegraphics[width=\unitlength,page=7]{figdimsym3NEW.pdf}}%
    \put(0.64081068,0.35602116){\makebox(0,0)[lt]{\lineheight{1.25}\smash{\begin{tabular}[t]{l}$c_1$\end{tabular}}}}%
    \put(0.80982785,0.35515419){\makebox(0,0)[lt]{\lineheight{1.25}\smash{\begin{tabular}[t]{l}$d_1$\end{tabular}}}}%
    \put(0,0){\includegraphics[width=\unitlength,page=8]{figdimsym3NEW.pdf}}%
    \put(0.72435174,0.3283182){\makebox(0,0)[lt]{\lineheight{1.25}\smash{\begin{tabular}[t]{l}$e_2$\end{tabular}}}}%
    \put(0,0){\includegraphics[width=\unitlength,page=9]{figdimsym3NEW.pdf}}%
    \put(0.72435174,0.27701925){\makebox(0,0)[lt]{\lineheight{1.25}\smash{\begin{tabular}[t]{l}$e_1$\end{tabular}}}}%
    \put(0,0){\includegraphics[width=\unitlength,page=10]{figdimsym3NEW.pdf}}%
    \put(0.63831476,0.28988221){\makebox(0,0)[lt]{\lineheight{1.25}\smash{\begin{tabular}[t]{l}$R$\end{tabular}}}}%
  \end{picture}%
\endgroup%